\documentclass{elsarticle}

\usepackage[english]{babel}

\usepackage[letterpaper,top=2cm,bottom=2cm,left=3cm,right=3cm,marginparwidth=1.75cm]{geometry}

\usepackage{amsmath}
\usepackage{graphicx}
\usepackage{subfigure}
\usepackage[colorlinks=true, allcolors=blue]{hyperref}
\usepackage{algorithm}
\usepackage{algpseudocode}

\usepackage{psfrag}
\usepackage{rotating}
\usepackage{color}
\usepackage{amssymb}
\usepackage{amsmath}
\usepackage{amsthm}
\usepackage{verbatim}
\usepackage{subfigure}
\usepackage{tabularx} 
\usepackage{tikz}
\usepackage{pgfplots}

\usepackage{soul}
\usepackage{latexsym}
\usepackage{amsfonts}
\usepackage{algorithm}
\usepackage{algorithmicx}
\usepackage{algpseudocode}
\usepackage{stmaryrd}
\usepgfplotslibrary{colorbrewer,external}
\newcommand{\f}[1]{\mathbf{#1}}
\newcommand{\R}{\mathbb R}

\newcommand{\Du}{\partial_{1}}
\newcommand{\Dv}{\partial_{2}}

\DeclareMathOperator{\supp}{supp}

\newcommand{\indOmega}{\mathcal{I}_{\Omega}}
\newcommand{\indSigma}{\mathcal{I}_{\Sigma}}
\newcommand{\indSigmaI}{\mathcal{I}_{\Sigma}^{\circ}}

\newcommand{\indChi}{\mathcal{I}_{\chi}}
\newcommand{\indChiI}{\mathcal{I}_{\chi}^{\circ}}

\newcommand{\US}[2]{\mathbb{S}_{#1}^{#2}}

\newcommand{\UW}{\mathbb{A}}
\newcommand{\PhiB}{\Phi}

\newcommand{\trunc}{\mathop{\mathrm{trunc}}}
\newcommand{\Trunc}{\mathop{\mathrm{Trunc}}}

\newcommand{\GG}{G} 

\newtheorem{theorem}{Theorem}[section]

\newtheorem{remark}{Remark}

\definecolor{pinegreen}{rgb}{0.0, 0.47, 0.44}
\definecolor{darkgreen}{rgb}{0.12, 0.35, 0.09}

\begin{document}

\begin{frontmatter}

\title{Adaptive isogeometric phase-field modeling of the Cahn-Hilliard equation: Suitably graded hierarchical refinement and coarsening on multi-patch geometries}

\author[unifi]{Cesare Bracco}
\author[unifi]{Carlotta Giannelli}
\author[unipv,imati]{Alessandro Reali\corref{cor1}}
\author[unipv]{Michele Torre}
\author[epfl,imati]{Rafael V\'azquez}
\cortext[cor1]{Corresponding author. Email: alessandro.reali@unipv.it}

\address[unifi]{Dipartimento di Matematica e Informatica ``U.Dini'', Università degli Studi di Firenze, Viale Morgagni 67/A, 50134 Firenze, Italy}
\address[unipv]{Dipartimento di Ingegneria Civile e Architettura, Universit\`a degli Studi di Pavia, Via Ferrata 3, 27100 Pavia, Italy}
\address[imati]{Istituto di Matematica Applicata e Tecnologie Informatiche ``E. Magenes'' del CNR, Via Ferrata 5, 27100 Pavia, Italy}
\address[epfl]{Institute of Mathematics, \'Ecole Polytechnique F\'ed\'erale de Lausanne, Station 8, 1015 Lausanne, Switzerland}

\begin{abstract}
We present an adaptive scheme for isogeometric phase-field modeling, to perform suitably graded hierarchical refinement and coarsening on both single- and multi-patch geometries by considering truncated hierarchical spline constructions which ensures $C^1$ continuity between patches. We apply the proposed algorithms to the Cahn-Hilliard equation, describing the time-evolving phase separation processes of immiscible fluids. We first verify the accuracy of the hierarchical spline scheme by comparing two classical indicators usually considered in phase-field modeling, for then demonstrating the effectiveness of the grading strategy in terms of accuracy per degree of freedom. A selection of numerical examples confirms the performance of the proposed scheme to simulate standard modes of phase separation using adaptive isogeometric analysis with smooth THB-spline constructions.
\end{abstract}


\end{frontmatter}

\section{Introduction}

Over the last fifty years, Tom Hughes has been the main actor of numerous groundbreaking advancements in almost every aspect of Computational Mechanics. This is the reason why it is not even clear whether primarily classifying him as a world expert of Solid Mechanics, Structural and Fluid Dynamics, Fluid-Structure Interaction (or, more in general, Multi-Physics) problems, Mechanics of Materials, Scientific Computing, or Applied Mathematics. His impact has been so significant and deep to make his peers consider him among the great scientific leaders of the Computational Mechanics community. Additionally, he has constantly been an example for several generations of researchers, always mentoring, motivating, and supporting young scholars. 

Isogeometric Analysis (IGA) \cite{hughes2005isogeometric, cottrell2009isogeometric} is only the latest of Tom Hughes' many seminal contributions. It constituted indeed a huge scientific success, going beyond the classical finite element method and actively engaging many researchers from different branches of Engineering, Numerical Analysis, and Computer-Aided Geometric Design (CAGD), while also attracting very significant interest from the industrial world.

The original idea of IGA is to bridge the gap between CAGD and numerical simulation by adopting the spline functions used to describe geometry in CAGD as shape functions to build Galerkin-based simulation methods within the isoparametric framework. Besides the potential advantages in geometry representation, the inherent smoothness of spline functions has been shown to offer various benefits including, among others, an improved accuracy-to-computational-cost ratio \cite{morganti2015patient}, enhanced spectral properties \cite{cottrell2006isogeometric,evans2013isogeometric}, high robustness \cite{lipton2010robustness}, also opening the door to the solution of higher order partial differential equations (PDE) in primal form \cite{gomez2008isogeometric, kiendl2009isogeometric}. 

In this context, the present paper deals with the simulation, via phase-field modeling, of the time-evolving phase separation processes of immiscible fluids described by the Cahn-Hillard equations \cite{CH1,CH2}. The problem is governed by a time-dependent PDE that is of fourth order in space, thus requiring at least $C^1$-continuity to be solved in primal form by means of a Galerkin method. In general, this is not easy to achieve with classical finite elements, while IGA has been proven to be an ideal framework for this kind of problems \cite{gomez2008isogeometric}. One of the issues arising with phase-field modeling is that, given its nature, it typically requires very fine meshes to properly resolve the high gradients present at the interfaces between the two different phases and, therefore, the possibility of using adaptive schemes able to locally refine and coarsen the adopted meshes would be highly desirable. 

We herein propose to achieve this by taking advantage of the local refinement properties of Truncated Hierarchical B-splines (THB-splines) \cite{giannelli2012} within the adaptive approach described in \cite{THB-refinement-coarsening}. To this end, we consider and compare the performance of two classical indicator choices for phase-field modeling problems \cite{hennig2018, nagaraja2019phase, proserpio2021phase}, and carefully study the effectiveness of our grading strategy in terms of accuracy per degree of freedom, showing the clear advantages that can be provided by the combination of THB-splines and suitably-graded adaptivity.

Moreover, since the dawn of IGA, it has been clear that dealing with complex geometries might involve the necessity of considering assemblies of multiple patches \cite{cottrell2007studies}, which in the context of high order PDEs might not be simple, in general, due to the need of enforcing $C^1$-continuity across patch boundaries. In this direction, many solutions have been proposed in the literature with different levels of complexity and accuracy (see, e.g., \cite{horger2019hybrid, KaSaTa19b, leonetti2020robust, farahat2023isogeometric, takacstoshniwal2023} and references therein). We follow here the recent extension to the (truncated) hierarchical spline framework \cite{BrGiKaVa23} of the multi-patch strategy presented in \cite{KaSaTa19b}, and we successfully test its potential to be a reliable tool for the adaptive phase-field modeling simulation of phase separation problems over non-trivial geometries comprising multiple patches.

The outline of the paper is as follows. {Section~\ref{sec:phase_field_discretization} introduces the phase-field model of the Cahn-Hilliard equation, while Section~\ref{sec:ada_mp} presents adaptive isogeometric methods with smooth (truncated) hierarchical splines on multi-patch geometries, after briefly recalling the definition of $C^1$ splines on analysis suitable geometries.} Suitably graded refinement and coarsening are also discussed. The adaptive isogeometric method for phase-field modeling of the Cahn-Hiliard equation is then presented in Section~\ref{sec:aigm-ch}, and a selection of numerical examples is reported in Section~\ref{sec:num_exe}. Conclusions are finally drawn in Section~\ref{sec:closure}.

\section{Phase-field modeling of the Cahn-Hilliard equation}
\label{sec:phase_field_discretization}

The Cahn-Hilliard equation \cite{CH1,CH2} is a model for the phase separation of immiscible fluids, assuming a simplified isotropic and isothermal framework. The  thermodynamic state of the mixture is described by a function $u\left(\mathbf{x}, t\right)$ of space $\mathbf{x}$ and time $t$, representing an order parameter of the mass fraction, typically referred to as the  ``phase field''. Under the isothermal hypothesis, the thermodynamic potential for the two-phase immiscible mixture is the Ginzburg-Landau free energy, defined as the following functional $\mathcal{G}:H^1(\Omega)\mapsto\R$ (with $\Omega$ being an open subset of $\mathbb{R}^d$, where $d$ is the spatial dimension):
\begin{equation*}
\mathcal{G}[u]=\int_\Omega\left(F(u) +  \frac{\lambda}{2}|\ensuremath{\nabla} u|^2\right)\hbox{\rm d} \mathbf{x}.
\end{equation*}
The gradient term in the equation above accounts for the interface free energy, while $F$ is the free energy of the homogeneous system, assumed to be
\begin{equation}\label{dw}
F(u)=\frac{\sigma}{4}\left(u^2-\frac{\nu}{\sigma}\right)^2,
\end{equation}
which is a non-convex function with a double-well structure characterized by two local minima (the so-called binodal points) for $u=\pm\sqrt{\nu/\sigma}$. 

The Cahn-Hilliard equation is then obtained imposing mass conservation, i.e., $\cfrac{\partial u}{\partial t} + \ensuremath{\nabla}\cdot \mathbf{j} = 0$,
assuming a mass flux $\mathbf{j}=-\ensuremath{\nabla}\left(\cfrac{\delta\mathcal{G}}{\delta u}\right)$ which implies a free energy decreasing in time, with $\cfrac{\delta\mathcal{G}}{\delta u}$ being the variational derivative of the free energy with respect to the variations of the phase field, $\delta u$. This yields:
\begin{equation*}
\frac{\partial u}{\partial t}=\ensuremath{\Delta}\left(F'(u)-\lambda\ensuremath{\Delta} u\right),
\end{equation*}
being 
\begin{equation*}
    F'(u) = \sigma u^3 - \nu u.
\end{equation*}

The minimization of the free energy then drives the dynamics of the solution of the Cahn-Hilliard equation above, with the phase field $u$ approaching the binodal points while creating regions of each pure phase which progressively expand over time. Such regions are separated by thin layers with a thickness that can be shown to have a characteristic length in the order of $\sqrt{\lambda}$.

The Cahn-Hilliard equation needs to be complemented by proper initial and boundary conditions. Assuming a smooth boundary of the domain decomposed into two complementary parts on which Dirichlet and Neumann boundary conditions are specified ($\partial \Omega=\overline{\partial \Omega_D\cup\partial \Omega_N}$), the initial/boundary-value problem over the spatial domain $\Omega$ and the time interval $[0,T]$ can be stated as follows: Given $u_0:\overline{\Omega}\mapsto\R$ and $u_D:\partial \Omega_D\mapsto\R$, find $u:\overline{\Omega}\times[0,T]\mapsto\R$ such that
\begin{eqnarray}
\label{ch_ibv}
 & \cfrac{\partial u}{\partial t}= \Delta \left( F'\left( u \right)  -\lambda\Delta u\right) & \qquad\hbox{in}\quad  \Omega\times[0,T], \\ \label{ch_dbc}   
 & u = u_D          & \qquad                             \hbox{on}\quad\partial \Omega_D\times[0,T], \\[.15cm]\label{ch_flux-condition}
 & \ensuremath{\nabla}\left(F'\left( u \right)-\lambda\Delta u\right)\cdot\mathbf{n}=0 & \qquad\hbox{on}\quad\partial \Omega_N\times[0,T], \\[.15cm]\label{ch_neumann-condition} 
 & \ensuremath{\nabla}\ u \cdot \mathbf{n}=0 & \qquad                 \hbox{on}\quad\partial \Omega\times[0,T],   \\[.15cm]
\label{ch_ibvf}
 & u(\mathbf{x},0) = u_0(\mathbf{x}) & \qquad   \hbox{in}  \quad  \overline{\Omega}.
\end{eqnarray}

In this setting, equations \eqref{ch_dbc} and \eqref{ch_neumann-condition} are Dirichlet (essential) boundary conditions requiring a proper definition of the functional spaces for trial and test functions if strongly imposed, while \eqref{ch_flux-condition} is a Neumann (natural) boundary condition. 

The set of equations \eqref{ch_ibv}-\eqref{ch_ibvf} is discretized numerically in space and time as described in the following by means of a Galerkin isogeometric approach and a generalized-$\alpha$ method, respectively. 

\subsection{Spatial discretization}
In this work, without loss of generality, we focus on problems where Neumann boundary conditions \eqref{ch_flux-condition} are enforced on the whole boundary $\partial \Omega$ of the domain, i.e., $\partial \Omega_N = \partial \Omega$ and $\partial \Omega_D = \emptyset$, and consequently condition \eqref{ch_dbc} is not enforced. Introducing the essential boundary condition \eqref{ch_neumann-condition}, a weak formulation of the initial/boundary-value problem reads as
\begin{equation}
\int_{\Omega} \delta u \dfrac{\partial u}{\partial t} \, \hbox{\rm d} \mathbf{x} =
- \int_{\Omega} \nabla \delta u \cdot \nabla F'(u) \, \hbox{\rm d} \mathbf{x}
- \int_{\Omega} \Delta \delta u \, \lambda \, \Delta u  \, \hbox{\rm d} \mathbf{x}
+ \int_{\partial \Omega} \nabla \delta u  \cdot \mathbf{n} \, \lambda \, \Delta u \, \hbox{\rm d} \mathbf{x},
\label{eq:weak_cont}
\end{equation}
where the integral defined on the boundary is the consistency term arising from the divergence theorem.

To discretize the weak form, $u$ is then approximated as
\begin{equation*}
    u\left( \mathbf{x}, t\right) = \mathbf{N}\left(\mathbf{x}\right) \hat{\mathbf{u}} \left( t \right),
\end{equation*}
being $\mathbf{N}\left(\mathbf{x}\right)$ the row vector of the adopted basis functions and $\hat{\mathbf{u}}\left( t \right)$ the column vector of the (time-dependent) control variables. Differently from classical second-order PDEs, the fourth-order Cahn-Hilliard equation requires $C^1$-continuous basis functions for approximation in a primal Galerkin fashion, which can be easily achieved using, e.g., B-splines on a single patch \cite{gomez2008isogeometric,kastner2016isogeometric}. However, in a multi-patch geometry, the continuity of standard basis at the patch interface is usually $C^0$. To overcome this, in the present work we select the functional spaces -- and consequently the basis functions -- to be $C^1$ at the interfaces, as described in the dedicated Section~\ref{sec:ada_mp}, and extend the discretization of the primal weak form to multi-patch domains.

Considering a Bubnov-Galerkin approach and thus exploiting the same functional space to approximate the field $u$ and its variation $\delta u$, the semi-discrete version of \eqref{eq:weak_cont} results in
\begin{equation*}
  \mathbf{M} \dot{ \hat{\mathbf{u}} } = - \overline{\mathbf{F}} - \mathbf{K}_\Delta \hat{\mathbf{u}} + \mathbf{K}_\partial \hat{\mathbf{u}},
\end{equation*}
where
\begin{equation*}
    \mathbf{M} = \int_{\Omega} \mathbf{N}^\top  \mathbf{N} \, \hbox{\rm d} \mathbf{x}
\end{equation*}
is the mass matrix and
\begin{equation*}
     \dot{ \hat{\mathbf{u}} } = \dfrac{\partial \hat{\mathbf{u}} }{\partial t}
\end{equation*}
is the velocity of the control variables. The gradient of $F'(u)$ is computed by means of the chain rule as
\begin{equation*}
    \nabla F'(u) = F''(u) \, \nabla \mathbf{N} \hat{\mathbf{u}},
\end{equation*}
leading to the (nonlinear) vector:
\begin{equation*}
    \overline{\mathbf{F}}(\hat{\mathbf{u}}) = \int_{\Omega}  F''(u) \, \nabla \mathbf{N}^\top \cdot \nabla \mathbf{N}  \hat{\mathbf{u}} \,  \hbox{\rm d} \mathbf{x},
\end{equation*}
with
\begin{equation}
    F''(u) = 3 \sigma \left( \mathbf{N} \hat{\mathbf{u}} \right)^2 - \nu.
    \label{eq:F''}
\end{equation}
Eventually, the matrices $ \mathbf{K}_\Delta$ and $\mathbf{K_\partial}$ are computed as
\begin{equation*}
    \mathbf{K}_\Delta = \int_{\Omega} \Delta \mathbf{N}^\top \, \lambda \, \Delta \mathbf{N} \, \hbox{\rm d} \mathbf{x},
\end{equation*}
and
\begin{equation*}
    \mathbf{K_\partial} = \int_{\Omega} \left(\nabla \mathbf{N} \cdot \mathbf{n} \right)^\top \, \lambda \, \Delta \mathbf{N} \, \hbox{\rm d} \mathbf{x}.
\end{equation*}

Essential boundary conditions \eqref{ch_neumann-condition} are weakly imposed by means of Nitsche's method \cite{zhao2017variational}, adding two vanishing terms to the left-hand side of \eqref{eq:weak_cont}:
\begin{equation*}
   - \int_{\partial \Omega}  \Delta \delta u  \, \lambda \,     \nabla  u  \cdot \mathbf{n} \, \hbox{\rm d} \mathbf{x}
\end{equation*}
and 
\begin{equation*}
   \int_{\partial \Omega}  \nabla \delta u  \cdot \mathbf{n} \, h \, \varepsilon_{N} \,     \nabla  u  \cdot \mathbf{n}  \, \hbox{\rm d} \mathbf{x},
\end{equation*}
where $\varepsilon_N$ is a positive constant and $h$ is a characteristic length, herein assumed to be equal to the element size.
Coherently, the augmented semi-discrete form leads to the following residual equation:
\begin{equation}
    \mathbf{R} \left( \hat{\mathbf{u}}, \dot{ \hat{\mathbf{u}} }  \right) =
    \mathbf{M} \dot{ \hat{\mathbf{u}} } + \overline{\mathbf{F}} (\hat{\mathbf{u}}) + \left( \mathbf{K}_\Delta  - \mathbf{K}_\partial - \mathbf{K}_\partial^\top + \mathbf{M}_N \right) \hat{\mathbf{u}} = \mathbf{0},
    \label{eq:residual_CH}
\end{equation}
with
\begin{eqnarray*}
  \mathbf{M}_N=\int_{\partial \Omega} \left( \nabla \mathbf{N}  \cdot \mathbf{n} \right)^\top \, h \,  \varepsilon_{N} \,   \left( \nabla \mathbf{N}  \cdot \mathbf{n} \right) \, \hbox{\rm d} \mathbf{x}. 
\end{eqnarray*}
Such a residual equation is advanced in time by means of the generalized-$\alpha$ method, as presented in what follows.

\subsection{Time marching scheme}
The generalized-$\alpha$ method was successfully employed in previous studies \cite{gomez2008isogeometric, kastner2016isogeometric} to time integrate the Cahn-Hilliard equation and is therefore adopted also herein. Given the control variables and their velocities at the previous time step $t_n$, denoted as $\hat{\mathbf{u}}_n$ and $\dot{ \hat{\mathbf{u}} }_n$, the solution at the current time step $t_{n+1} = t_n + \Delta t$ is computed by solving
\begin{equation}
    \mathbf{R} \left( \hat{\mathbf{u}}_{n+\alpha_{f}}, \dot{ \hat{\mathbf{u}}}_{n+\alpha_{m}}   \right) = \mathbf{0}.
    \label{eq:res_alpha}
\end{equation}
The vectors $\hat{\mathbf{u}}_{n+\alpha_{f}}$ and $\dot{ \hat{\mathbf{u}}}_{n+\alpha_{m}}$ represent the solution and its velocity at intermediate $\alpha$-levels between $t_n$ and $t_{n+1}$, and are respectively computed as
\begin{equation*}
    \hat{\mathbf{u}}_{n+\alpha_{f}} = \hat{\mathbf{u}}_n + \alpha_f \left( \hat{\mathbf{u}}_{n+1} - \hat{\mathbf{u}}_n \right)
\end{equation*}
and
\begin{equation*}
    \dot{\hat{\mathbf{u}}}_{n+\alpha_{f}} = \dot{\hat{\mathbf{u}}}_n + \alpha_m \left( \dot{\hat{\mathbf{u}}}_{n+1} - \dot{\hat{\mathbf{u}}}_n \right),
\end{equation*}
being the solution at $t_{n+1}$ related to the current velocity and to the solution at the previous time step by the following condition:
\begin{equation*}
    \hat{\mathbf{u}}_{n+1} = \hat{\mathbf{u}}_{n} + \Delta t \,  \dot{\hat{\mathbf{u}}}_{n} + \gamma \Delta t \left( \dot{\hat{\mathbf{u}}}_{n+1}   - \dot{\hat{\mathbf{u}}}_{n} \right).
\end{equation*}

The choice of the values of the parameters governs the accuracy, the numerical damping, and the stability of the method. We follow the classical strategy by selecting
\begin{equation*}
    \gamma = \dfrac{1}{2} + \alpha_m - \alpha_f, \hspace{10pt} \alpha_m = \dfrac{1}{2} \left( \dfrac{3-\rho_\infty}{1+\rho_\infty} \right), \hspace{10pt} \alpha_f =  \dfrac{1}{1+\rho_\infty} ,
\end{equation*}
and assuming that the spectral radius of the amplification matrix is $\rho_\infty = 0.5$.

The iterative solution process to solve \eqref{eq:res_alpha} starts from the predictions
\begin{equation*}
    \hat{\mathbf{u}}_{n+1} = \hat{\mathbf{u}}_{n},
\end{equation*}
\begin{equation*}
    \dot{\hat{\mathbf{u}}}_{n+1} = \dfrac{\gamma - 1 }{\gamma} \dot{\hat{\mathbf{u}}}_{n},
\end{equation*}
which are successively corrected till the norm of the residual is within a predefined tolerance by means of the Newton-Raphson method. The generic update reads as:
\begin{equation*}
    \dot{\hat{\mathbf{u}}}_{n+1} = \dot{\hat{\mathbf{u}}}_{n+1} +\mathbf{b},
\end{equation*}
\begin{equation*}
    \hat{\mathbf{u}}_{n+1} = \hat{\mathbf{u}}_{n+1} + \gamma \Delta t \, \mathbf{b},
\end{equation*}
being $\mathbf{b}$ the solution of the linear system
\begin{equation*}
    \mathbf{A} \mathbf{b} = - \mathbf{R},
\end{equation*}
where the tangent matrix $\mathbf{A}$ is the linearization of the residual with respect to $ \dot{\hat{\mathbf{u}}}_{n+1}$, given by
\begin{equation*}
    \mathbf{A} = \alpha_m \mathbf{M} + \alpha_f \gamma \Delta t \left( \mathbf{K}_\Delta -  \mathbf{K}_\partial -  \mathbf{K}_\partial^\top +  \mathbf{M}_N + \mathbf{K}_F \right).
\end{equation*}

In the equation above, $\mathbf{K}_F$ is the linearization of $\overline{\mathbf{F}}( \hat{\mathbf{u}})$ with respect to $\hat{\mathbf{u}}_{n+1}$, which reads as
\begin{equation*}
    \mathbf{K}_F = 
    \int_{\Omega} \nabla \mathbf{N}^\top \, F''(u_{n+ \alpha_f}) \, \nabla \mathbf{N}   \, d\mathbf{x}
    + \int_{\Omega} \nabla \mathbf{N}^\top \cdot \left( F'''(u_{n+ \alpha_f}) \, \nabla \mathbf{N} \hat{\mathbf{u}}_{n+ \alpha_f} \right) \, \mathbf{N} \, \hbox{\rm d}\mathbf{x},
\end{equation*}
where  $F''(u_{n+ \alpha_f})$ is computed by substituting $\hat{\mathbf{u}}_{n+ \alpha_f}$ in \eqref{eq:F''} and
\begin{equation*}
     F'''(u_{n+ \alpha_f}) = 6 \, \sigma \, \mathbf{N} \hat{\mathbf{u}}_{n+ \alpha_f}.
\end{equation*}

\section{Adaptive isogeometric methods on multi-patch geometries}
\label{sec:ada_mp} 

In this section, we consider the construction of $C^1$ splines on analysis suitable multi-patch geometries proposed in \cite{KaSaTa19a,KaSaTa19b} and its extension to the hierarchical spline model, recently introduced in \cite{BrGiKaVa20,BrGiKaVa23}. We also present algorithms to perform refinement and coarsening which ensure the properties of linear independence and admissibility on analysis suitable multi-patch configurations. The refinement algorithm for (truncated) hierarchical splines here considered was presented in \cite{BrGiKaVa23} together with complexity estimates, while the coarsening module on multi-patch geometries is newly introduced by extending the algorithm for the single-patch case proposed in \cite{THB-refinement-coarsening}.

\subsection{$C^1$ splines on analysis suitable geometries}

\begin{figure}[!th]
\centering
\includegraphics[width=0.35\textwidth]{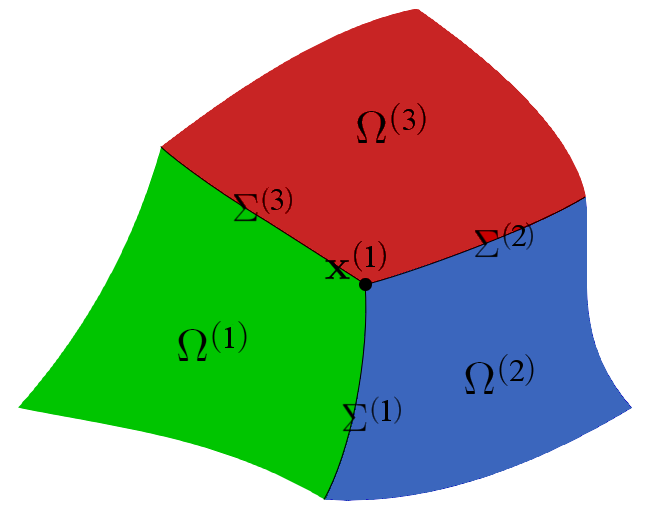}
\caption{Example of open three-patch domain, composed by quadrilateral patches $\Omega^{(i)}$, $i \in \indOmega$, inner edges~$\Sigma^{(i)}$, $i \in \indSigmaI$, and inner vertices~$\f{x}^{(i)}$, $i \in \indChiI$.}
\label{fig:exgeo}
\end{figure}

Let $\Omega \subset \R^2$ be an open multi-patch domain defined by the set of disjoint quadrilateral patches $\Omega^{(i)}$, $i \in \indOmega$, inner edges~$\Sigma^{(i)}$, $i \in \indSigmaI$, and inner vertices~$\f{x}^{(i)}$, $i \in \indChiI$ (see Figure \ref{fig:exgeo}). The boundary 
 of the multi-patch domain is the union of boundary edges 
 and boundary vertices. 
 For any $i \in \indOmega$, the open patch $\Omega^{(i)}$ is the image of a bijective and regular map $\f{F}^{(i)}: [0,1]^2 \rightarrow \overline{\Omega^{(i)}}$, which defines the spline parameterization $\f{F}^{(i)} \in (\US{p}{r} \otimes \US{p}{r})^2$, where $\US{p}{r}$ is the space of univariate splines  of degree~$p \geq 3$ and regularity~$1 \leq r \leq p-2$ in $[0,1]$. 
The geometry of the multi-patch domain $\Omega$ is described in terms of the single spline parameterizations $\f{F}^{(i)}$, $i \in \indOmega$.

By following \cite{CoSaTa16,KaSaTa19a}, we assume the multi-patch geometry to be analysis suitable (AS) $G^1$ and define $C^1$ isogeometric spaces with optimal 
polynomial reproduction properties \cite{CoSaTa16,KaSaTa17b}. For each inner edge~$\Sigma^{(i)}$, $i \in \indSigmaI$, of an analysis suitable $G^1$ geometry there exist linear functions
\begin{equation}\label{eq:gdata}
\alpha^{(i,0)},\quad \alpha^{(i,1)},\quad \beta^{(i,0)}, \quad \beta^{(i,1)},
\end{equation}
with $\alpha^{(i,0)}$ and $\alpha^{(i,1)}$ relatively prime, such that 
\begin{equation*}
 \alpha^{(i,0)}  (\xi) \alpha^{(i,1)}  (\xi) > 0, \quad
 \forall\, \xi \in [0,1],
\end{equation*}
and
\begin{equation*}
\alpha^{(i,0)} (\xi)  \Dv \f{F}^{(i_1)}(\xi,0)  +
        \alpha^{(i,1)}(\xi) \Du  \f{F}^{(i_0)}(0,\xi) + \beta^{(i)} (\xi)
        \Dv  \f{F}^{(i_0)}(0,\xi)  =\boldsymbol{0},
\end{equation*}
with
 $\beta^{(i)} (\xi)= \alpha^{(i,0)} (\xi) \beta^{(i,1)}(\xi)+ \alpha^{(i,1)}(\xi)\beta^{(i,0)}(\xi)$.
The computation of the gluing data \eqref{eq:gdata} is necessary for the construction of the $C^1$ basis functions, and it is detailed in \cite[Appendix A.2]{BrGiKaVa23}, see also \cite{CoSaTa16}. The simplest cases of analysis-suitable $G^1$ multi-patch geometries are piecewise bilinear parameterizations \cite{CoSaTa16,KaBuBeJu16,KaViJu15}, but geometries of higher degree can also be defined, and it is possible to re-parameterize non-AS $G^1$ multi-patch geometry into analysis-suitable $G^1$ configurations  \cite{KaSaTa17b,KaSaTa19b}.

\begin{figure}[!t]
\begin{subfigure}[Patch interior basis function in $\PhiB_\Omega$.]{
\includegraphics[trim=0cm 0cm 0cm 0cm, clip, width=0.46\textwidth]{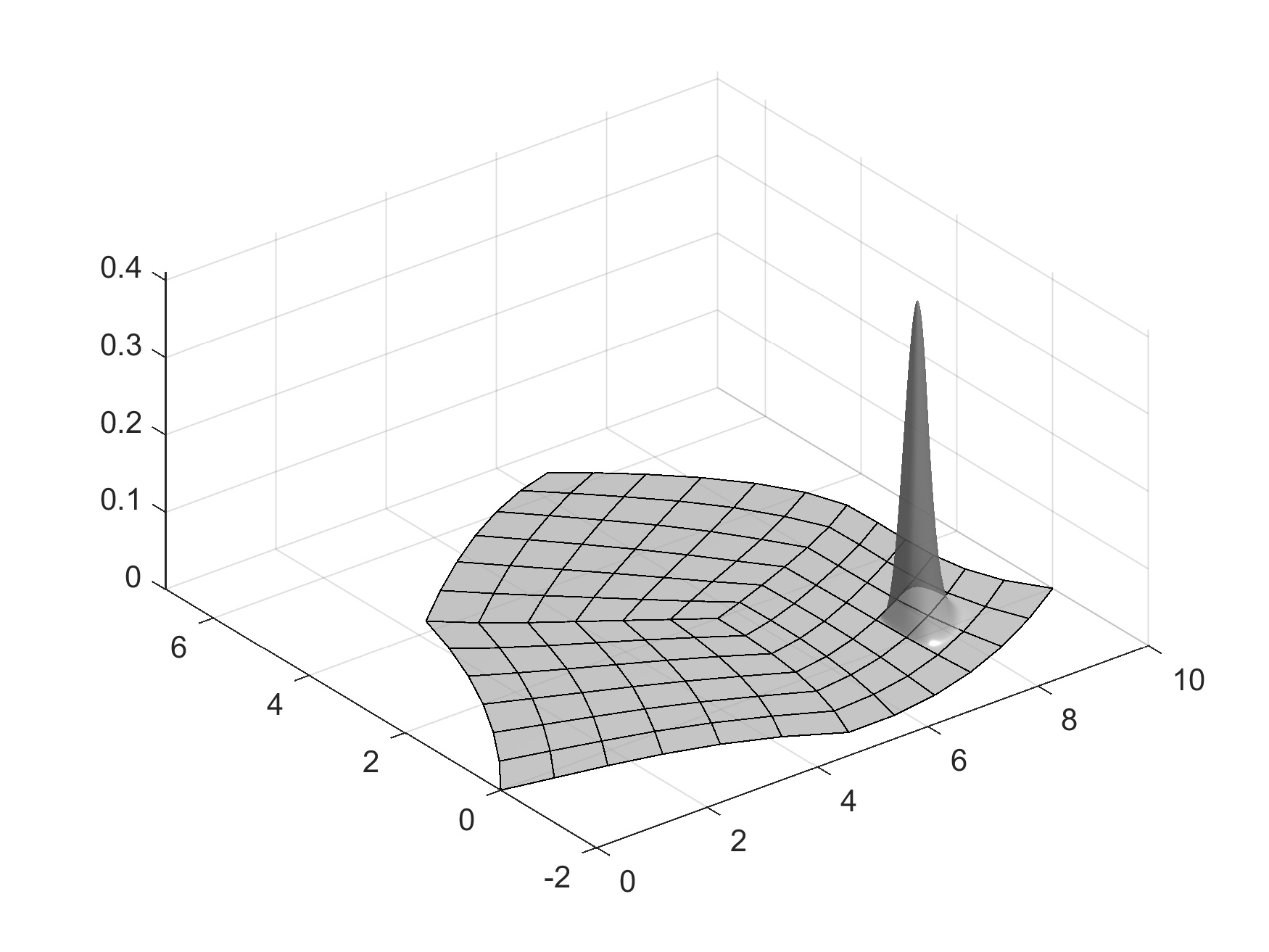}
}
\end{subfigure}
\begin{subfigure}[Edge basis function in $\PhiB_\Sigma$.]{
\includegraphics[trim=0cm 0cm 0cm 0cm, clip, width=0.46\textwidth]{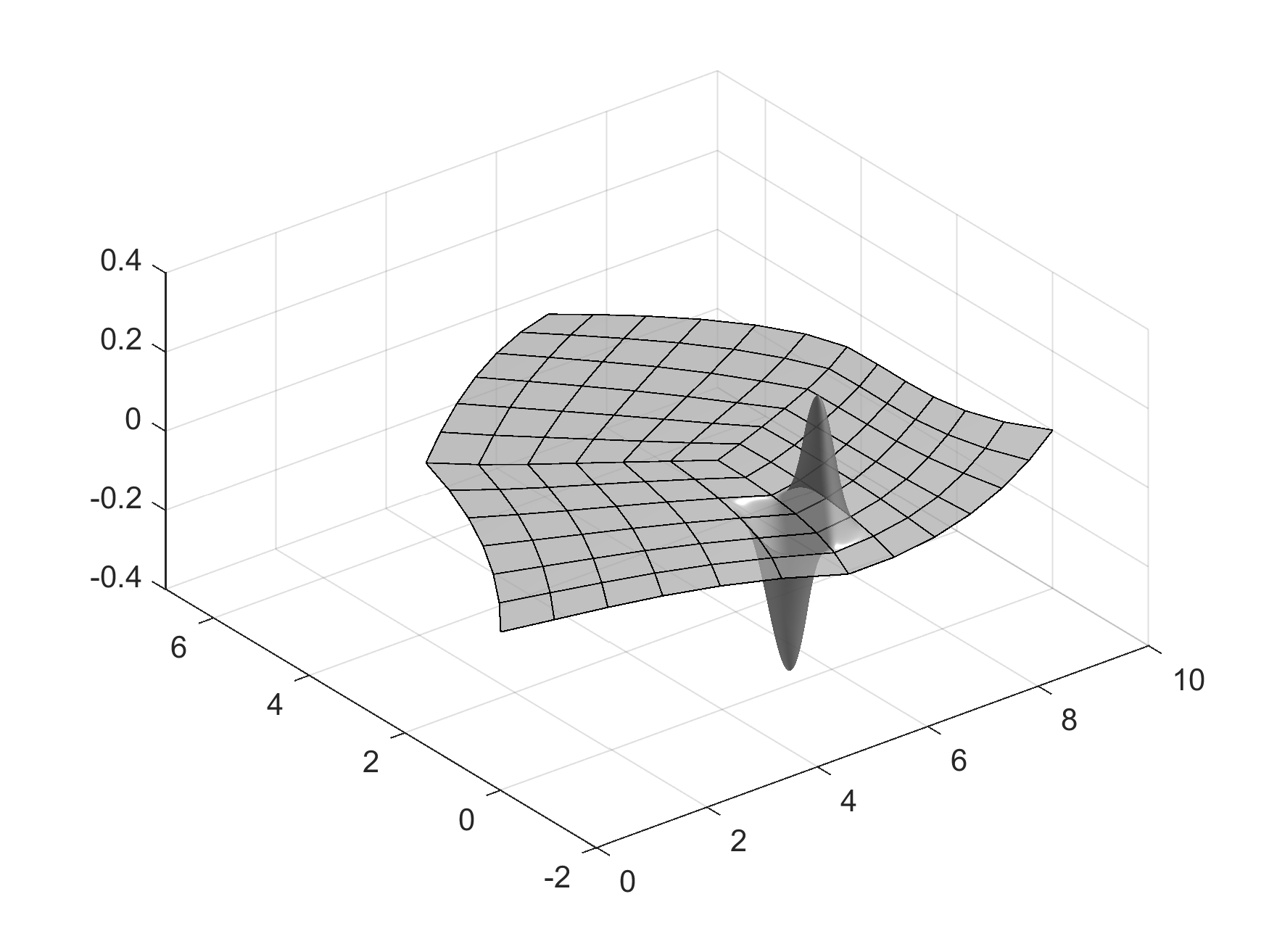}
}
\end{subfigure}
\begin{subfigure}[Edge basis function in $\PhiB_\Sigma$.]{
\includegraphics[trim=0cm 0cm 0cm 0cm, clip, width=0.46\textwidth]{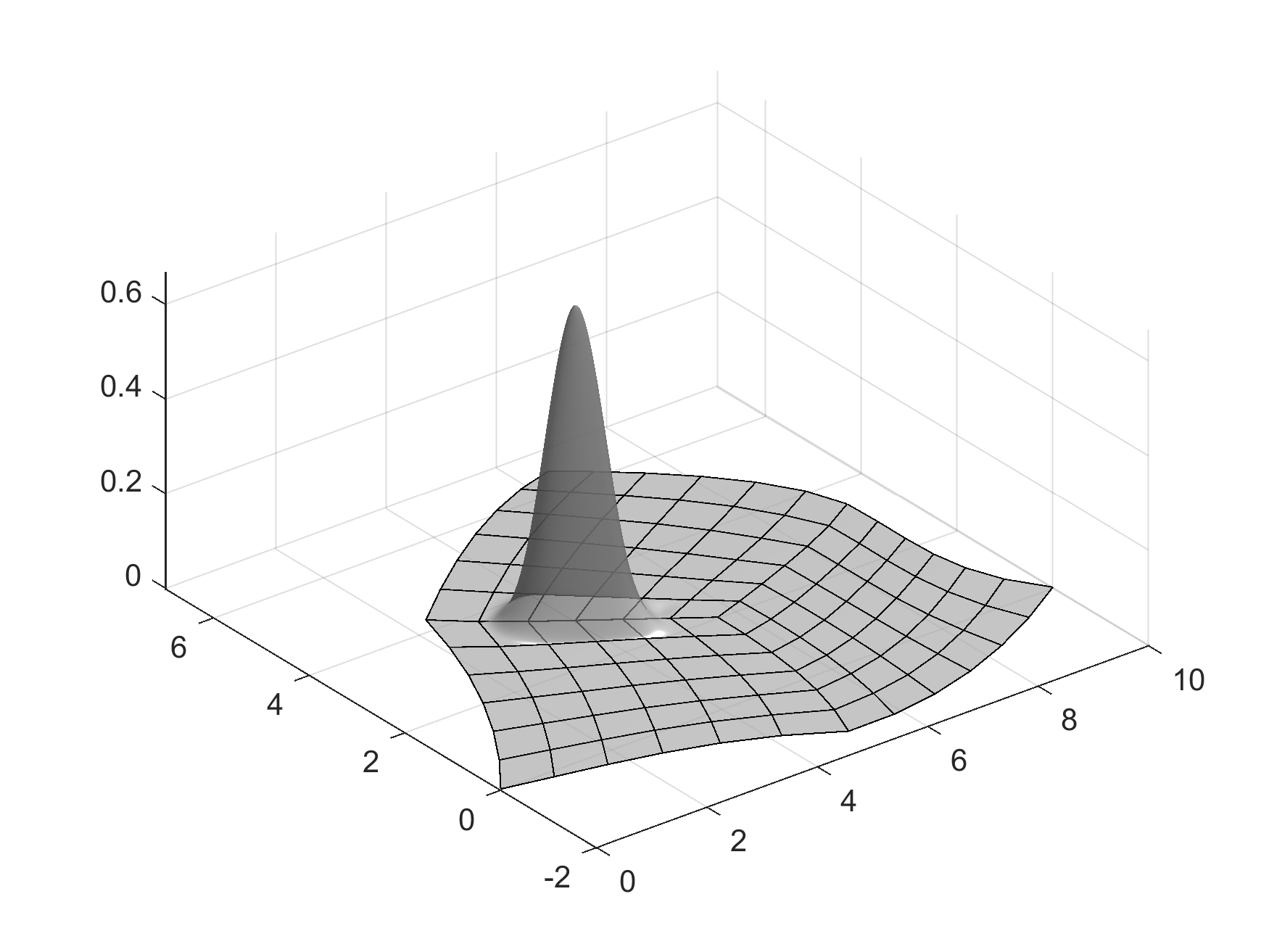}
}
\end{subfigure}
\begin{subfigure}[Vertex basis function in $\PhiB_\chi$.]{
\includegraphics[trim=0cm 0cm 0cm 0cm, clip, width=0.46\textwidth]{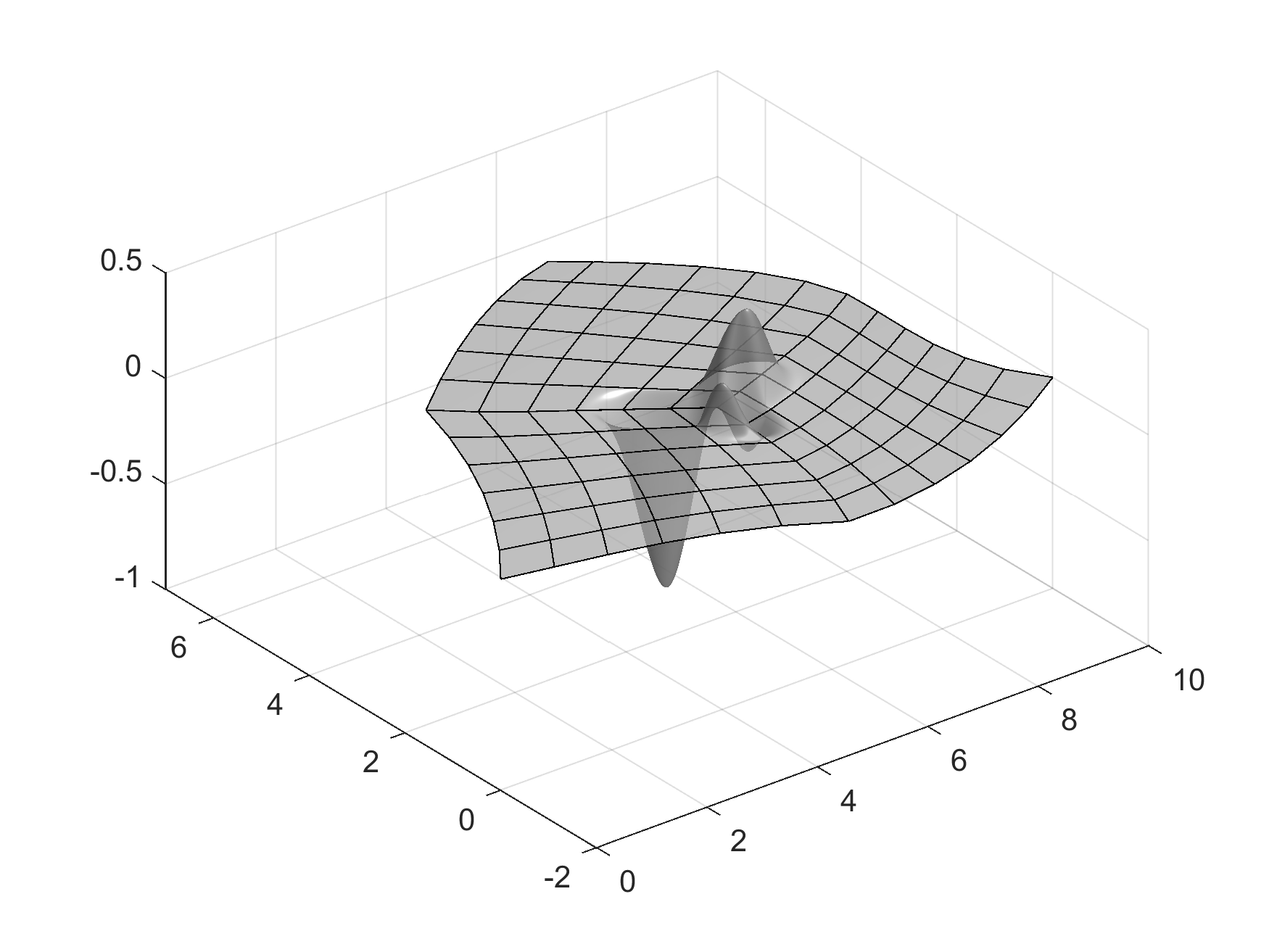}
}
\end{subfigure}

\caption{Examples of elements of the $C^1$ basis $\Phi$ defined on the geometry of Figure \ref{fig:exgeo}.
}
\label{fig:c1basis}
\end{figure}

We consider the discrete space~$\UW = \mathrm{span}\, \PhiB$, with
\begin{equation} 
\PhiB = \PhiB_\Omega \cup \Phi_\Sigma \cup \Phi_\chi, 
\quad\text{ and } \quad
\PhiB_\Omega = \bigcup_{i \in \indOmega} \PhiB_{\Omega^{(i)}}, \quad
\PhiB_\Sigma = \bigcup_{i \in \indSigma} \PhiB_{\Sigma^{(i)}}, \quad
\PhiB_\chi = \bigcup_{i \in \indChi} \PhiB_{\f{x}^{(i)}}, \label{eq:basis_A}
\end{equation}
where the $\PhiB_\Omega$, $\PhiB_\Sigma$, and $\PhiB_\chi$ identify patch interior, edge, and vertex basis functions, respectively. The basis functions in any of these sets are globally $C^1$-smooth on the multi-patch domain $\Omega$. In addition, the vertex functions in $\PhiB_\chi$ are also $C^2$-smooth at the corresponding 
vertex~$\f{x}^{(i)}$. The patch interior basis functions in $\PhiB_\Omega$ are ordinary B-splines with zero values and zero derivatives on every edge and vertex.  The edge basis functions in $\PhiB_\Sigma$ are nonzero on mesh elements of two adjacent patches of a certain edge, or of a single patch for boundary edges. Finally, vertex basis functions in $\PhiB_\chi$ are nonzero on mesh elements of the patches surrounding a certain vertex. Examples for the three kinds of functions are given in Figure~\ref{fig:c1basis}. Note that, in view of the $C^2$ interpolation condition considered in the construction for any vertex, there is always a fixed number of six vertex functions associated to a vertex, independently of the vertex valence. We refer to \cite{KaSaTa19a} for the details on the basis functions construction by only highlighting that the three different kinds of basis functions restricted to the patch $\Omega^{(i)}$ can be expressed as linear combinations of the (mapped) B-splines $\US{p}{r} \otimes \US{p}{r}$, see also \cite{BrGiKaVa20,BrGiKaVa23} where the representation of the edge and vertex functions in terms of standard B-splines is further detailed. Following the notation of Section~\ref{sec:phase_field_discretization}, this gives the existence of a matrix $\mathbf{C}$
\[
\mathbf{N}^\top = \mathbf{C} \mathbf{B}^\top,
\]
where $\mathbf{N}$ is the row vector of $C^1$ basis functions, and $\mathbf{B}$ is the row vector of standard B-splines for the whole domain, and discontinuous between patches. This expression allows to assemble the matrices from single patch contributions. Obviously, it is possible to add a further step to obtain the B\'ezier extraction \cite{BoScEvHu11}, that is, to express the $C^1$ basis functions in terms of Bernstein polynomials on each element.

\subsection{$C^1$ hierarchical splines on analysis suitable geometries}

We consider a nested sequence of $C^1$ spline spaces $\UW^0\subset \UW^1\subset\ldots\subset\UW^{N-1}$ defined on the multi-patch domain $\Omega$ described with an analysis-suitable $G^1$ parameterization. For $\ell=0,1,\ldots,N-1$, each space $ \UW^\ell$ is spanned by the basis $\PhiB^\ell$ composed by patch interior, edge, and vertex functions, as in \eqref{eq:basis_A} with respect to the tensor-product grid $\GG^\ell$ of level $\ell$. Note that, when the domain is composed of one single patch, the construction considers as underlying basis functions at each hierarchical level only patch interior B-splines \cite{vuong2011}.

\begin{figure}[!t]
\begin{subfigure}[Tensor-product grid $G^0$.]{
\includegraphics[trim=0cm 0cm 0cm 0cm, clip, width=0.30\textwidth]{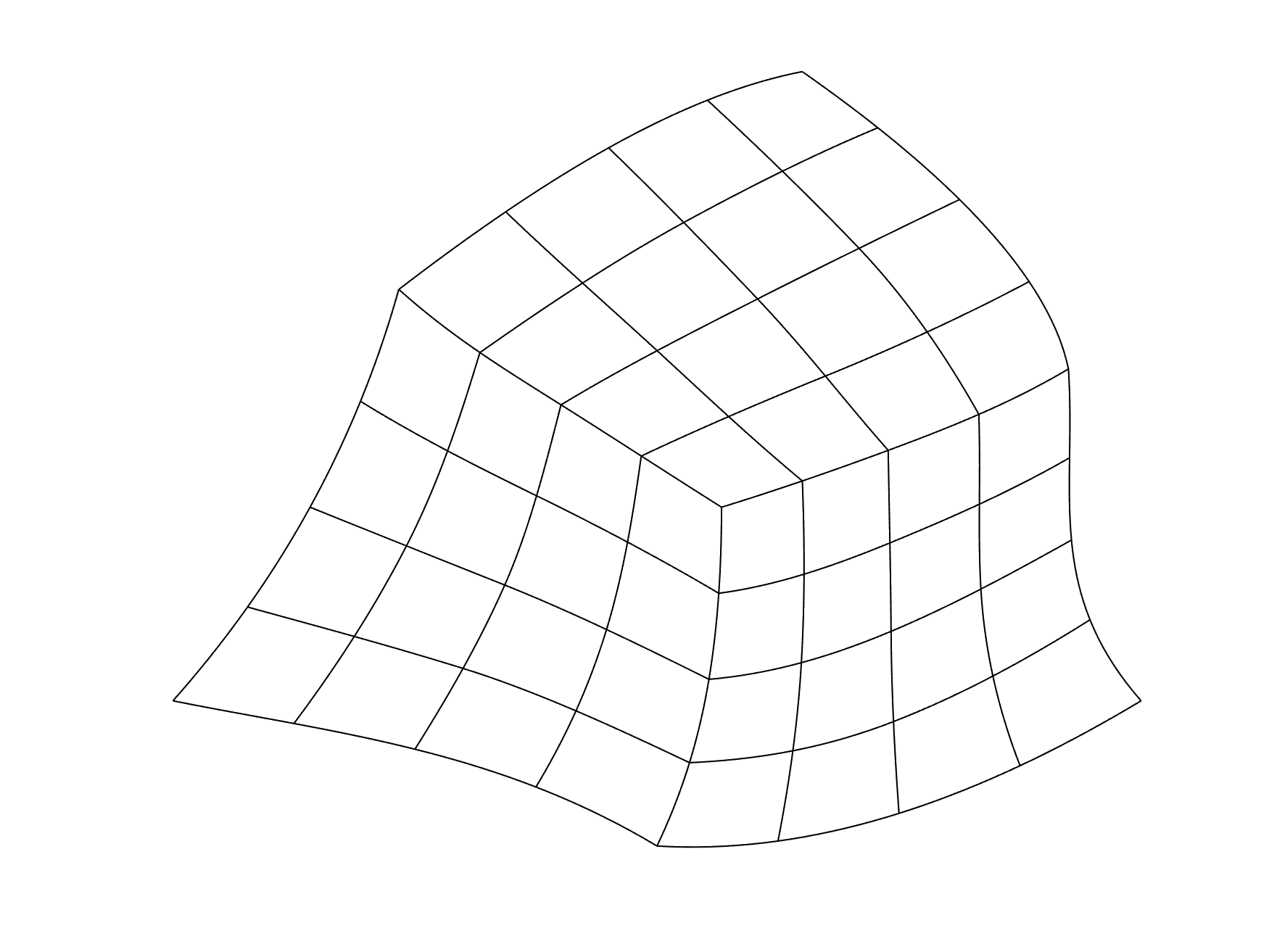}
}
\end{subfigure}
\begin{subfigure}[Tensor-product grid $G^1$.]{
\includegraphics[trim=0cm 0cm 0cm 0cm, clip, width=0.30\textwidth]{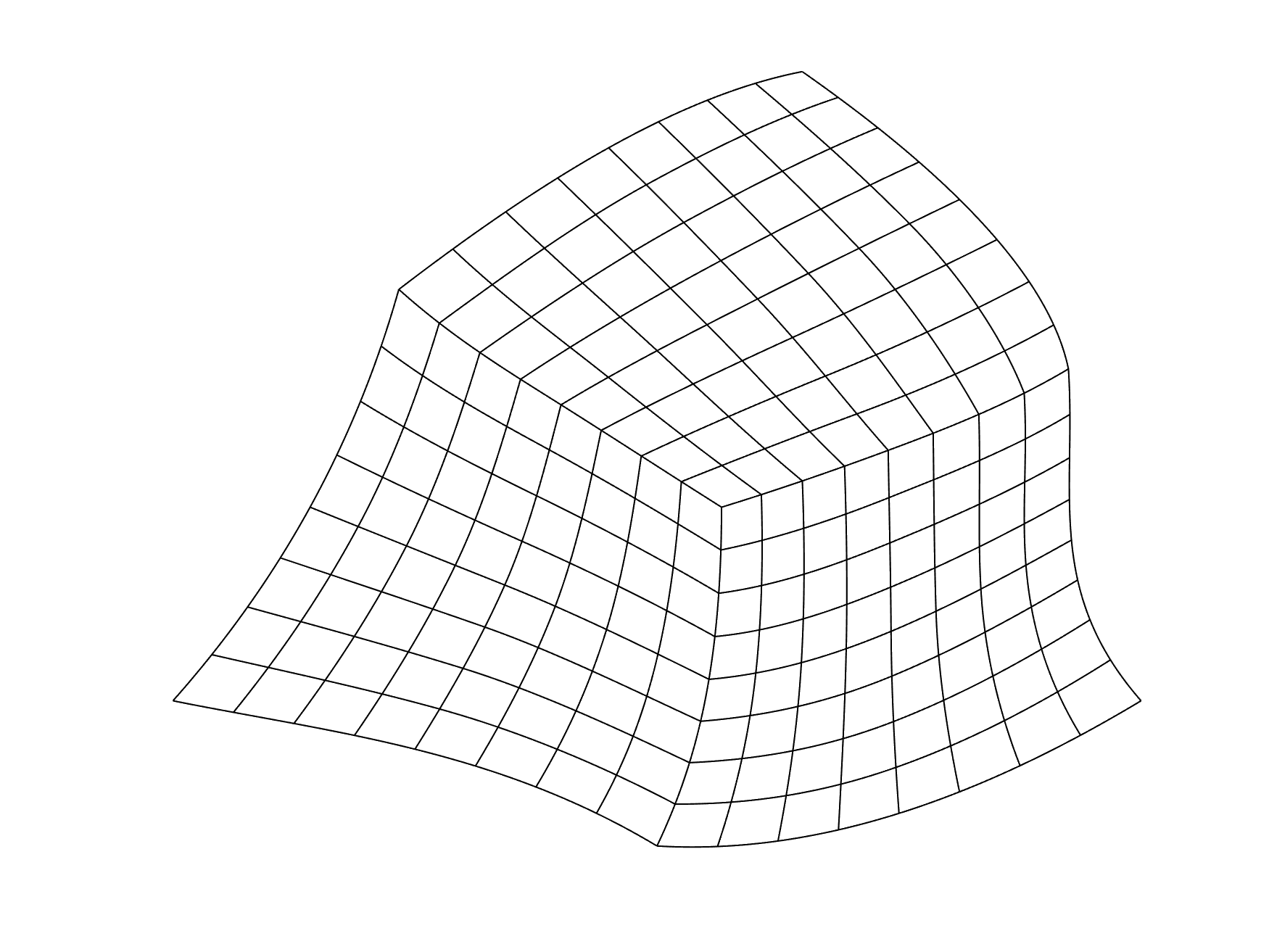}
}
\end{subfigure}
\begin{subfigure}[Tensor-product grid $G^2$.]{
\includegraphics[trim=0cm 0cm 0cm 0cm, clip, width=0.30\textwidth]{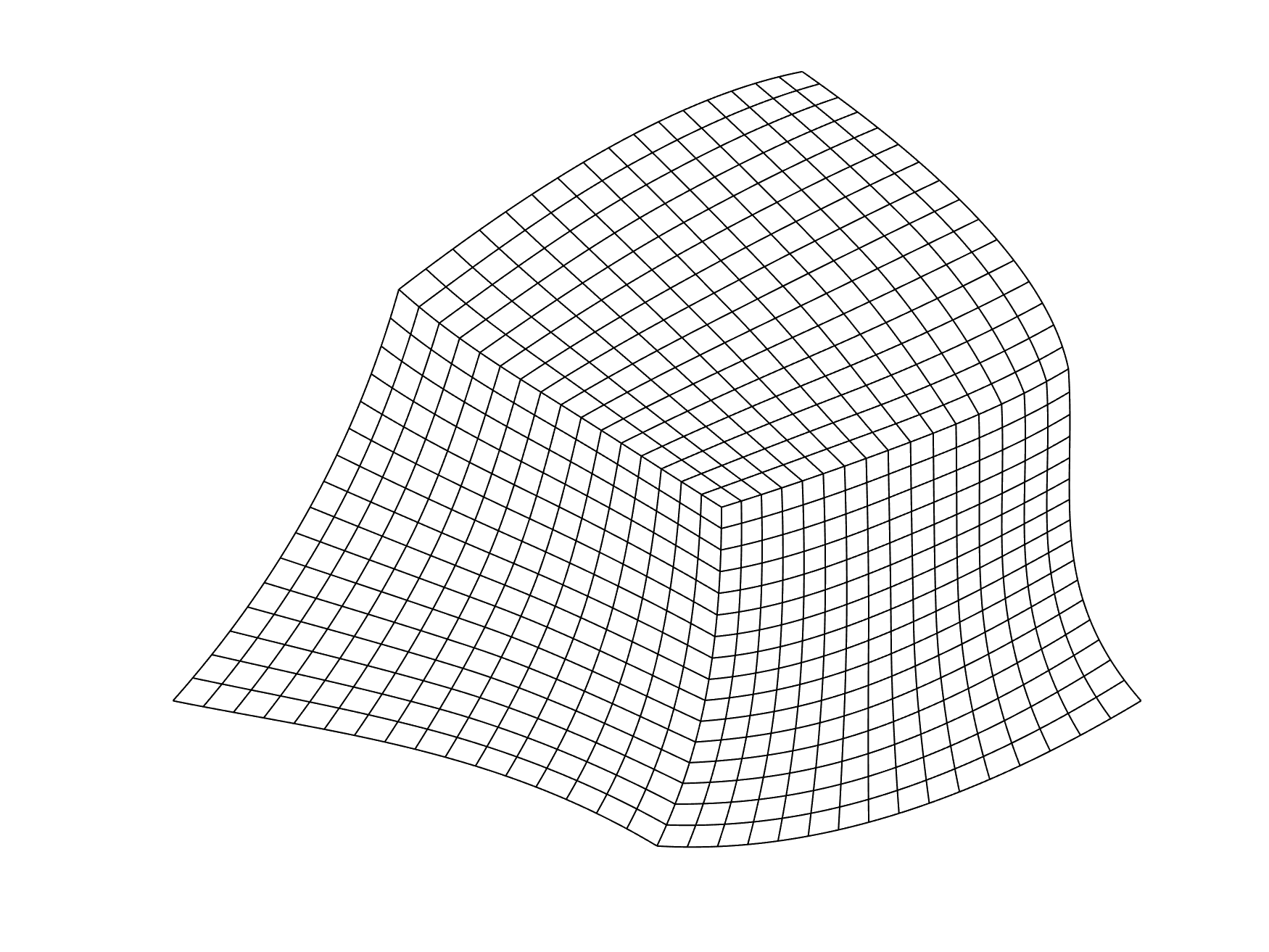}
}
\end{subfigure}
\begin{subfigure}[Domain $\Omega^0$.]{
\includegraphics[trim=0cm 0cm 0cm 0cm, clip, width=0.30\textwidth]{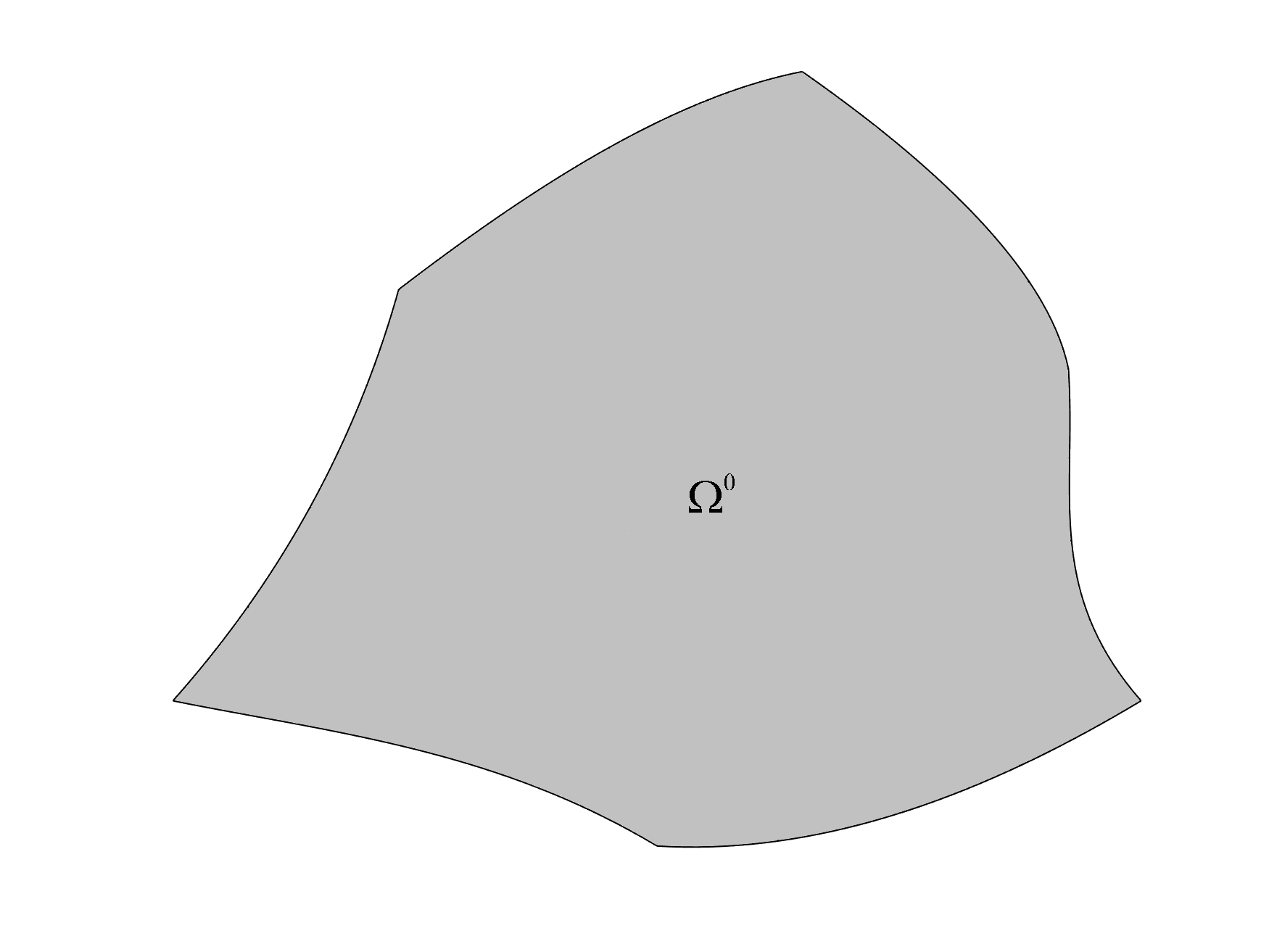}
}
\end{subfigure}
\begin{subfigure}[Domain $\Omega^1$.]{
\includegraphics[trim=0cm 0cm 0cm 0cm, clip, width=0.30\textwidth]{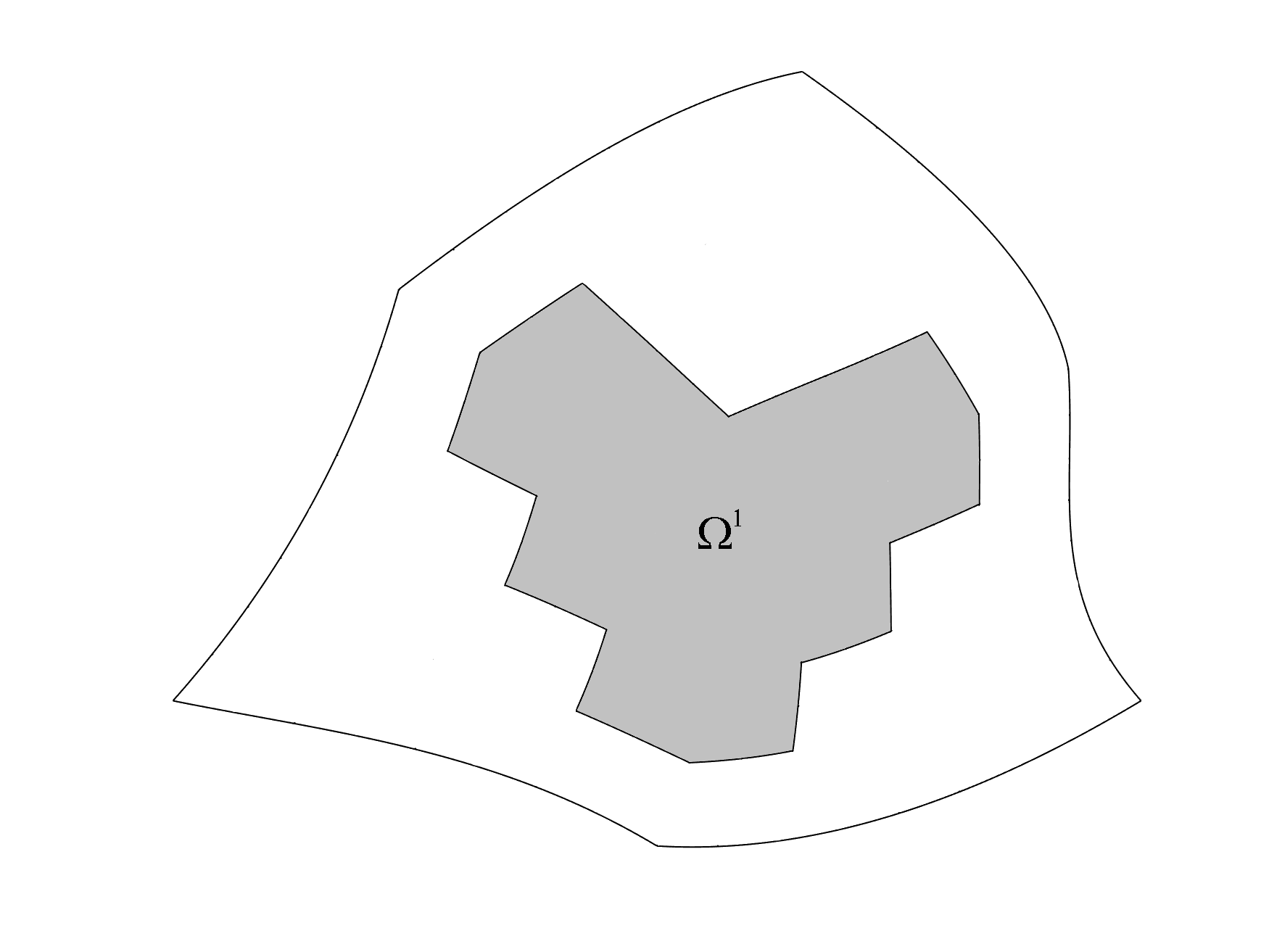}
}
\end{subfigure}
\begin{subfigure}[Domain $\Omega^2$.]{
\includegraphics[trim=0cm 0cm 0cm 0cm, clip, width=0.30\textwidth]{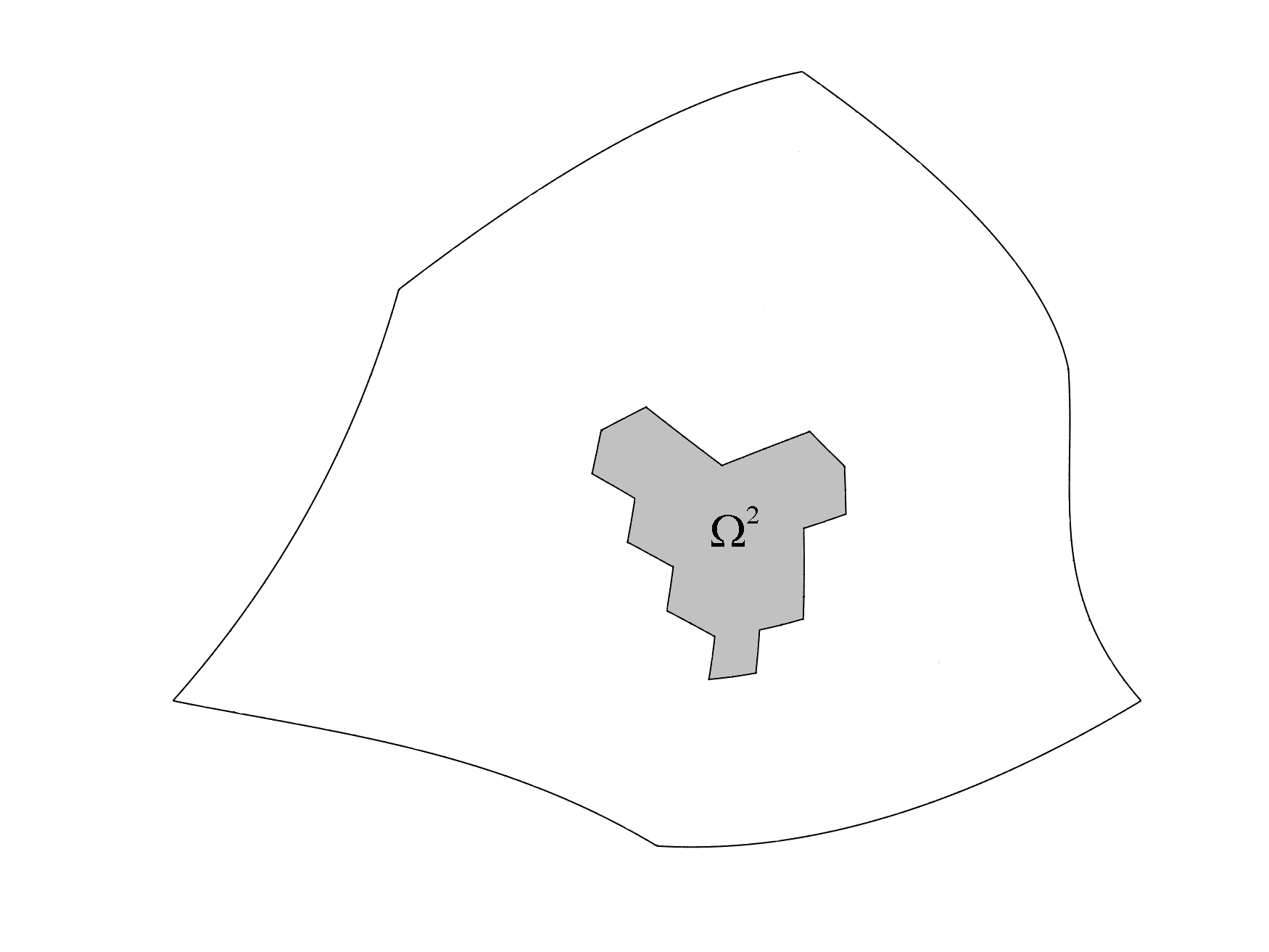}
}
\end{subfigure}
\begin{subfigure}[Tensor-product grid $G^0$ on $\Omega^0$.]{
\includegraphics[trim=0cm 0cm 0cm 0cm, clip, width=0.30\textwidth]{figures/ex_hmesh_tp_a.png}
}
\end{subfigure}
\begin{subfigure}[Tensor-product grid $G^1$ on $\Omega^1$.]{
\includegraphics[trim=0cm 0cm 0cm 0cm, clip, width=0.30\textwidth]{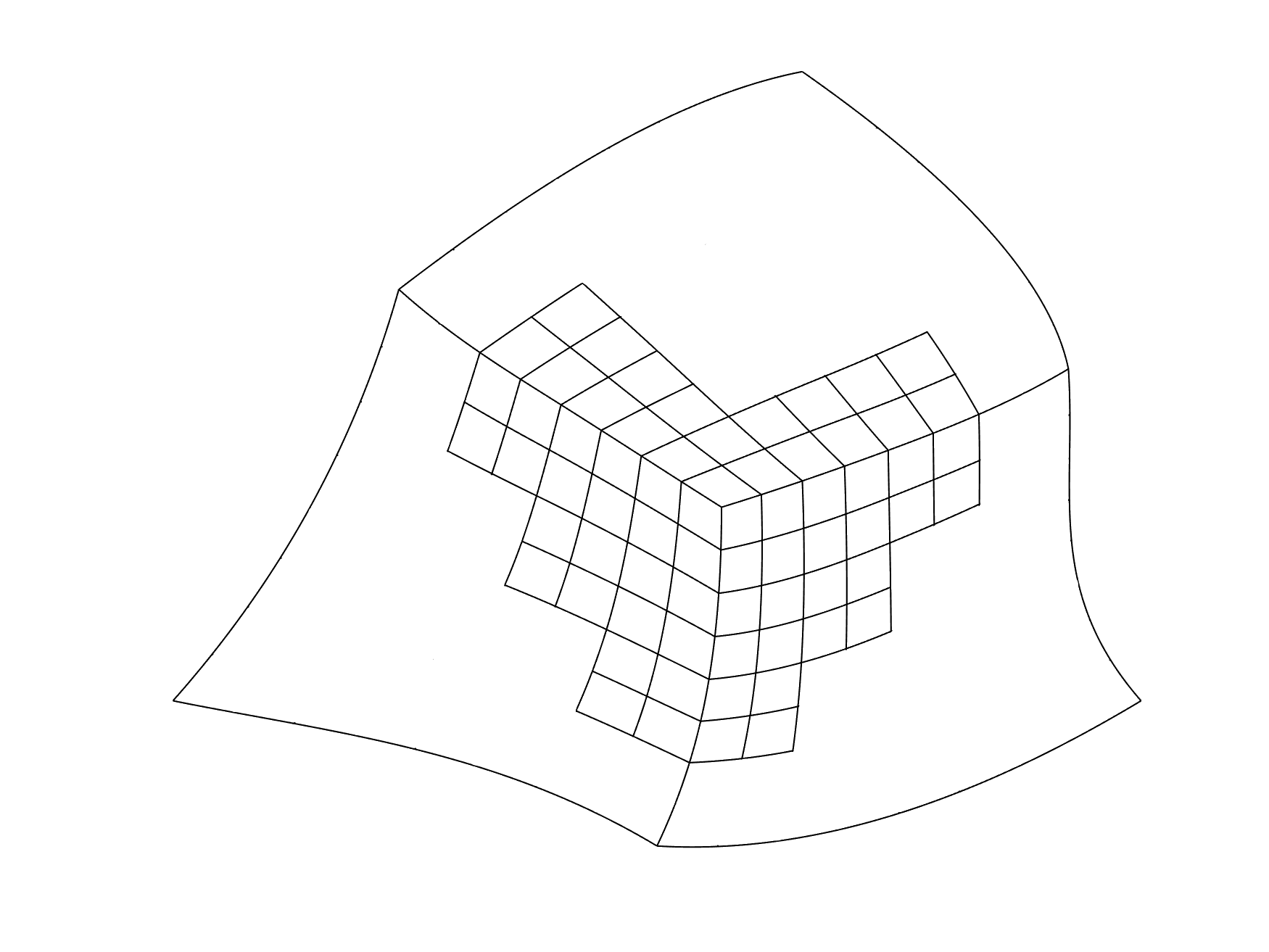}
}
\end{subfigure}
\begin{subfigure}[Tensor-product grid $G^2$ on $\Omega^2$.]{
\includegraphics[trim=0cm 0cm 0cm 0cm, clip, width=0.30\textwidth]{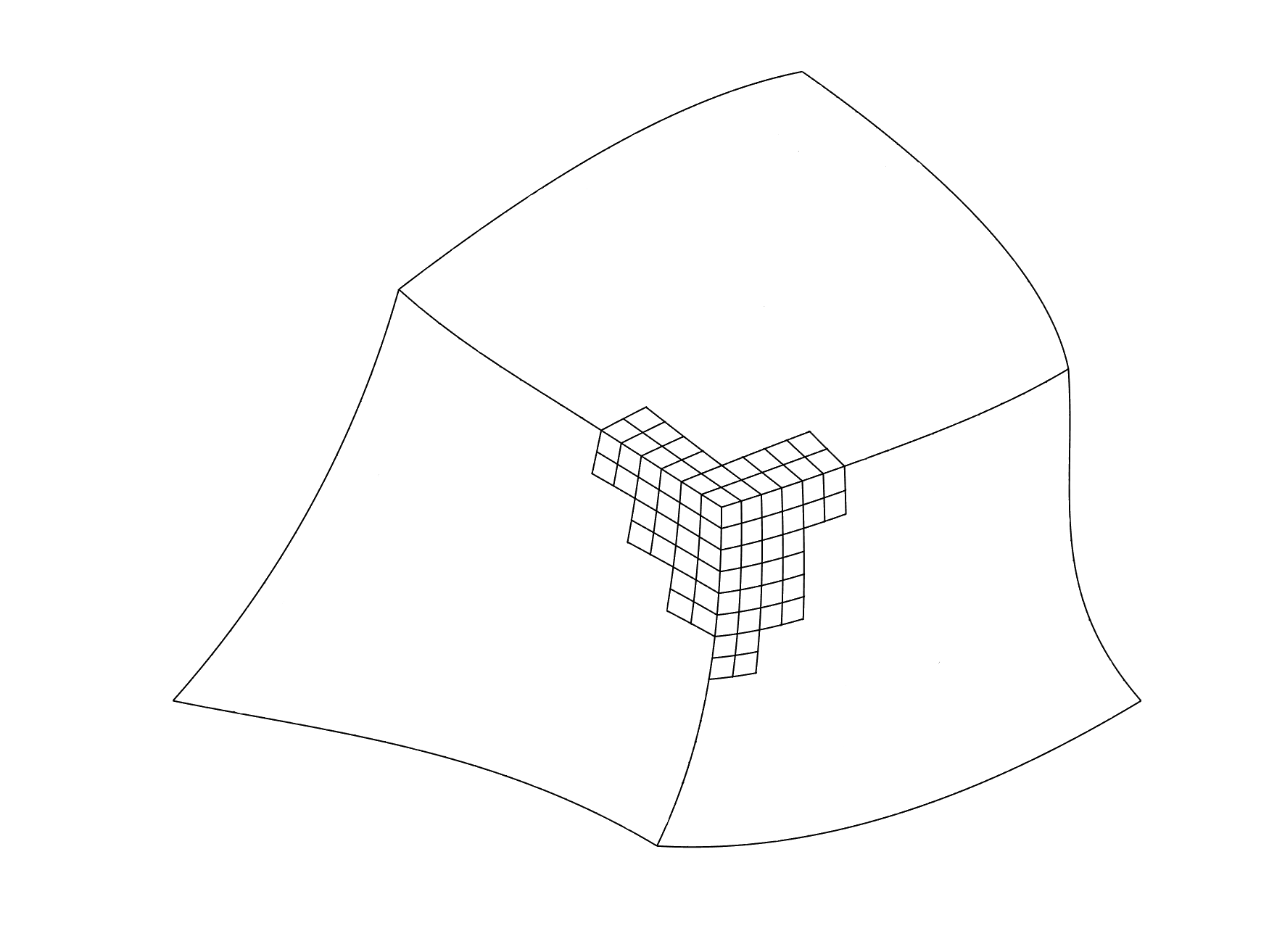}
}
\end{subfigure}
\begin{center}
\begin{subfigure}[Resulting hierarchical mesh.]{
\includegraphics[trim=0cm 0cm 0cm 0cm, clip, width=0.30\textwidth]{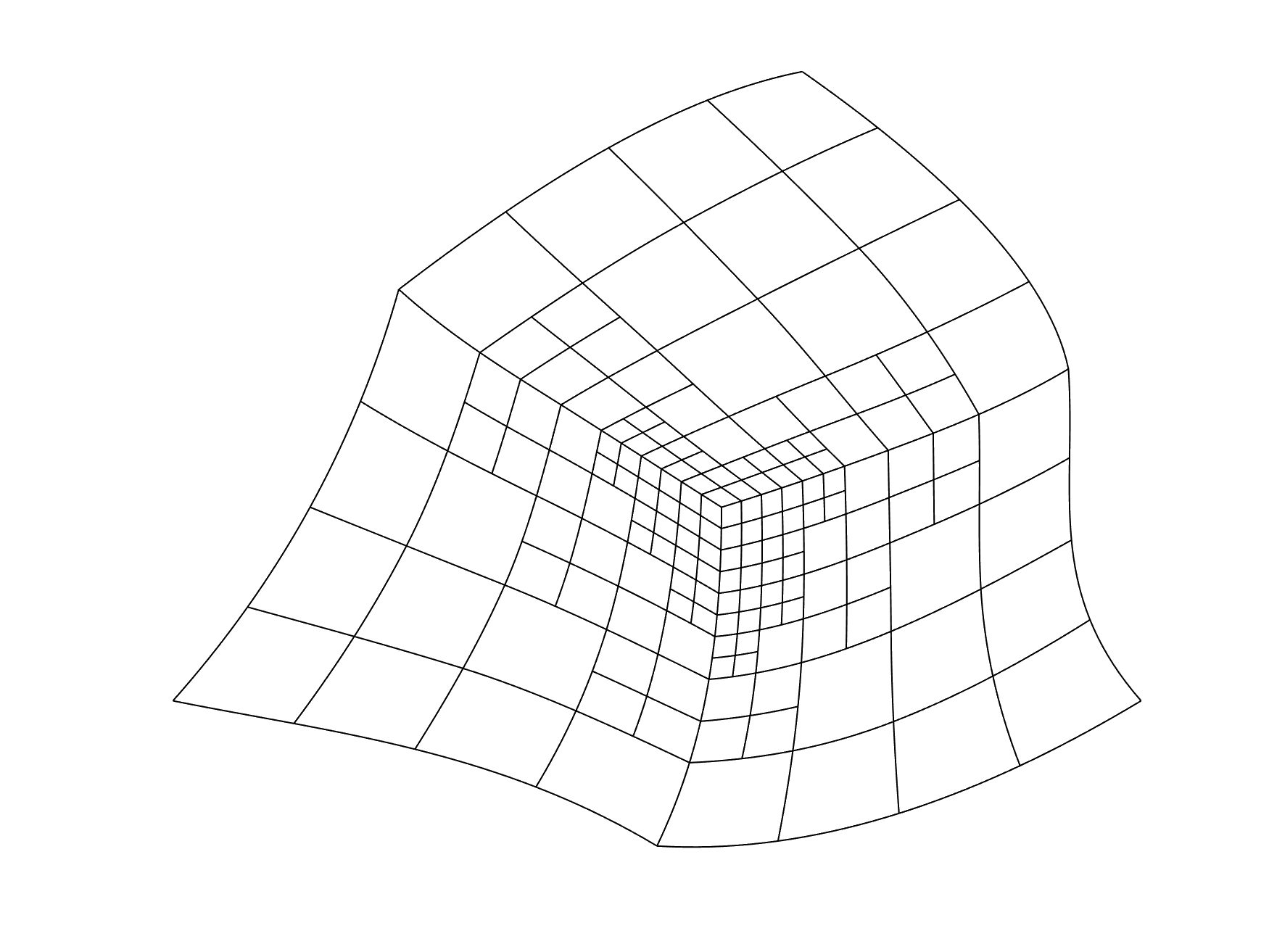}
}
\end{subfigure}
\end{center}

\caption{Construction of a hierarchical mesh with three levels.}
\label{fig:hierarchical_mesh}
\end{figure}

The set of hierarchical splines $\cal H$ is defined by iteratively activating basis functions at refined levels, while simultaneously de-activating coarse basis functions whose support is fully covered by the newly introduced hierarchical splines. The definition is the following:
\begin{equation*}\label{eq:hbasis}
{\cal H}:=\left\{\phi\in\Phi^\ell: {\supp}^0\phi\subseteq \Omega^\ell \wedge
{\supp}^0\phi\not\subseteq\Omega^{\ell+1}, 
\ell=0,\ldots,N-1
\right\},
\end{equation*}
where  $\supp^0 \phi := \supp\phi \cap \Omega^0$ and ${\Omega}^0\supseteq{\Omega}^1 \supseteq \ldots \supseteq{\Omega}^{N-1}$, $\Omega^0\subseteq \Omega$, is a sequence of closed nested domains which identifies the refinement regions at different resolution levels, see Figure~\ref{fig:hierarchical_mesh}.

In order to reduce the interaction between hierarchical functions introduced at different refinement levels and recover the partition of unity property of standard B-splines, the truncated basis for hierarchical B-splines was introduced in \cite{giannelli2012}, see also \cite{giannelli2014,giannelli2016}. We here apply the same idea to define the set of truncated hierarchical splines
\begin{equation*}\label{eq:thbasis}
{\cal T}:=\left\{{\Trunc}^{\ell+1}(\phi): \phi \in \Psi^\ell \cap {\cal H}, 
\ell=0,\ldots,N-1
\right\},
\end{equation*}
with
\begin{equation*}\label{eq:succtrunc}
{\Trunc}^{\ell+1}(\phi) := {\trunc}^{N-1}\left(\ldots\left(
{\trunc}^{\ell+1}\left(\phi\right)\right)\ldots\right)
\end{equation*}
and ${\Trunc}^{N}(\phi)=\phi$, where, for any spline $s\in\Phi^{\ell}$, we exploit the refinement mask with respect to the refined basis $\Phi^{\ell+1}$ presented in \cite{BrGiKaVa23} to express $s$ as	
\[
s = \sum_{\phi\in \Phi^{\ell+1}} c_{\phi}^{\ell+1}(s) \phi,
\]
and define the truncation of $s$ with respect to level $\ell+1$ as
\begin{equation*}
{\trunc}^{\ell+1} (s) := \sum_{\phi\in \Phi^{\ell+1},\, \supp^0 \phi\not\subseteq\Omega^{\ell+1}} c_{\phi}^{\ell+1}(s) \phi. 
\end{equation*}

In the single-patch case and in the multi-patch case with $C^0$ continuity, THB-splines are linearly independent due to the local linear independence of B-splines of one level, see \cite{giannelli2014}. Unfortunately, the property of local linear independence is not valid for the $C^1$ splines of \cite{KaSaTa19b} that we described in the previous setting, and the standard hypothesis of local linear independence for the basis $\Phi^\ell$ of each level, usually assumed in the hierarchical spline model, is substituted in \cite{BrGiKaVa23} by a weaker condition, which still ensures the linear independence of (truncated) hierarchical splines. For $C^1$ hierarchical multi-patch splines, the condition of linear independence imposes a constraint on the mesh in the vicinity of the vertices, that was introduced in \cite{BrGiKaVa23} and we reproduce here for the sake of completeness.
\begin{theorem}\label{thm:li}
If for every active vertex function of level $\ell$, associated to the vertex $\mathbf{x}^{(i)}$, there exists an active element of level $\ell$ adjacent to the vertex $\mathbf{x}^{(i)}$, both $\mathcal{H}$ and $\mathcal{T}$ are linearly independent.
\end{theorem}
In other words, for every active vertex function of level $\ell$ there must be an active element of level $\ell$ adjacent to the same vertex.In the next section we present a refinement algorithm that respects this condition by suitably refining adjacent elements near any vertex, and a coarsening algorithm that also respects it by checking when an element close to a vertex can be reactivated. For this, we need to introduce for every level $\ell$ the following set of elements:
\begin{equation} \label{eq:set_of_unrefined}
{\cal A_{\ell}}:=\bigcup_{k=0,...,\ell}\{Q\in G^k:\, Q\not\subset \Omega_{\ell+1}\},
\end{equation}
which consists of elements up to level $\ell$ that have never been refined, that is, they are either active or have not been activated.



\subsection{Suitably graded hierarchical spline refinement and coarsening}
For the use of adaptive methods it is important to maintain a certain grading of the mesh, to guarantee the stability and the good conditioning of the numerical method, and also to preserve the accuracy of the approximation during coarsening. Refinement algorithms that preserve a suitable grading of the mesh for THB-splines were introduced in \cite{buffa2016c,buffa2017b}, under the name of admissible refinement, by requiring that on any active element only functions of $\mu$ different levels can be active, with $\mu > 1$ an integer parameter. In practice, whenever an element of level $\ell$ is marked for refinement, some elements of level $\ell-\mu+1$ on its neighborhood are also marked. In \cite{bgmp16} the algorithm was proven to have linear complexity. The same grading condition for THB-splines was then extended to a coarsening algorithm in \cite{THB-refinement-coarsening}, and applied to the analysis of the transient heat equation.

\begin{figure}[!t]
\begin{subfigure}{
\includegraphics[trim=0cm 0cm 0cm 0cm, clip, width=0.46\textwidth]{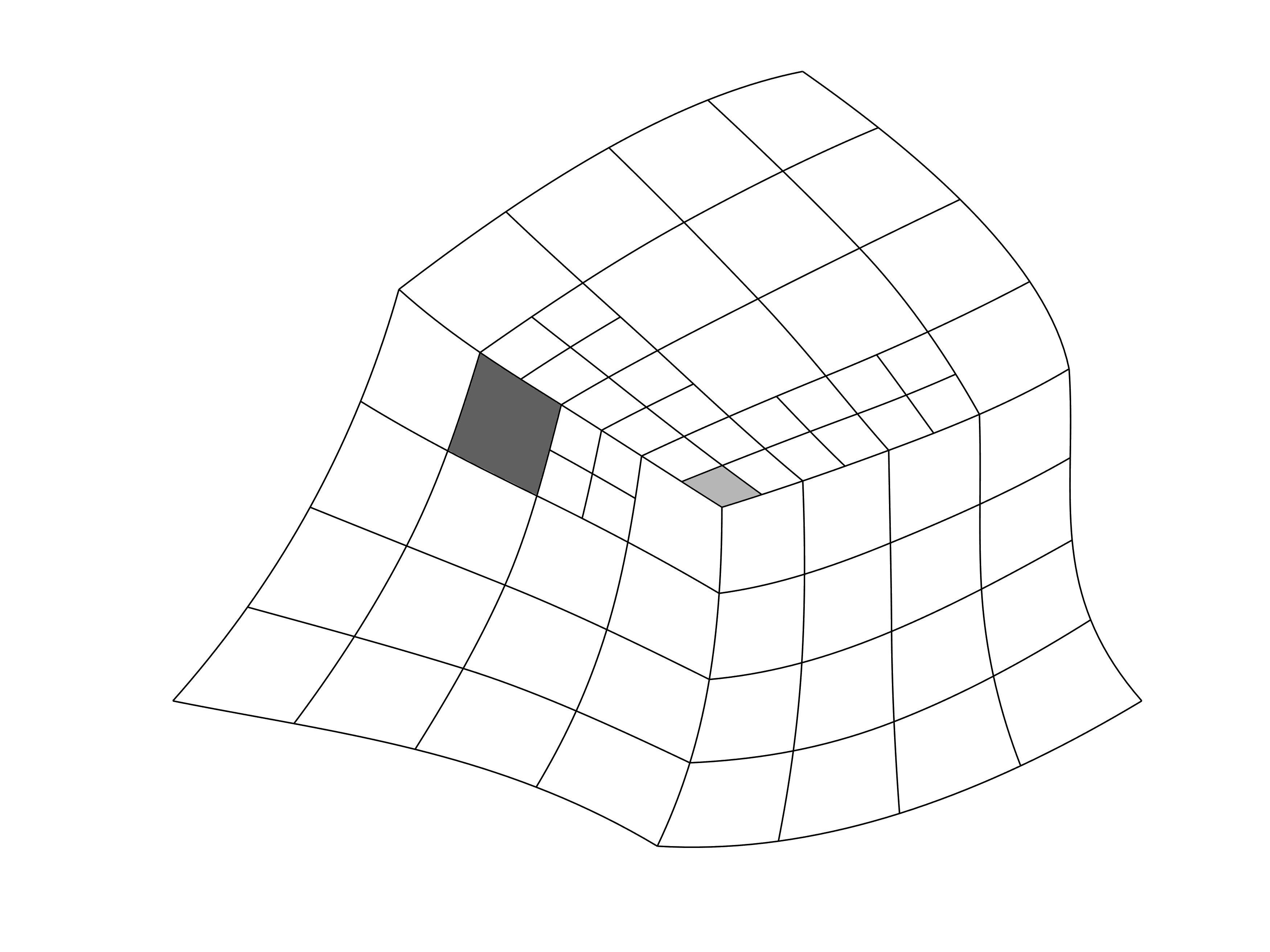}
}
\end{subfigure}
\begin{subfigure}{
\includegraphics[trim=0cm 0cm 0cm 0cm, clip, width=0.46\textwidth]{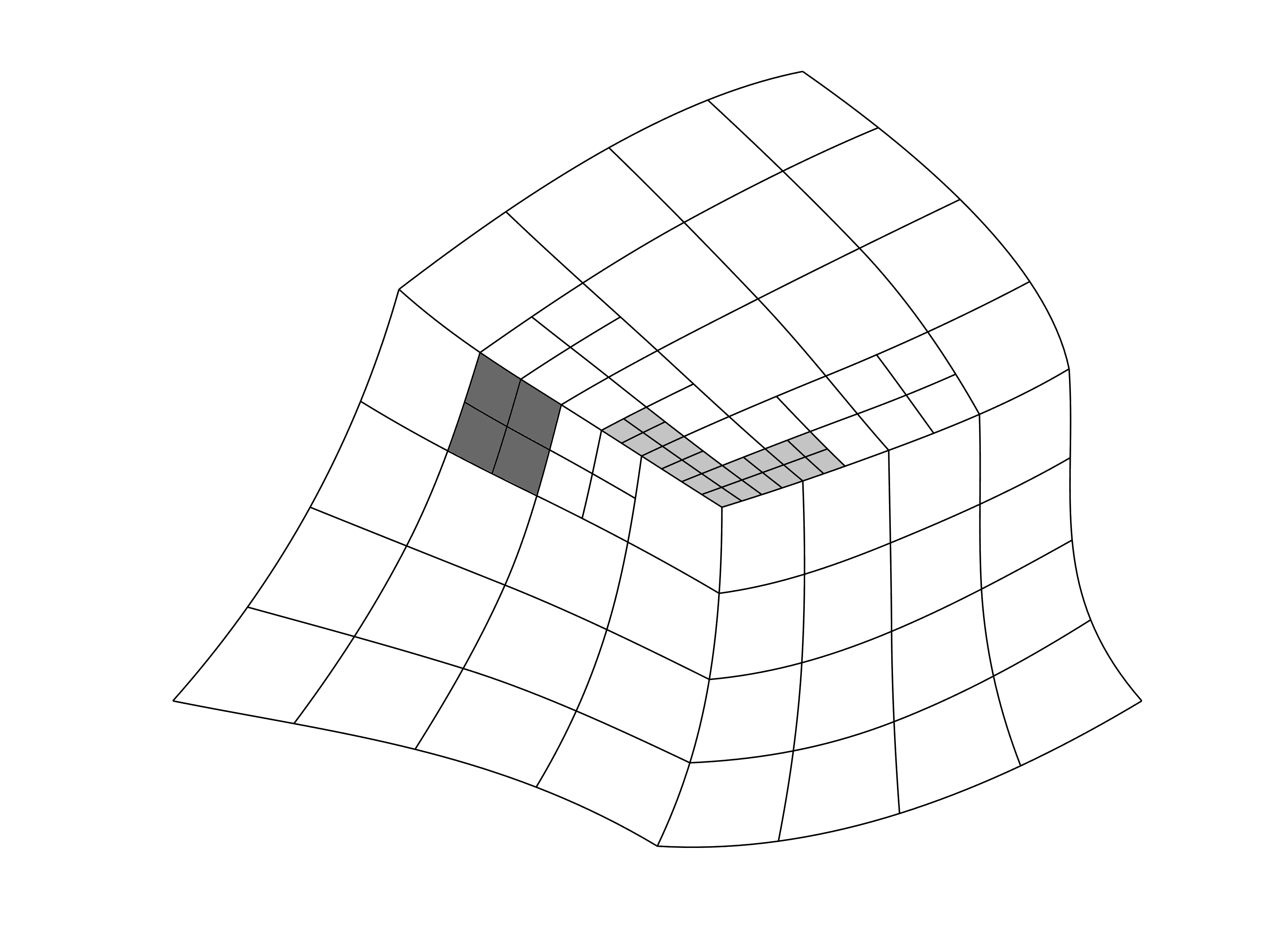}
}
\end{subfigure}
\caption{Example of refinement: when an element adjacent to a vertex is marked for refinement (light grey on left figure), some other elements have to be refined (light grey on right figure); when the marked element is away from the vertex (dark grey), no further elements have to be refined to ensure linear independence.}
\label{fig:refinement}
\end{figure}

In the case of $C^1$ hierarchical splines on multi-patch geometries, refinement and coarsening algorithms must not only maintain the grading of the mesh, but also respect the condition of linear independence of the basis functions. A refinement algorithm that fulfills both properties has been recently introduced in \cite{BrGiKaVa23}, where it is also proven to have linear complexity. To ensure linear independence, whenever an element adjacent to a vertex is marked for refinement, some other elements of the same level on its neighborhood are also marked. More precisely, for an element $Q$ of level $\ell$ which is adjacent to a vertex $\mathbf{x}^{(i)}$, we define the vertex-patch neighborhood $\mathcal{N}_\chi(Q)$ as the set of active elements of level $\ell$ belonging to the same patch as $Q$, and contained in the support of vertex functions associated to $\mathbf{x}^{(i)}$. Then, if $Q$ is marked for refinement, the elements in $\mathcal{N}_\chi(Q)$ are also marked. The suitable grading of the mesh is obtained as for THB-splines by marking some other elements of a coarser level $\ell-\mu +1$, in a neighborhood that we denote by $\mathcal{N}_r(Q,\mu)$, see \cite[Section 5.4]{BrGiKaVa23} for details. For completeness, we describe in Algorithm~\ref{REFINE_ADMISSIBLE} the way we implemented this, which differs from previous references in that the recursive algorithm is replaced by a loop starting from the finest level. It is important to note that the two neighborhoods $\mathcal{N}_\chi(Q)$ and $\mathcal{N}_r(Q,\mu)$ must be added in that precise order, to guarantee that the admissibility check is done for all elements. It is also worth to note that the maximum number of elements in $\mathcal{N}_\chi(Q)$ is five, counting $Q$ itself, in the case of regularity $r = p-2$. This number is reduced for low regularity, and becomes equal to one for regularity $r < p-3$, in which case no further refinement is necessary. An example of how the vertex-patch neighborhood affects the refinement is given in Figure~\ref{fig:refinement}: when an element adjacent to a vertex is refined, some other elements of its level are also refined; instead, if the element is not adjacent to a vertex, no further elements are refined, even if it is contained in the support of a vertex basis function.

\begin{algorithm}[h]
\caption{\texttt{Refine} \label{REFINE_ADMISSIBLE}}
\begin{algorithmic}[1]
\Require{Hierarchical mesh $\mathcal{Q}$, marked elements $\mathcal{M} = \cup_{\ell=0}^{N-1} \mathcal{M}_\ell$, admissibility class $\mu$} 
\Ensure{Refined hierarchical mesh $\mathcal{Q}$}
\For {$\ell = N-1, \ldots, 0$}
\For {$Q \in \mathcal{M}_\ell$ and $Q$ adjacent to a vertex}
\State{$\mathcal{M}_\ell \gets \mathcal{M}_\ell \cup \mathcal{N}_\chi(Q)$}
\EndFor
\For {$Q \in \mathcal{M}_\ell$}
\State{$\mathcal{M}_{\ell-\mu+1} \gets \mathcal{M}_{\ell-\mu+1} \cup \mathcal{N}_r(Q,\mu)$}
\EndFor
\State{Update $\mathcal{Q}$ by replacing $\mathcal{M}_\ell$ with their children}
\EndFor
\end{algorithmic}
\end{algorithm}

Finally, we present in Algorithm~\ref{COARSEN_ADMISSIBLE} a new coarsening strategy that extends the one introduced in \cite{THB-refinement-coarsening} to the case of $C^1$ multi-patch hierarchical splines, in the sense that it not only maintains the grading of the mesh, but also ensures the linear independence of the multi-patch basis functions. The first part of the algorithm, from lines 3 to 9, is essentially the same as the one presented in \cite{THB-refinement-coarsening}. We are given a hierarchical mesh $\mathcal{Q}$, and a list of active elements $\mathcal{M}$ that have been marked to be coarsened, from which we compute a list of coarser elements $\mathcal{M}_c$ to be reactivated. As in the mentioned reference, we choose a conservative strategy in which coarse elements are reactivated only if all their children are marked, and this is ensured by the first condition in line 6. The second condition of the same line checks whether the element can be reactivated without violating the admissibility condition. An element $Q$ of level $\ell-1$ can be reactivated only if all the elements of level $\ell+\mu - 1$ in a certain neighborhood, that we denote by $\mathcal{N}_c(Q,\mu)$, are not active. The definition of this neighborhood is based on the support of the multi-patch $C^1$ splines of different levels, and is completely analogous to the standard case of THB-splines in a single patch.

\begin{figure}[!t]
\begin{subfigure}[Marking for coarsening (light grey, grey and dark grey elements).]{
\includegraphics[trim=0cm 0cm 0cm 0cm, clip, width=0.46\textwidth]{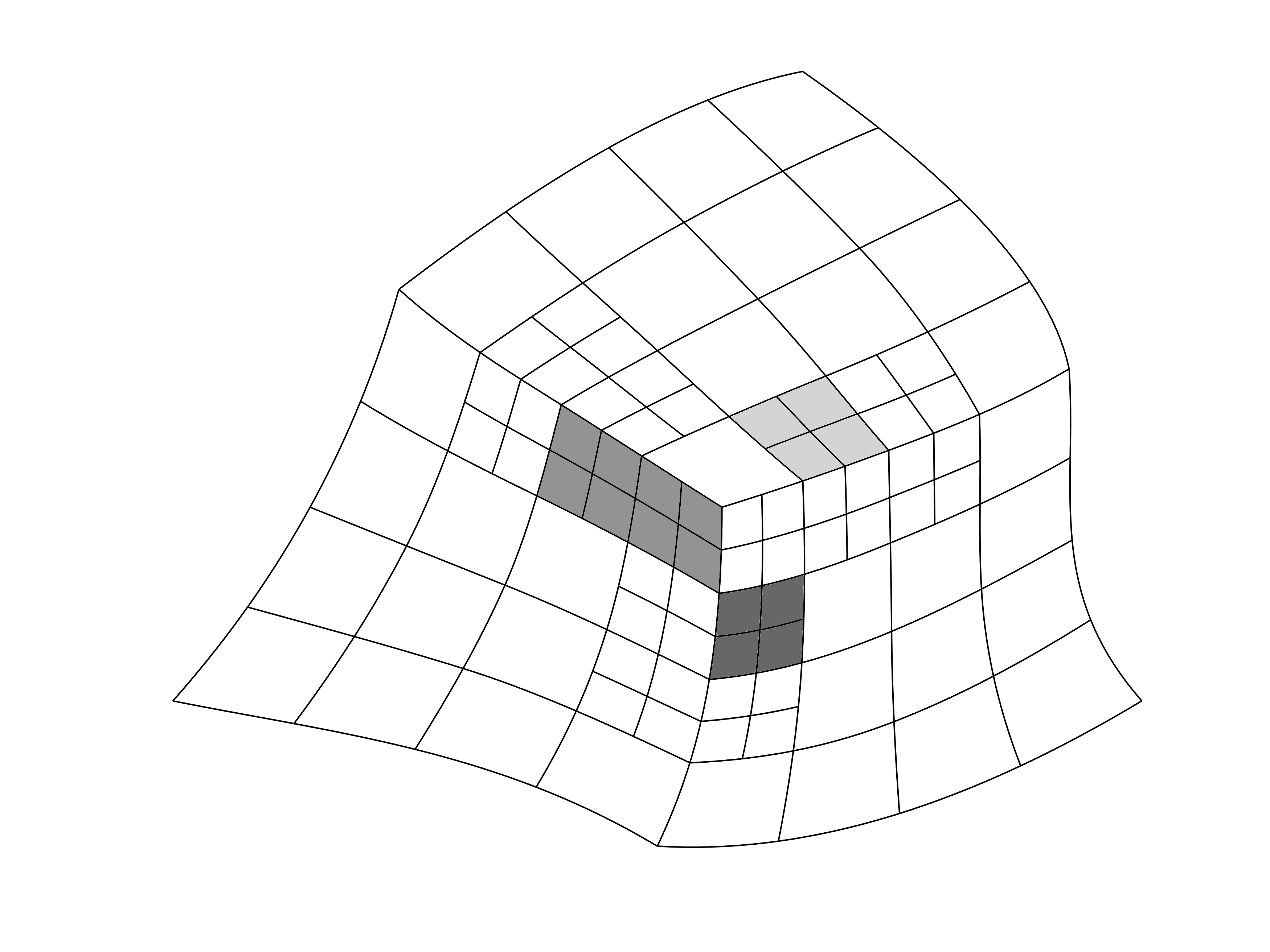}
}
\end{subfigure}
\begin{subfigure}[Resulting mesh after coarsening.]{
\includegraphics[trim=0cm 0cm 0cm 0cm, clip, width=0.46\textwidth]{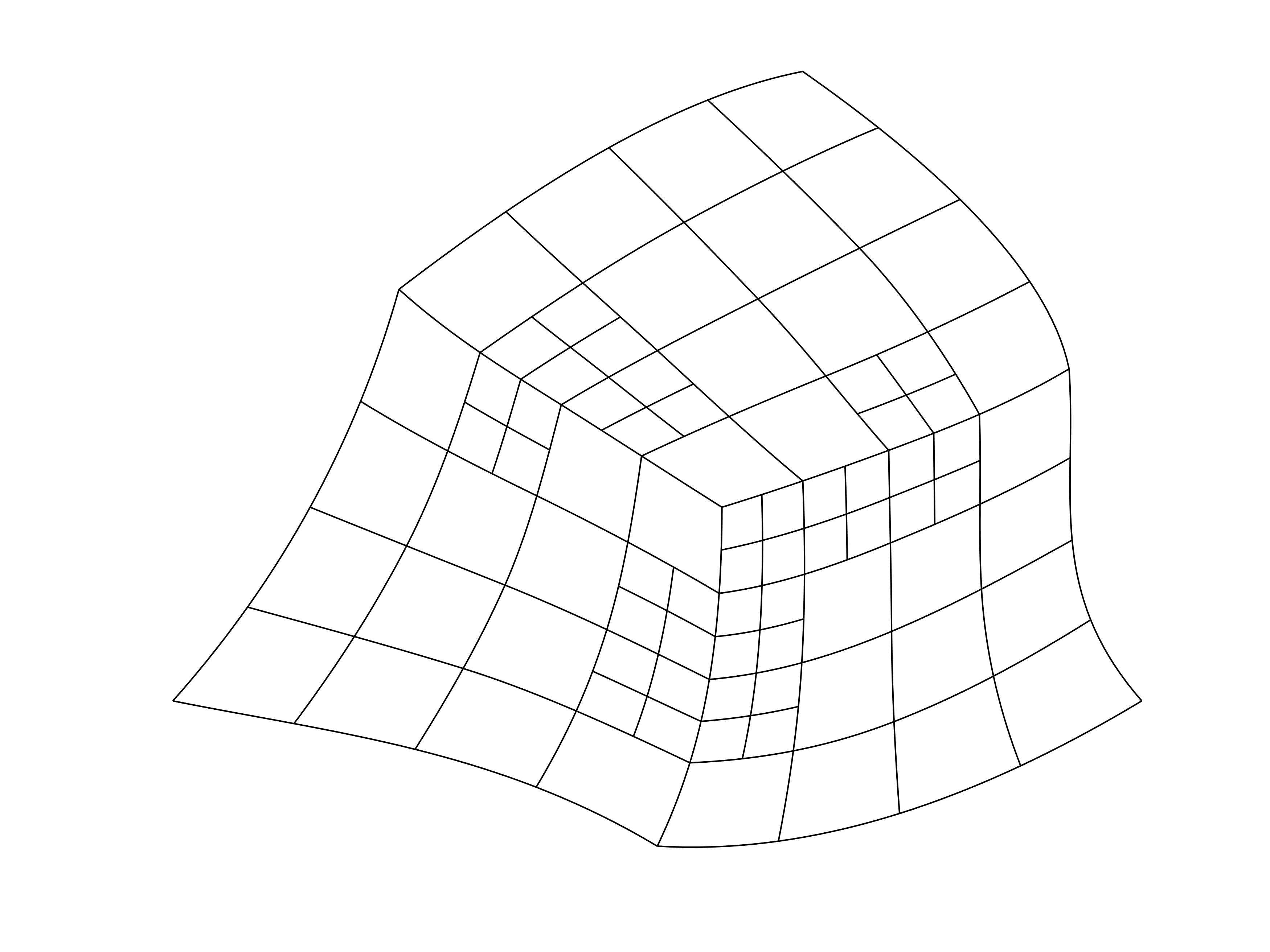}
}
\end{subfigure}
\caption{Example of marked elements in vertex patches (a) and of the corresponding mesh obtained after coarsening (b). An element (whose children are all marked) is reactivated if the elements adjacent to the vertices in the same patch are coarser (light grey) or are also marked for reactivation (grey). If neither of these conditions is satisfied, the element is not reactivated (dark grey).}
\label{fig:coarsening}
\end{figure}

The last steps of the algorithm, from lines 10 to 14, check that the condition for linear independence is not violated. In the refinement step, whenever an element adjacent to a vertex was marked, the vertex-patch neighborhood $\mathcal{N}_\chi(Q)$ was also marked. In the coarsening step we have to do an inverse check: before reactivating an element in the vertex-patch neighborhood, we have to make sure that the element adjacent to the vertex has not been refined yet or it will be reactivated. To do so we define, for an element $Q$ of level $\ell$, the neighborhood $\mathcal{N}_{\chi,c}(Q)$ as the inverse of the vertex-patch neighborhood that we used for refinement. The set is defined as
\[
\mathcal{N}_{\chi,c}(Q) = \{ Q' : Q \in \mathcal{N}_{\chi}(Q') \},
\]
and it collects elements $Q'$ of level $\ell$, and adjacent to a vertex, such that $Q$ belongs to their vertex-patch neighborhood. 
Then, the condition in line 11 checks that the elements marked at that point are either not refined yet or will be reactivated, otherwise the element $Q$ must remain refined. Examples of the effect of this check can be seen in Figure \ref{fig:coarsening}: an element $Q$, whose children are all marked, can be reactivated if the elements in $\mathcal{N}_{\chi,c}(Q)$ are coarser (light grey) or are also marked (grey). If neither of these conditions is satisfied, the element is not reactivated (dark grey). As before, the order of the two steps is important, and the admissibility check must be done before the check for linear independence close to the vertices. 
\begin{algorithm}[h]
\caption{\texttt{Coarsen} \label{COARSEN_ADMISSIBLE}}
\begin{algorithmic}[1]
\Require{Hierarchical mesh $\mathcal{Q}$, marked active elements $\mathcal{M} = \cup_{\ell=1}^{N-1} \mathcal{M}_\ell$, admissibility class $\mu$} 
\Ensure{Coarsened hierarchical mesh $\mathcal{Q}$}
\For {$\ell = N-1, \ldots, 1$}
\State{$\mathcal{M}_c = \emptyset$}
\State{$\mathcal{R}_{c} \gets \texttt{get\_parents}(\mathcal{M}_\ell)$}
\For {$Q \in \mathcal{R}_{c}$}
\State{$Q_c \gets \texttt{get\_children}(Q)$}
\If {$Q_c \subset \mathcal{M}_\ell$ \textbf{and} $ \mathcal{N}_c(Q,\mu) = \emptyset$}
\State{$\mathcal{M}_c \gets \mathcal{M}_c \cup Q$}
\EndIf
\EndFor
\For{$Q \in \mathcal{M}_c$}
\If $\mathcal{N}_{\chi,c}(Q) \not \subset (\mathcal{M}_c\cup\mathcal{A}_{\ell-1})$
\State{$ \mathcal{M}_c \gets \mathcal{M}_c \setminus Q$}
\EndIf
\EndFor
\State{Update $\mathcal{Q}$ by activating $\mathcal{M}_c$ and removing their children}
\EndFor
\end{algorithmic}

\end{algorithm}

\begin{remark}
In the refinement algorithm we have implicitly assumed that the coarsest mesh is at least $4 \times 4$ on each patch. For coarser meshes, the check in line 2 of Algorithm~\ref{REFINE_ADMISSIBLE} has to be repeated once, because the neighborhood of an element adjacent to a vertex can contain elements adjacent to other vertices. Similarly, for coarse meshes an element can be close to two vertices, and the neighborhood $\mathcal{N}_{\chi,c}(Q)$ may contain two elements, for which the check in line 11 of Algorithm~\ref{COARSEN_ADMISSIBLE} has to be done for each of them. Due to this condition, if the mesh of one patch is formed by $4 \times 4$ elements, coarsening will only occur if all the elements of the patch are marked.
\end{remark}

%
%
%
%
%
%
%
%
\section{Adaptive isogeometric algorithm for phase-field modeling of the Cahn-Hilliard equation} \label{sec:aigm-ch}
The adaptive isogeometric methods introduced in Section \ref{sec:ada_mp} are here specialized to solve the Cahn-Hilliard equation, presented in Section \ref{sec:phase_field_discretization}.  
Analyses of phase separation are sensitive to the discretization of the initial condition, requiring a fine mesh at the initial steps, that can then be coarsened and adapted to track the phase interface during the simulations. We implement this rationale according to the algorithms herein described.

\begin{algorithm}[th]
\caption{\texttt{Solve\_Cahn\_Hilliard\_equation} \label{alg:solve_CH}}
\begin{algorithmic}[1]
\Require{Coarsest level mesh $\mathcal{Q}_0$, admissibility class $\mu$, maximum number of hierarchical levels $N$, and initial condition $ u_0 \left( \mathbf{x} \right)$ } 
\Ensure{Control variables $\mathbf{U}:= \{ \hat{\mathbf{u}}_0 , \dots ,  \hat{\mathbf{u}}_{T} \}$, and hierarchical meshes $\mathbf{Q}:= \{ \mathcal{Q}_0 , \dots , \mathcal{Q}_{T} \}$ at every time step}
\State{$t \gets 0$}
\State{$ \mathcal{Q}_t \gets G^N $} \Comment{Initialize as a uniform mesh of level $N$}
\State{$\left( \hat{\mathbf{u}}_t, \dot{\hat{\mathbf{u}}}_t \right) \gets L^2$-projection $\left( \mathcal{Q}_t, u_0 \left(\mathbf{x} \right) \right)$}  \Comment{Compute initial control variables}
\While{$t<T$}
\State{$\left( \mathcal{Q}_{t^*},  \hat{\mathbf{u}}_{t^*}, \dot{\hat{\mathbf{u}}}_{t^*}, \mathcal{E}_{\mathcal{Q}} \right) \gets $ \texttt{advance\_C\_H\_adaptive} $\left( \mathcal{Q}_{t},  \hat{\mathbf{u}}_{t}, \dot{\hat{\mathbf{u}}}_{t}, \mu, N \right)$} \Comment{Algorithm \ref{alg:CH_step_adaptive}}
\State{$t \gets t +  \Delta t$}
\State{$\left( \mathcal{Q}_{t},  \hat{\mathbf{u}}_{t}, \dot{\hat{\mathbf{u}}}_{t} \right) \gets $ \texttt{coarsening\_algorithm} $\left( \mathcal{Q}_{t^*},  \hat{\mathbf{u}}_{t^*}, \dot{\hat{\mathbf{u}}}_{t^*}, \mu, \mathcal{E}_{\mathcal{Q}} \right)$}  \Comment{Algorithm \ref{alg:CH_coarsening}}
\EndWhile
\end{algorithmic}
\end{algorithm}

In a simulation, as shown in Algorithm \ref{alg:solve_CH}, we initialize the mesh with a predefined maximum hierarchical level $N$, such that the element size $h$ is compatible with the characteristic length scale of the problem, presented in \cite{gomez2008isogeometric}. Differently from other schemes for adaptivity, the maximum hierarchical depth used in the refinement area is fixed a priori and does not depend on error estimates performed during the analysis \cite{geopdes-hierarchical, THB-refinement-coarsening}. 

The initial conditions are then projected onto the finest mesh by means of a classical $L^2$-projection while the velocity can be recovered, for instance, by solving \eqref{eq:residual_CH} providing the initial solution for time integration. Within a single time step, the Cahn-Hilliard equation is initially solved refining the mesh, as shown in Algorithm \ref{alg:CH_step_adaptive}, and then the control variables are projected on a coarse discretization, as presented in Algorithm \ref{alg:CH_coarsening}.

\begin{algorithm}[th]
\caption{\texttt{Advance\_C\_H\_adaptive} \label{alg:CH_step_adaptive}}
\begin{algorithmic}[1]
\Require{$\mathcal{Q}_{t}$,  $\hat{\mathbf{u}}_{t}$, $\dot{\hat{\mathbf{u}}}_{t}$, $\mu$, $N$} 
\Ensure{$ \mathcal{Q}_{t^*} ,  \hat{\mathbf{u}}_{t^*}, \dot{\hat{\mathbf{u}}}_{t^*}, \mathcal{E}_{\mathcal{Q}}$ }
\State{$\left( \mathcal{Q}_{n} ,  \hat{\mathbf{u}}_{n}, \dot{\hat{\mathbf{u}}}_{n} \right) \gets \left(  \mathcal{Q}_{t} ,  \hat{\mathbf{u}}_{t}, \dot{\hat{\mathbf{u}}}_{t}\right) $}
\State{$ \left( \hat{\mathbf{u}}_{n+1}, \dot{\hat{\mathbf{u}}}_{n+1} \right) \gets $ nonlinear generalized-$\alpha$ step $\left( \mathcal{Q}_{n} ,  \hat{\mathbf{u}}_{n}, \dot{\hat{\mathbf{u}}}_{n} \right)$ } \Comment{See Section \ref{sec:phase_field_discretization}}
\State{$ \mathcal{E}_{\mathcal{Q}} \gets$ compute indicator $\left(\mathcal{Q}_{n} ,   \hat{\mathbf{u}}_{n+1} \right)$ } \Comment{See  \eqref{eq:est_field} and \eqref{eq:est_grad}}
\State{$ \mathcal{M} \gets $ mark elements for refinement $\left( \mathcal{E}_{\mathcal{Q}}, \mathcal{Q}_{n} \right)$ }   \Comment{Check maximum refinement level $N$}
\While{$\mathcal{M} \neq  \emptyset$}
\State{$ \mathcal{Q}_{n^*}  \gets $ \texttt{refine} $\left( \mathcal{Q}_n,  \mathcal{M}, \mu  \right)$ } \Comment{Algorithm \ref{REFINE_ADMISSIBLE}}
\State{$\left( \hat{\mathbf{u}}_{n^*}, \dot{\hat{\mathbf{u}}}_{n^*} \right) \gets $ represent on a refined space $\left( \mathcal{Q}_{n^*},  \mathcal{Q}_n,  \hat{\mathbf{u}}_{n}, \dot{\hat{\mathbf{u}}}_{n} \right) $  } \Comment{Described in \cite{geopdes-hierarchical}}
\State{$\left(  \mathcal{Q}_{n},  \hat{\mathbf{u}}_{n}, \dot{\hat{\mathbf{u}}}_{n}  \right) \gets \left(  \mathcal{Q}_{n^*},  \hat{\mathbf{u}}_{n^*}, \dot{\hat{\mathbf{u}}}_{n^*} \right)$}
\State{$ \left( \hat{\mathbf{u}}_{n+1}, \dot{\hat{\mathbf{u}}}_{n+1} \right) \gets $ nonlinear generalized-$\alpha$ step $\left( \mathcal{Q}_{n} ,  \hat{\mathbf{u}}_{n}, \dot{\hat{\mathbf{u}}}_{n} \right)$ } \Comment{See Section \ref{sec:phase_field_discretization}}
\State{$ \mathcal{E}_{\mathcal{Q}} \gets$ compute indicator $\left(\mathcal{Q}_{n} ,   \hat{\mathbf{u}}_{n+1} \right)$ } \Comment{See  \eqref{eq:est_field} and \eqref{eq:est_grad}}
\State{$ \mathcal{M} \gets $ mark elements for refinement $\left( \mathcal{E}_{\mathcal{Q}}, \mathcal{Q}_{n} \right)$ }   \Comment{Check maximum refinement level $N$}
\EndWhile
\State{$\left( \mathcal{Q}_{t^*},  \hat{\mathbf{u}}_{t^*}, \dot{\hat{\mathbf{u}}}_{t^*} \right) \gets \left( \mathcal{Q}_{n}, \hat{\mathbf{u}}_{n+1}, \dot{\hat{\mathbf{u}}}_{n+1} \right)$}
\end{algorithmic}
\end{algorithm}

The Cahn-Hilliard equation is solved adapting the mesh to the solution till the phase interface is discretized with elements of the finest hierarchical level, as shown in Algorithm~\ref{alg:CH_step_adaptive}. It is an iterative process in which the nonlinear Cahn-Hilliard equation is initially solved advancing the solution in time by means of generalized-$\alpha$ method, iterating up to convergence of the residual (see Section~\ref{sec:phase_field_discretization}). Such a solution is successively used to compute a posteriori an indicator of the mesh quality $\mathcal{E}_{\mathcal{Q}}$ for which alternative definitions of the indicator are provided and compared in Section \ref{sec:num_estim}. Afterwards, the elements of the mesh are marked for refinement verifying that the maximum hierarchical depth is not exceeded. If elements are marked, then the mesh is updated refining such elements (see Algorithm~\ref{REFINE_ADMISSIBLE}), and the initial conditions of the analyzed time step are represented in the finer space to restart the process with the new mesh.

\begin{algorithm}[th]
\caption{\texttt{Coarsening\_algorithm} \label{alg:CH_coarsening}}
\begin{algorithmic}[1]
\Require{ $ \mathcal{Q}_{t^*} , \hat{\mathbf{u}}_{t^*}, \dot{\hat{\mathbf{u}}}_{t^*},   \mu, \mathcal{E}_{\mathcal{Q}}$} 
\Ensure{$ \mathcal{Q}_{t} ,  \hat{\mathbf{u}}_{t}, \dot{\hat{\mathbf{u}}}_{t} $ }
\State{$ \mathcal{M} \gets $ mark elements for coarsening $\left( \mathcal{Q}_{t^*}, \mathcal{E}_{\mathcal{Q}} \right) $}
\State{$  \mathcal{Q}_{t} \gets$ \texttt{coarsen} $\left( \mathcal{Q}_{t^*}, \mathcal{M} , \mu \right) $} \Comment{Algorithm \ref{COARSEN_ADMISSIBLE}}
\State{$\left(  \hat{\mathbf{u}}_{t}, \dot{\hat{\mathbf{u}}}_{t}  \right)  \gets $ penalized $L^2$-projection $\left(  \mathcal{Q}_{t} ,   \mathcal{Q}_{t^*} , \hat{\mathbf{u}}_{t^*}, \dot{\hat{\mathbf{u}}}_{t^*}  \right)$ }
\end{algorithmic}
\end{algorithm}

Once the control variables are advanced in time on a refined mesh, elements away from the transition region are coarsened, as presented in Algorithm \ref{alg:CH_coarsening}. Using the mesh indicator $\mathcal{E}_{\mathcal{Q}}$, several elements may be marked for coarsening and the mesh updated as presented in Algorithm \ref{COARSEN_ADMISSIBLE}. In this process, the solution is approximated on a coarser space by means of the $L^2$-projection. However, to guarantee that the flux through the boundary is null even after projection, we add a penalty term on the boundary. For instance, the projection of the control variables ($\hat{\mathbf{u}}_{t^*}$) from the refined space (basis functions $\mathbf{N}_f$) onto the coarse set of basis functions ($\mathbf{N}_c$) results in
\begin{equation*}
\left[ 
\int_{\bar{\Omega}} \mathbf{N}_c^\top  \mathbf{N}_c \, \hbox{\rm d} \mathbf{x} +
 \int_{\partial \bar{\Omega}}  \left(\nabla \mathbf{N}_c \cdot \mathbf{n} \right)^\top \varepsilon_P \, h \, \left(\nabla \mathbf{N}_c \cdot \mathbf{n} \right) \hbox{\rm d} \mathbf{x} 
\right] 
\hat{\mathbf{u}}_{t} = 
\int_{\bar{\Omega}} \mathbf{N}_c^\top \mathbf{N}_f \hat{\mathbf{u}}_{t^*} \hbox{\rm d} \mathbf{x},
\end{equation*}
where $\varepsilon_P$ is a penalty constant and $\hat{\mathbf{u}}_{t}$ is the coarser set of control variables.

As a final remark, we highlight that if the mesh is preserved between time steps, all the discrete terms that depend linearly on the solution or its time derivative (e.g., the mass matrix $\mathbf{M}$) are preserved too, reducing the overall computational effort.

%
%

%
%
%
%
%
%
\section{Numerical examples}
\label{sec:num_exe}
In this section, we study in 2D the proposed adaptive refinement and coarsening strategy for the Cahn-Hilliard equation by means of some numerical tests implemented in GeoPDEs \cite{geopdesv3, refinement, geopdes-hierarchical}. Firstly, in Section~\ref{sec:num_estim}, we compare two classical indicators adopted in the phase-field modeling literature, whereas in Section~\ref{sec:num_adm} we demonstrate the effectiveness of the grading strategy in terms of accuracy per degree of freedom. Exploiting these results, we then perform several simulations of standard phase separation modes using adaptive meshes in single-patch and multi-patch geometries in Section~\ref{sec:num_nuc_spin} and Section~\ref{sec:num_mpC1}, respectively.

The phase separation mode primarily depends on the average value of the field at the beginning of the analysis \cite{gomez2008isogeometric}; for instance, if the two phases are equally mixed, a spinodal decomposition is expected, while an unbalanced phase concentration leads to nucleation. In our tests, we model these phenomena assuming $\sigma = \nu = 1$ in the double-well function \eqref{dw} and imposing two different values of the initial average field $\bar{u}$, namely, $\bar{u}=0$  and $\bar{u}=0.4$, for spinodal decomposition and nucleation, respectively. Moreover, the initial field is perturbed by a small variation  $\delta\left(\mathbf{x}\right)$ randomly distributed in space, with $\delta \in \left[-0.005, 0.005\right]$, that results in the following initial condition:
\begin{equation}
    u_0\left(\mathbf{x}\right) = \bar{u} + \delta \left(\mathbf{x}\right), 
\end{equation}
interpolated on the basis functions by means of a standard $L^2$-projection, while the initial velocity is assumed to be null.

For both separation modes, we adopt the same set of boundary conditions presented in \cite{gomez2014accurate} and discretized as shown in Section \ref{sec:phase_field_discretization}. Essential boundary conditions are enforced by means of Nitsche's method \cite{zhao2017variational}, correlating the penalty constant $\varepsilon_N$ to the phase-field parameter $\lambda$ as $\varepsilon_N = 10^4 \lambda$, while natural boundary conditions are homogeneous. Moreover, during the projection of the results on a coarse mesh (see Algorithm \ref{alg:CH_coarsening}), the no-flux on the boundary is enforced by means of a penalty method with parameter $\varepsilon_P =10^3$.

In the following numerical examples, we use THB-splines as the base technology for single-patch adaptive analysis and their extension defined in Section \ref{sec:ada_mp} in the multi-patch case, with $\mu$ equal to the degree of the basis functions $p$ to enforce mesh admissibility.
Further, to assess the effect of the proposed adaptive scheme, comparisons are made according to the following definition of error at a given instant $\bar{t}$: 
\begin{equation}
err = \sqrt{\dfrac{\int_{\bar{\Omega}} \left( u \left( \mathbf{x}, \bar{t} \right) - u_{ref} \left( \mathbf{x}, \bar{t} \right)  \right)^2 \, \hbox{\rm d}\mathbf{x}} {\int_{\bar{\Omega}} \left( u_{ref} \left( \mathbf{x}, \bar{t} \right)   \right)^2 \, \hbox{\rm d}\mathbf{x}} },
\label{eq:L2_err_def}
\end{equation}
where the reference solution $u_{ref}$ is computed using standard tensor product B-splines. Results of the adaptive algorithm are compared to (i) overkill solutions and (ii) simulations that employ uniform elements of size $h$ equal to the size of the finest element in the hierarchical mesh, 
disentangling the error due to the coarsening from the representation of the solution close to the phase interfaces.
Considerations on the number of degrees of freedom (dofs) employed in the analysis are derived as well.

In performing our simulations, we set the absolute and the relative (with respect to the norm of the prediction step) tolerances on the norm of the residual for the convergence criterion of the Newton-Raphson method equal to $10^{-10}$, with one of these conditions always achieved in a maximum of five iterations in the tests herein presented. Eventually, we note that also the mesh is adapted within a maximum of four iterations at every time step.

\subsection{Comparison between classical indicators for phase-field modeling}
\label{sec:num_estim}
In adaptive phase-field modeling, classical indicators to select the areas where the mesh is refined or coarsened are based either on the value of the phase field \cite{nagaraja2019phase,proserpio2020framework} or of the phase-field gradient norm \cite{hennig2018}. Herein, we compare them specifically in the context of phase separations. 
The field-based adaptive strategy refines an element $\bar{\Omega}_e$ of the mesh, of area $\left|\bar{\Omega}_e\right|$, if it is in the transition zone between the two phases. For a formulation characterized by binodal points at $u=\pm 1$, this can be implemented according to the following element-wise definition for the estimation and marking strategy:
\begin{equation}
\left[ \mathcal{E}^v_{\mathcal{Q}} \right]_e = 1- \left| \dfrac{ \int_{\bar{\Omega}_e}  v  \, \hbox{\rm d}\mathbf{x}  }{  \left|\bar{\Omega}_e\right| }  \right| \quad\longrightarrow\quad
\begin{cases}
     \text{if  $ \left[ \mathcal{E}^v_{\mathcal{Q}} \right]_e > \alpha_q $ } & \text{marked for refinement,}\\
    \text{otherwise} & \text{marked for coarsening,}
 \end{cases}
\label{eq:est_field}
\end{equation}
where the parameter $\alpha_q$ is a constant value in the range $\alpha_q \in \left(0,1\right)$.  The indicator of the whole mesh $\mathcal{E}^v_{\mathcal{Q}}$ is obtained as the collection of the indicators of every element.
Similarly, the gradient-based strategy refines the zones where the phase field varies in space by checking if the value of the gradient norm exceeds a predefined parameter $\beta_q$:
\begin{equation}
\left[ \mathcal{E}^{\nabla}_{\mathcal{Q}} \right]_e = \dfrac{ \int_{\bar{\Omega}_e} \| \nabla v \| \, \hbox{\rm d}\mathbf{x}  }{  \left|\bar{\Omega}_e\right| }  \quad\longrightarrow\quad\
\begin{cases}
    \text{if $ \left[ \mathcal{E}^{\nabla}_{\mathcal{Q}} \right]_e  > \beta_q $} & \text{marked for refinement,} \\
    \text{otherwise} & \text{marked for coarsening,}
 \end{cases}
\label{eq:est_grad}
\end{equation}
 resulting in a similar definition of the whole mesh indicator $\mathcal{E}^\nabla_{\mathcal{Q}}$. In this second approach, the parameter $\beta_q$ is a positive constant, whose selection depends on the sharpness of the field variation, related to the parameter $\lambda$. Herein, we adopt alternatively one of these two definitions of the indicator to adapt the mesh, i.e.,  either $\mathcal{E}_{\mathcal{Q}} = \mathcal{E}^v_{\mathcal{Q}}$ or $\mathcal{E}_{\mathcal{Q}} = \mathcal{E}^\nabla_{\mathcal{Q}}$ in Algorithm~\ref{alg:CH_step_adaptive}.

\begin{figure}[th]
\begin{tabular}{cc}
(a)
\begin{minipage}[b]{.45\textwidth}
\includegraphics[width=\textwidth]{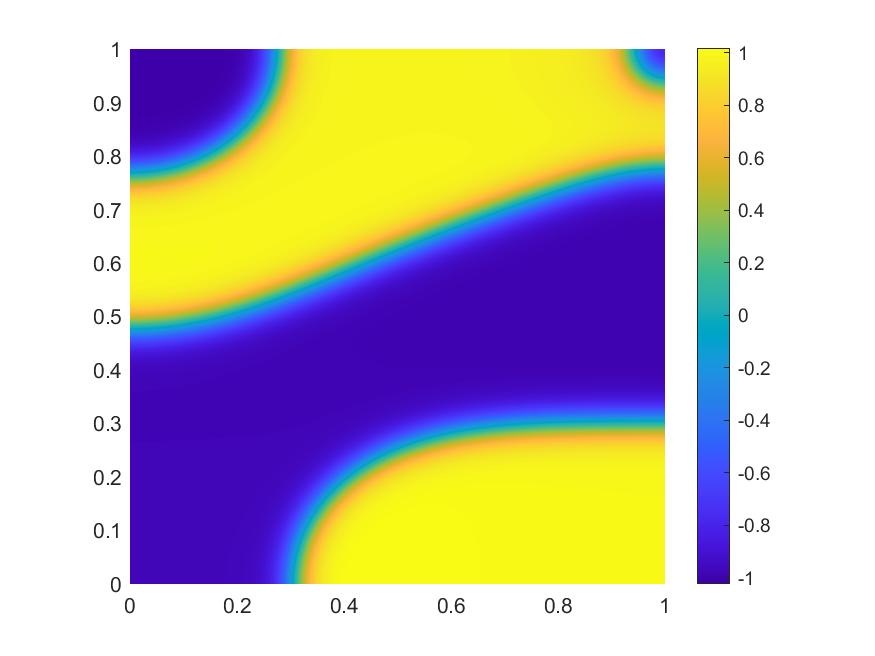} 
\end{minipage}
&
(b)
\begin{minipage}[b]{.45\textwidth}
\includegraphics[width=\textwidth]{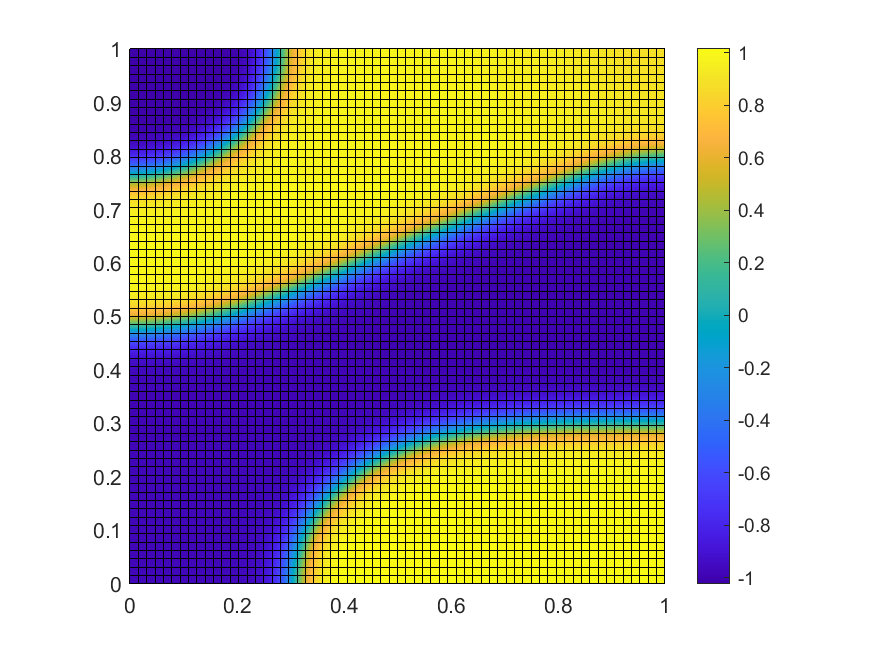} 
\end{minipage}
\end{tabular}
\caption{\label{fig:uniform_spin_square} Solution computed using tensor product meshes at the final simulation time. (a) Overkill solution computed with square elements of size $h=1/256$ (66,564 degrees of freedom (dofs)). (b) Solution computed with a uniform mesh composed of elements of size $h=1/64$ (4,356 dofs).}
\end{figure}

To compare the indicators, we model a spinodal decomposition on a unit square in the time interval $ t \in \left[0, 1 \right]$  discretized with uniform time steps of size $\Delta t  = 10^{-3}$. By assuming $\lambda= 6.15 \times 10^{-4}$, the element size $h$ to correctly simulate the field evolution can be determined as presented in \cite{gomez2008isogeometric} through the relation $h  \approx \sqrt{\lambda/2.5}$, that results in $h=1/64$. We use such a value to define a reference tensor-product mesh and the size of the finest element of the hierarchical mesh: the coarsest representation of the geometry, composed of a single square element, entails seven hierarchical levels ($N=7$) during the analysis. Furthermore, to assess the accuracy, we perform an overkill simulation using a mesh size four times smaller ($h=1/256$), as shown in Figure \ref{fig:uniform_spin_square}, and exactly the same initial condition obtained through the knot insertion algorithm, as done in \cite{gomez2014accurate}. For this single-patch geometry, we use quadratic ($p=2$) $C^1$-continuous B-splines and THB-splines to define tensor product and hierarchical spaces, respectively.

\begin{figure}[H]
\begin{tabular}{cc}
(a)
\begin{minipage}[b]{.45\textwidth}
\includegraphics[width=\textwidth]{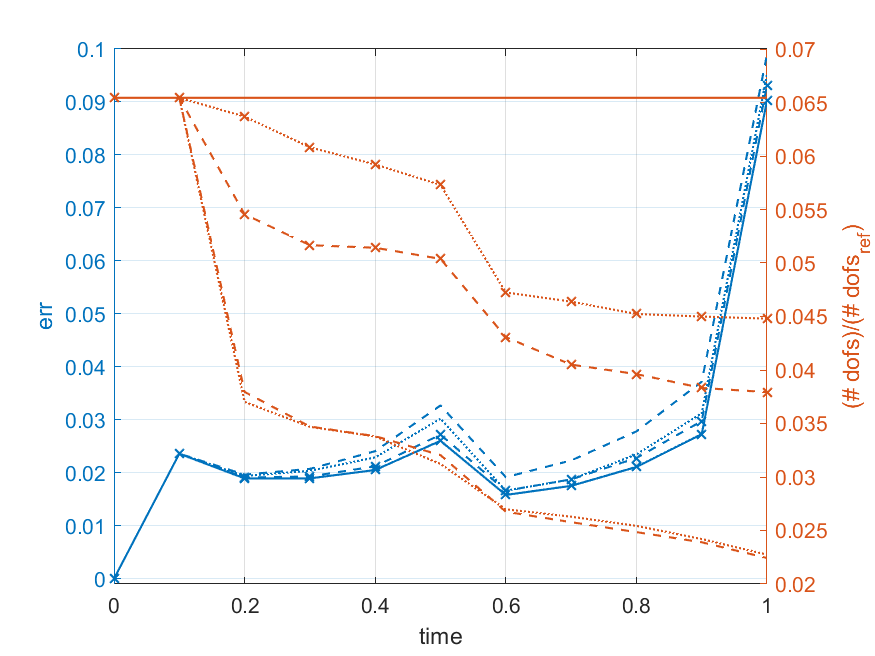} 
\end{minipage}
&
(b)
\begin{minipage}[b]{.45\textwidth}
\includegraphics[width=\textwidth]{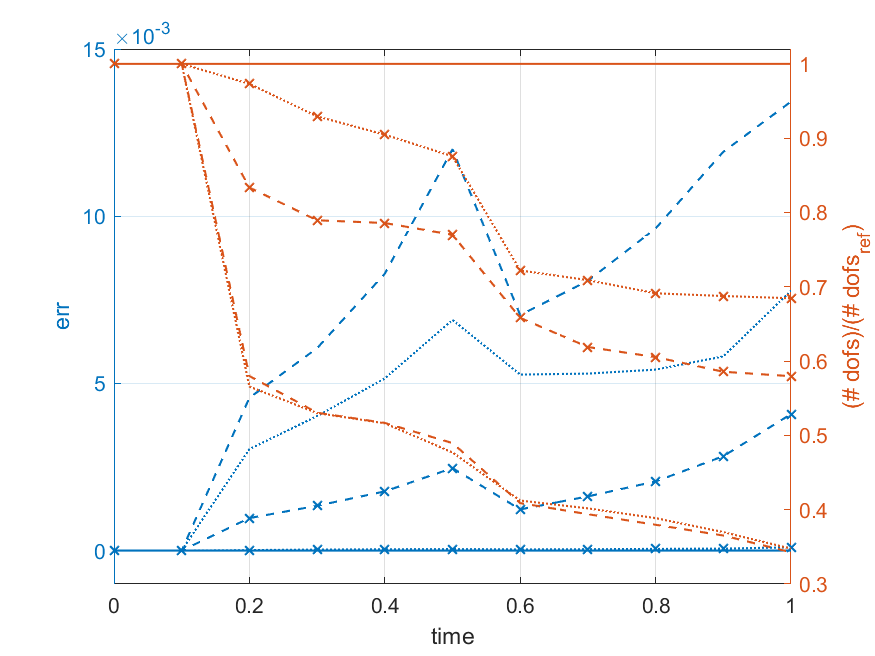} 
\end{minipage}
\\
\begin{minipage}[b]{.45\textwidth}
\centering
\includegraphics[width=\textwidth]{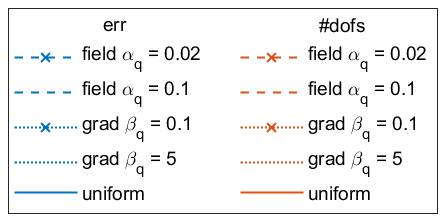}
\vspace{20pt}
\end{minipage}
&
(c)
\begin{minipage}[b]{.45\textwidth}
\centering
\includegraphics[width=\textwidth]{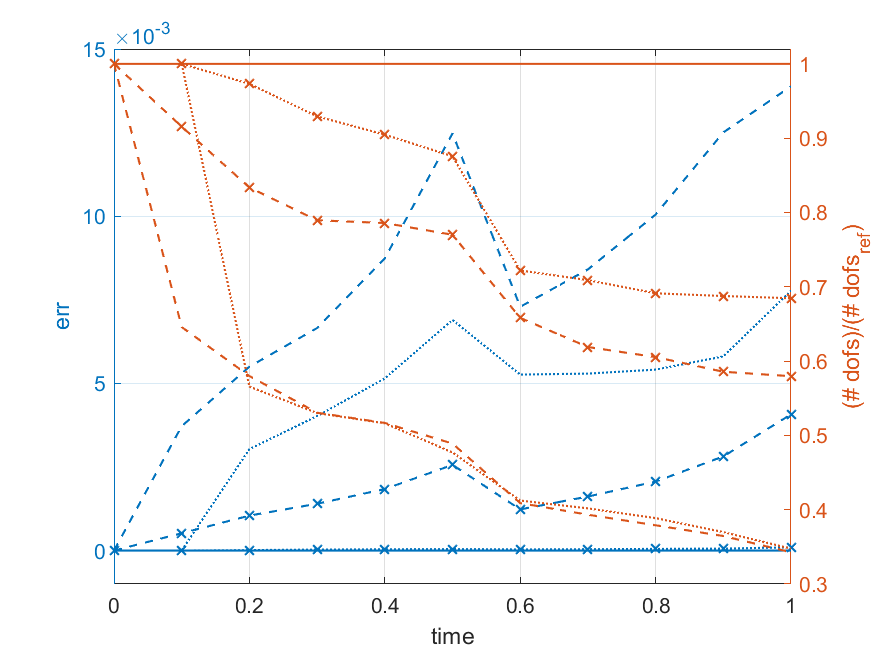}
\end{minipage}
\\
(d)
\begin{minipage}[b]{.45\textwidth}
\includegraphics[width=\textwidth]{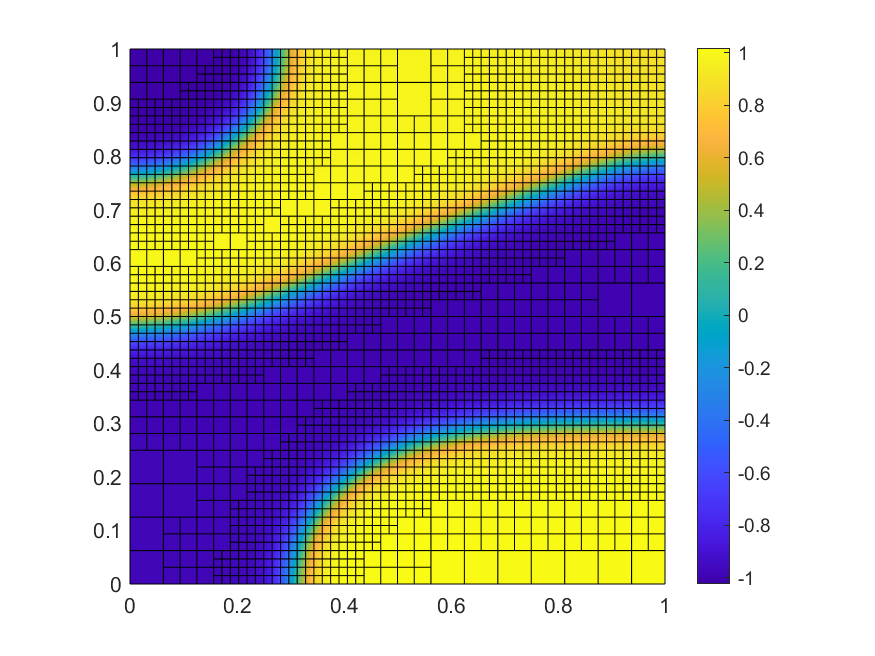} 
\end{minipage}
&
(e)
\begin{minipage}[b]{.45\textwidth}
\includegraphics[width=\textwidth]{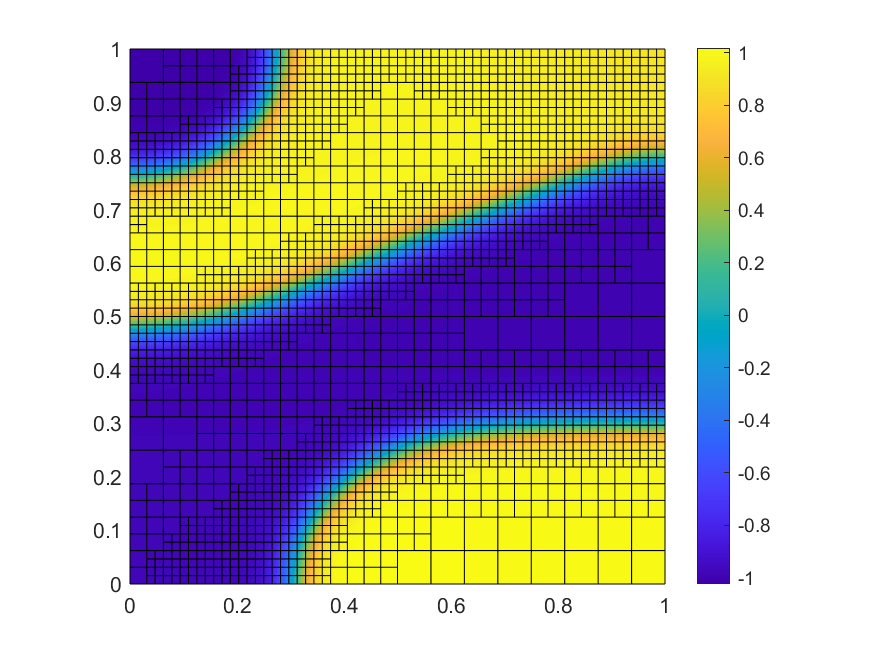} 
\end{minipage}
\end{tabular}
\caption{\label{fig:estim_comp} Comparison between gradient-based and field-based indicators for various time steps of the analysis. (a) Error ($err$, see \eqref{eq:L2_err_def}), in blue, and relative number of degrees of freedom ($\# dofs / \# dofs_{ref}$), in red, with respect to the overkill solution, the adaptivity algorithm is activated at time $t=0.1$. (b) The same quantities are presented using the uniform mesh in Figure \ref{fig:uniform_spin_square}(b) as the reference. (c) Field-based adaptivity activated at the initial time step compared to the uniform mesh, gradient-based adaptivity (activated at time $t=0.1$) is reported for comparison. (d)-(e) Solution and corresponding hierarchical mesh at the final simulation time: (d) gradient-based adaptivity ($\beta_q=0.1$, adaptivity activated at time $t=0.1$), and (e) field-based adaptivity ($\alpha_q=0.02$, adaptivity activated at time $t=0.1$).}
\end{figure}

In Figure \ref{fig:estim_comp}, the results of the analysis performed using the field-based ($\alpha_q = 0.02$ and $\alpha_q = 0.1$) and the gradient-based ($\beta_q = 0.1$ and $\beta_q = 5$) adaptivity, as well as uniform meshes, are presented. 

The gradient-based strategy fails if the coarsening algorithm is activated at the beginning of the analysis, at least for the analyzed values of $\beta_q$. When the phases are still mixed, the indicator marks most of the elements for coarsening, leading to unreliable results. This observation is in agreement with the findings presented in \cite{hennig2018}. For more restrictive values of $\beta_q$ the analysis may be performed successfully, however, a limited number of elements is coarsened during the simulation. To circumvent such a problem, in all the presented results for the gradient indicator, the coarsening algorithm (see Algorithm~\ref{alg:solve_CH}) is activated after a predefined time $t=0.1$, calibrated using the results of the tensor-product simulation to assess when the phases are separated. The automatic detection of such an instant is out of the scope of this paper. Moreover, we note that the existence of such an instant is not even guaranteed: in fact, if Dirichlet boundary conditions ($u_D=\bar{u} + \delta$) are applied for the whole duration of the analysis, in part of the domain the two phases do not separate.

In Figure \ref{fig:estim_comp}(a), the results of uniform and adaptive analyses, with the coarsening algorithm enabled after $t=0.1$, are compared to the overkill solution. All the simulation results are comparable and the error in time is in the order of 2\% for most of the time whereas the number of employed degrees of freedom is substantially reduced by adaptivity. At the final simulation time, the error increases because a coarse resolution of the transition area modifies the temporal evolution of the field (see Figure \ref{fig:uniform_spin_square}(a) and (b)). Figure \ref{fig:estim_comp}(b) presents the same simulations using the uniform mesh with $h=1/64$ as the reference to better highlight the differences with respect to the refinement in the interface area.

The adaptive strategies are comparable (see Figure \ref{fig:estim_comp}(d) and (e)) but the gradient-based method seems to be slightly more effective in terms of accuracy per degree of freedom, for this problem. However, the field-based indicator successfully completes the analysis activating the coarsening algorithm from the very beginning of the analysis, as shown in Figure \ref{fig:estim_comp}(c) where the coarsening algorithm is enabled at $t=0$, without any a-priori knowledge of the field evolution. Consequently, we select this approach for all the remaining simulations.

\subsection{Impact of mesh grading on simulation accuracy}
\label{sec:num_adm}
\begin{figure}[H]
\begin{tabular}{cc}
(a)
\begin{minipage}[b]{.45\textwidth}
\includegraphics[width=\textwidth]{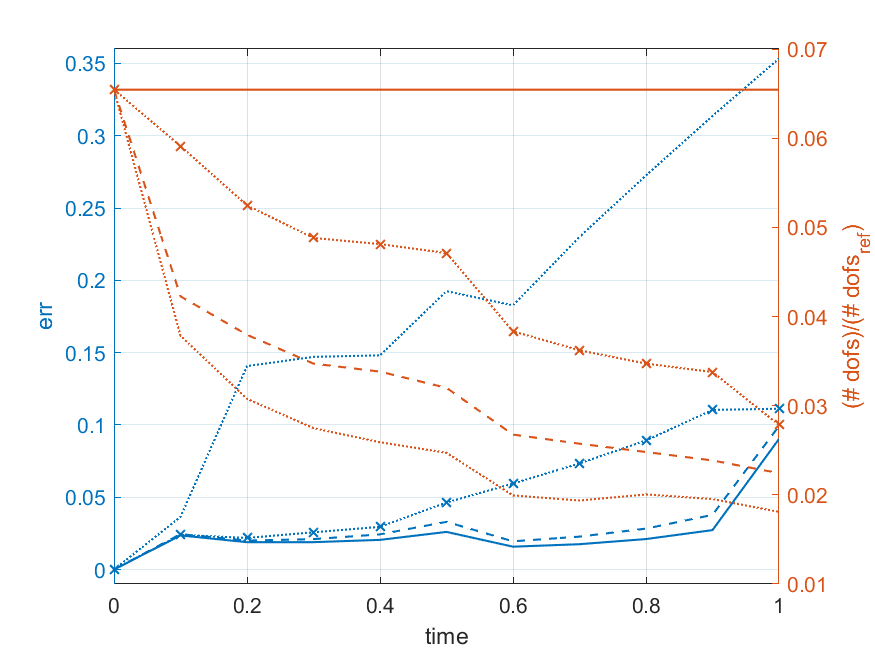} 
\end{minipage}
&
(b)
\begin{minipage}[b]{.45\textwidth}
\includegraphics[width=\textwidth]{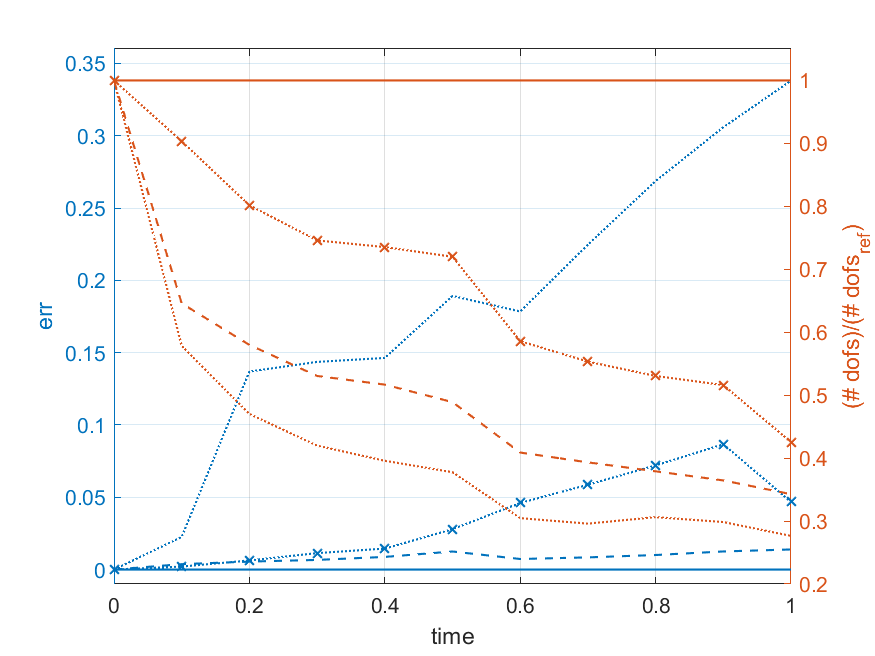} 
\end{minipage}
\\
\begin{minipage}[b]{.45\textwidth}
\centering
\includegraphics[width=\textwidth]{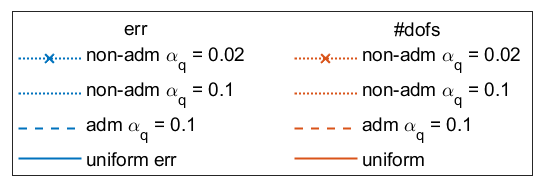}
\vspace{30pt}
\end{minipage}
&
(c)
\begin{minipage}[b]{.45\textwidth}
\centering
\includegraphics[width=\textwidth]{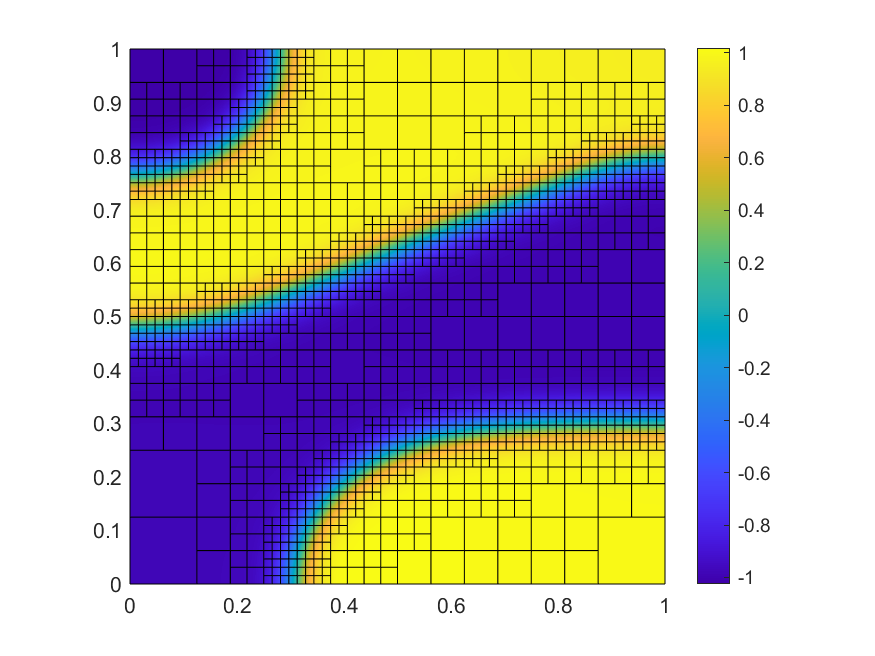}
\end{minipage}
\\
(d)
\begin{minipage}[b]{.45\textwidth}
\includegraphics[width=\textwidth]{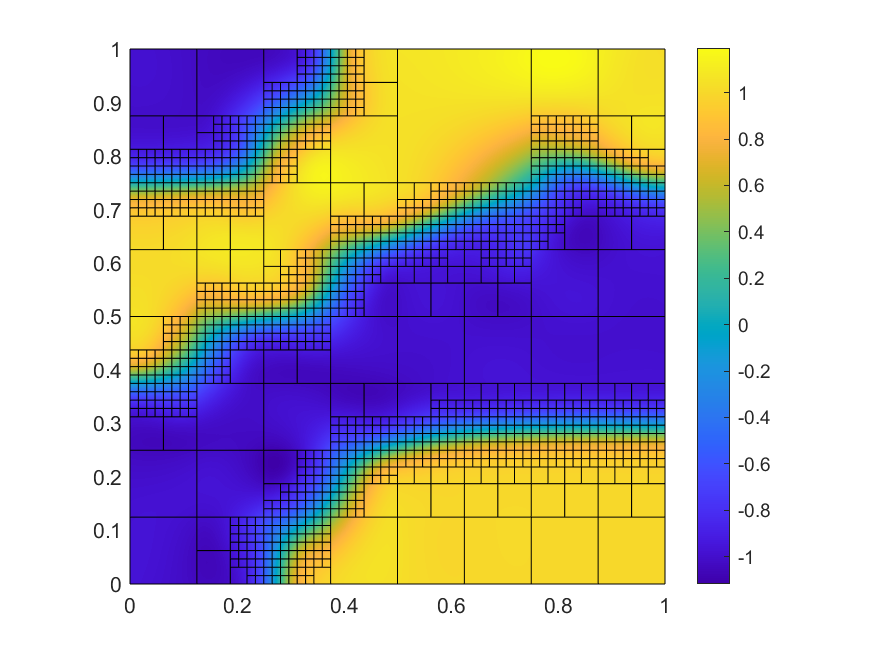} 
\end{minipage}
&
(e)
\begin{minipage}[b]{.45\textwidth}
\includegraphics[width=\textwidth]{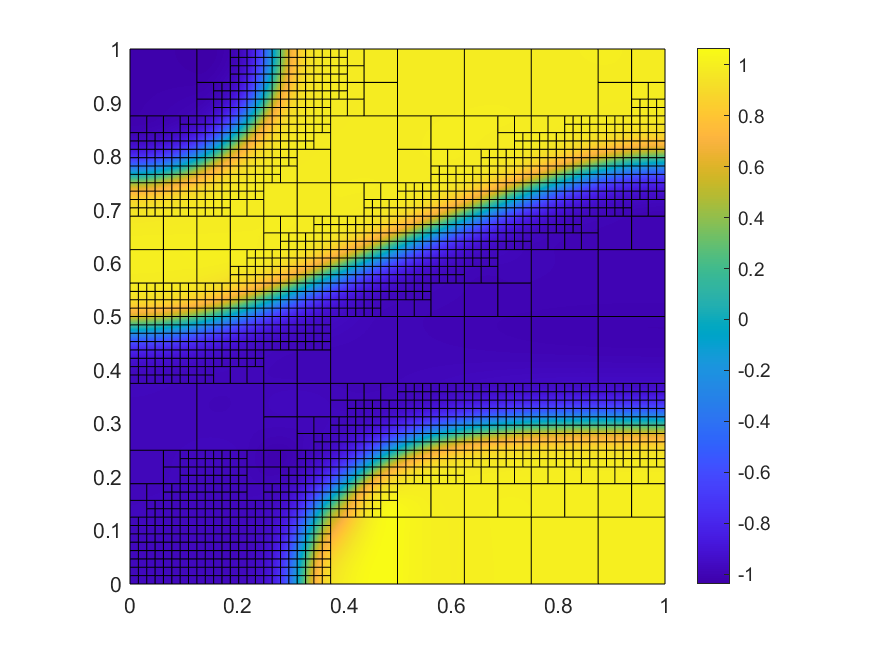} 
\end{minipage}
\end{tabular}
\caption{\label{fig:effect_adm} Comparison between suitably graded and non-admissible meshes for various time steps of the analysis. (a) Error ($err$  see \eqref{eq:L2_err_def}), in blue, and relative number of degrees of freedom ($\# dofs / \# dofs_{ref}$), in red, with respect to the overkill solution Figure \ref{fig:uniform_spin_square}(a). (b) The same quantities are presented using the uniform mesh in Figure \ref{fig:uniform_spin_square}(b) as the reference. (c)-(e) Solution and corresponding hierarchical mesh at the final simulation time computed with adaptive meshes: (c) admissible mesh ($\alpha_q =0.1$), (d) non-admissible mesh ($\alpha_q =0.1$), and (e) non-admissible mesh ($\alpha_q =0.02$).}
\end{figure}

In this second test, we re-simulate the same problem introduced before to assess the effect of a suitable mesh grading on the accuracy per degree of freedom. Starting from the same initial condition, spatial discretization, and time discretization ($\Delta t = 10^{-3}$), we analyze the results obtained by non-admissible meshes ($\alpha_q=0.02$ and $\alpha_q=0.1$) and a graded mesh ($\alpha_q=0.1$).

The results in Figure \ref{fig:effect_adm}(a) show that the admissible hierarchical and the uniform tensor product discretizations are in good agreement while the non-admissible ones are less accurate. These results are better highlighted in  Figure \ref{fig:effect_adm}(b), where the reference solution is the tensor product mesh $h=1/64$: fixing the parameter $\alpha_q=0.1$, the admissible mesh produces a result more accurate than the corresponding non-admissible one using a mildly finer mesh; a comparison between Figure \ref{fig:effect_adm}(c) and (d) visually quantifies such an effect. Similar spurious results for non-admissible meshes were observed in previous studies \cite{THB-refinement-coarsening}. To achieve the same level of accuracy of the graded mesh with a non-admissible discretization, the parameter $\alpha_q$ has to be reduced to enlarge the refined area. By setting $\alpha_q=0.02$ for the non-admissible mesh, a qualitatively correct result is obtained, as shown in Figure \ref{fig:effect_adm}(d). However, the number of degrees of freedom employed is larger than the one used in the admissible discretization while the solution is still less accurate. Consequently, the admissible strategy results in a more accurate solution per degree of freedom showing the benefits of mesh admissibility on the accuracy of the solution, especially when the marking strategy does not provide a graded transition of the element size.

\subsection{Adaptive simulations for nucleation and spinodal decompositions}
\label{sec:num_nuc_spin}
In this section, we study the accuracy of the adaptive scheme with respect to the value of the refinement parameter $\alpha_q$. Initially, we impose $\alpha_q = 0.02, 0.1, \, \rm{and} \,0.2$ for the previous spinodal test. Afterwards, we adopt the same data used before to model a nucleation process, varying only the initial conditions.

\begin{figure}[th]
\begin{tabular}{cc}
(a)
\begin{minipage}[b]{.45\textwidth}
\includegraphics[width=\textwidth]{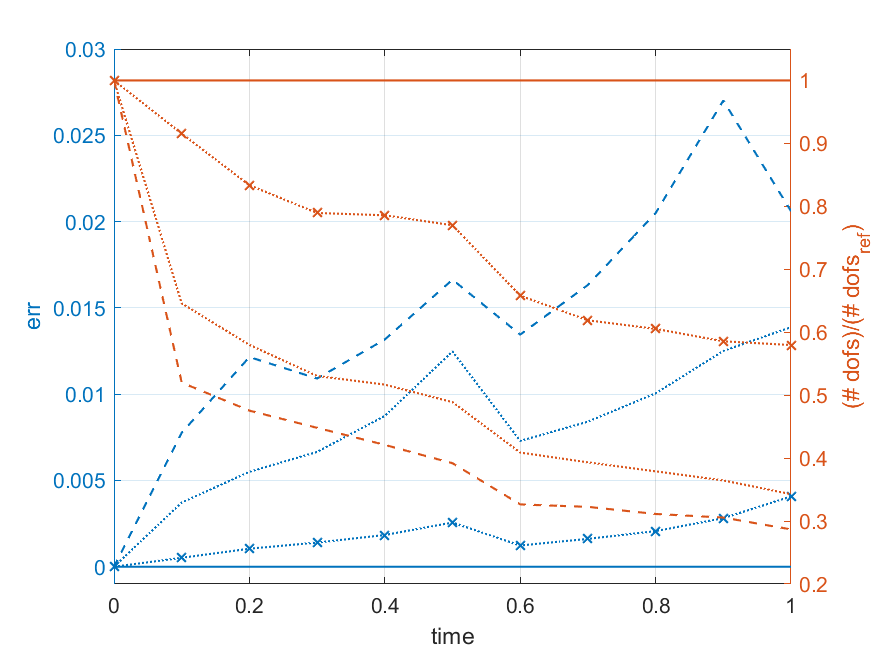} 
\end{minipage}
&
(b)
\begin{minipage}[b]{.45\textwidth}
\includegraphics[width=\textwidth]{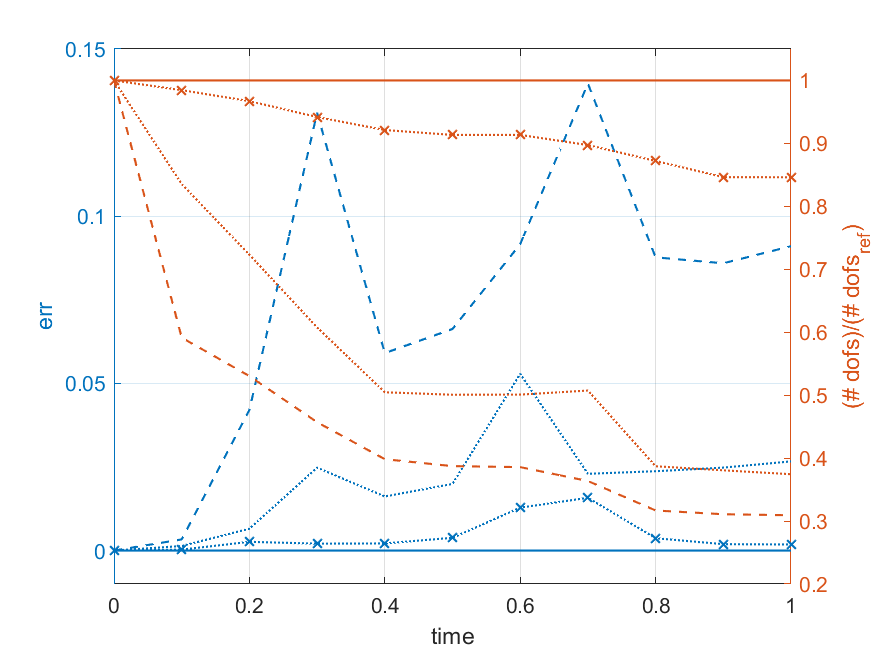} 
\end{minipage}
\\
\begin{minipage}[b]{.45\textwidth}
\centering
\includegraphics[width=\textwidth]{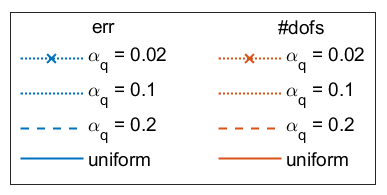}
\vspace{10pt}
\end{minipage}
&
(c)
\begin{minipage}[b]{.45\textwidth}
\centering
\includegraphics[width=\textwidth]{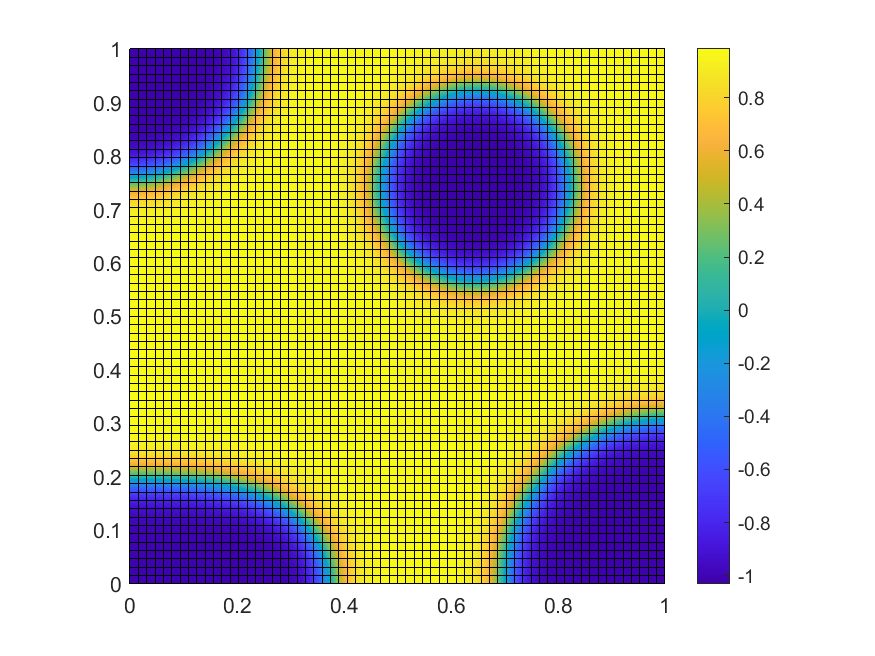}
\end{minipage}
\end{tabular}
\caption{\label{fig:parameter_study} Simulations of phase separations using different values of the parameter $\alpha_q$. (a) Error ($err$  see \eqref{eq:L2_err_def}), in blue, and relative number of degrees of freedom ($\# dofs / \# dofs_{ref}$), in red, for the spinodal decomposition. (b) The same quantities are presented for the nucleation process. (c) Reference solution for nucleation at the final simulation time (the corresponding reference for spinodal decomposition is shown in Figure \ref{fig:uniform_spin_square}(b)).}
\end{figure}

\begin{figure}[t!]
\begin{tabular}{cc}
(a)
\begin{minipage}[b]{.45\textwidth}
\includegraphics[width=\textwidth]{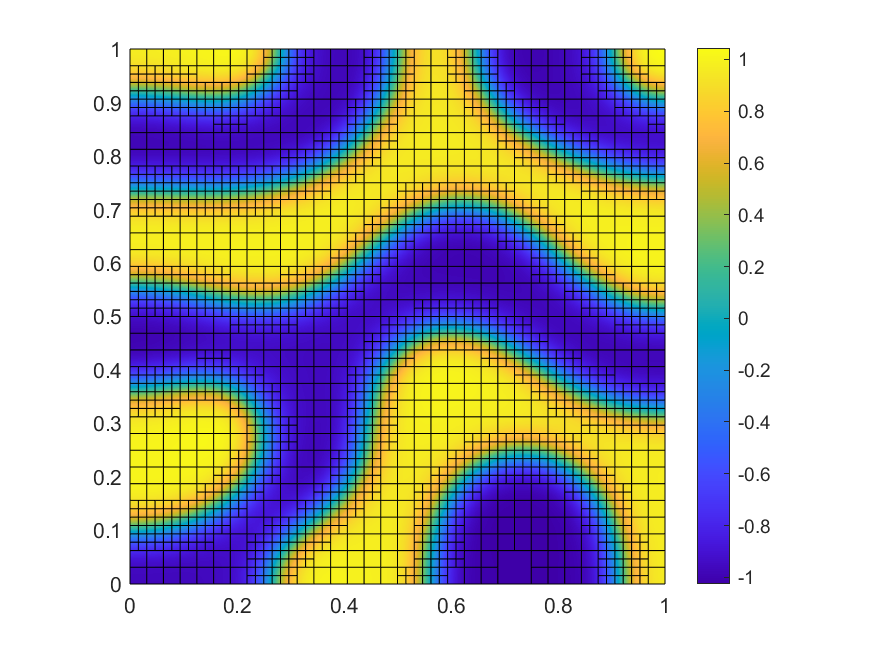} 
\end{minipage}
&
(b)
\begin{minipage}[b]{.45\textwidth}
\includegraphics[width=\textwidth]{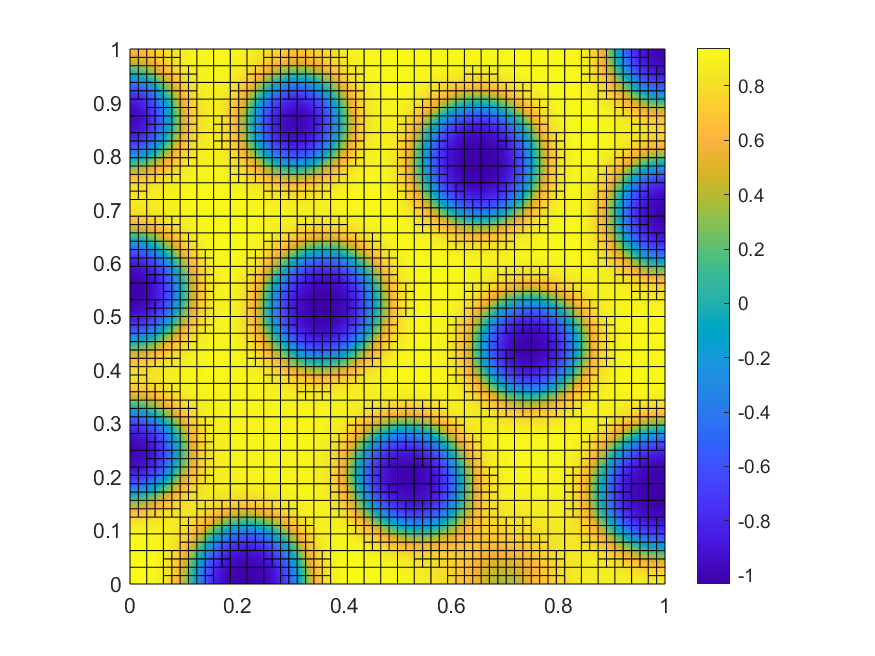} 
\end{minipage}
\\
(c)
\begin{minipage}[b]{.45\textwidth}
\centering
\includegraphics[width=\textwidth]{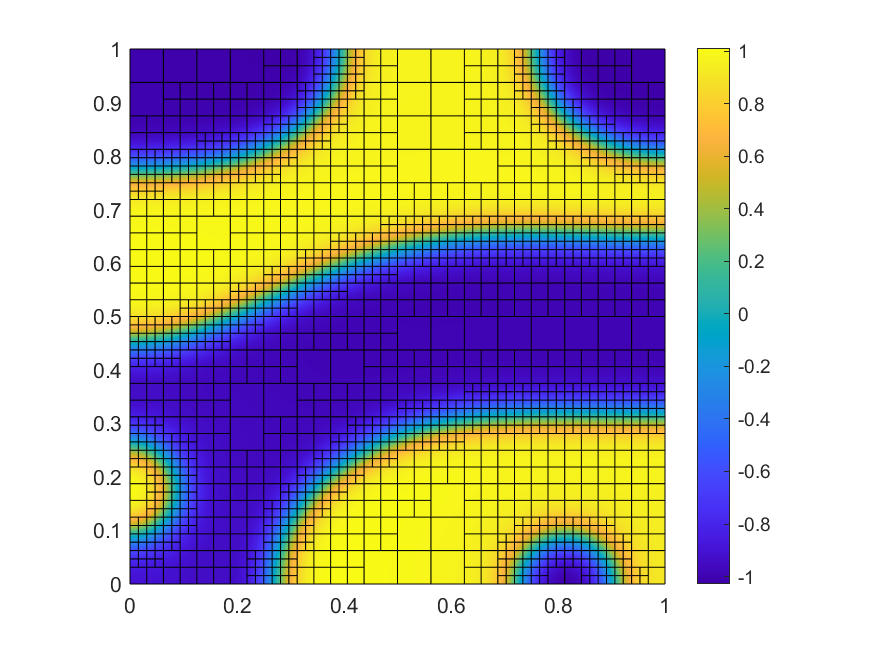} 
\end{minipage}
&
(d)
\begin{minipage}[b]{.45\textwidth}
\centering
\includegraphics[width=\textwidth]{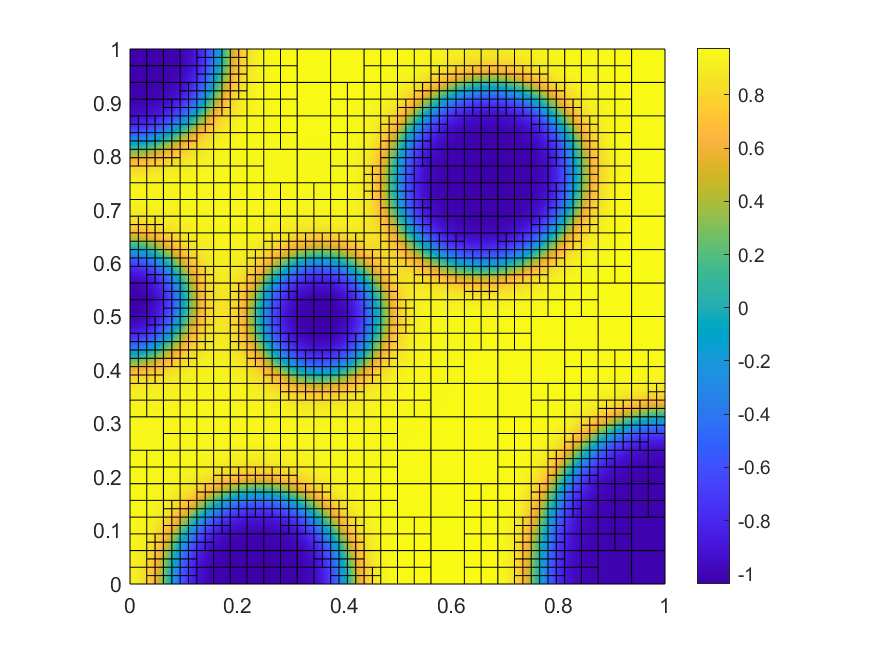} 
\end{minipage}
\\
(e)
\begin{minipage}[b]{.45\textwidth}
\centering
\includegraphics[width=\textwidth]{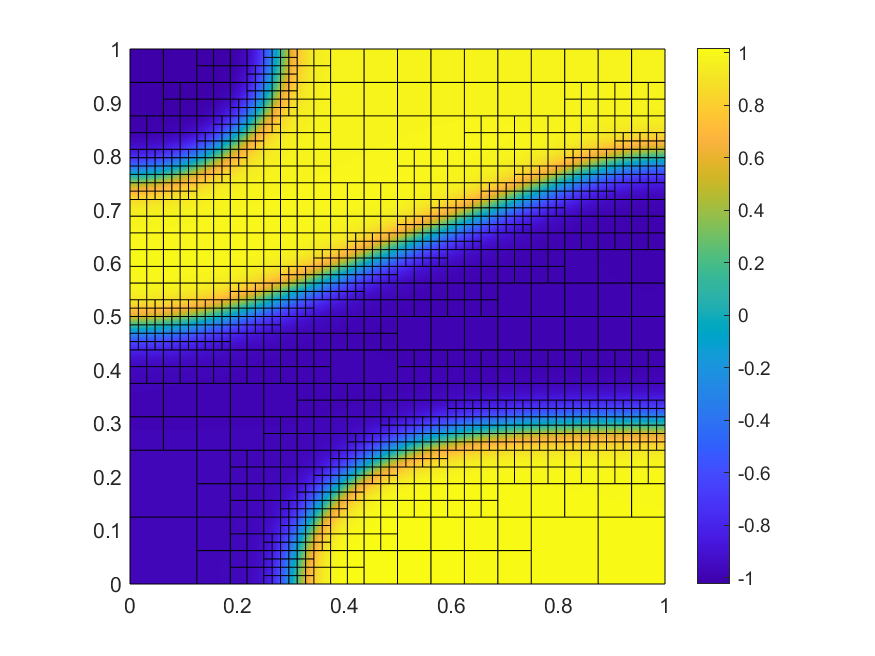} 
\end{minipage}
&
(f)
\begin{minipage}[b]{.45\textwidth}
\centering
\includegraphics[width=\textwidth]{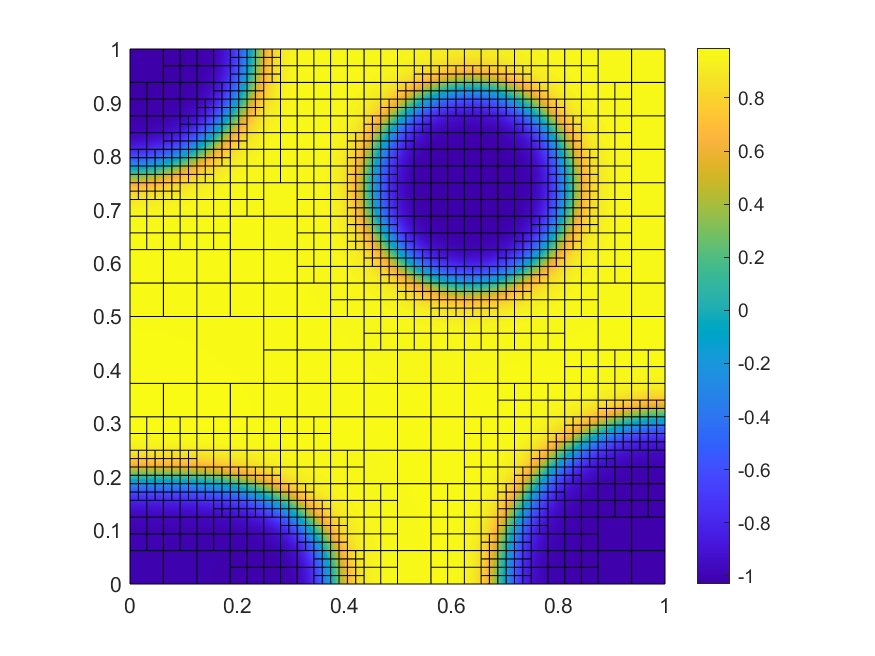}
\end{minipage}
\end{tabular}
\caption{\label{fig:sp_simulations} Contour plots at different instants representing spinodal decomposition ((a), (c), and (e)) and nucleation ((b), (d), and (f) in a single-patch geometry initially discretized using 4,356 degrees of freedom  (dofs) and setting $\alpha_q=0.2$. (a) 2,265 dofs at time $t=0.1$, (b) 2,576 dofs at time $t=0.1$, (c) 1,707 dofs at time $t=0.5$, (d) 1,687 dofs at time $t=0.5$, (e) 1,249 dofs at time $t=1$, and (f) 1,347 dofs at time $t=1$.
}
\end{figure}

The results in Figure \ref{fig:parameter_study}(a) confirm that, in a spinodal decomposition, increasing the value of $\alpha_q$ the number of degrees of freedom employed in the simulation is reduced at the cost of an increased error, as expected, and an equivalent behavior is obtained in a nucleation process, as shown in Figure~\ref{fig:parameter_study}(b). However, even using a coarse mesh, the results are qualitatively correct since the transition area is discretized using sufficiently fine elements ($h=1/64$), as demonstrated by the comparison between the simulation results in Figure \ref{fig:sp_simulations} and the reference solutions in Figs. \ref{fig:uniform_spin_square}(b) and \ref{fig:parameter_study}(c).

\subsection{Phase separations on multi-patch geometries}
\label{sec:num_mpC1}

In these last tests, we simulate phase separations on a geometry defined by three patches following the strategy introduced in Section \ref{sec:ada_mp}. Differently from previous tests on single-patch domains, splines with reduced continuity ($C^{p-2}$) are required to construct the hierarchical space with $C^1$-continuity across the patch interfaces, and therefore we select cubic ($p=3$) $C^1$-continuous splines.

\begin{figure}[t!]
%
%
\centering
\includegraphics[width=.5\textwidth]{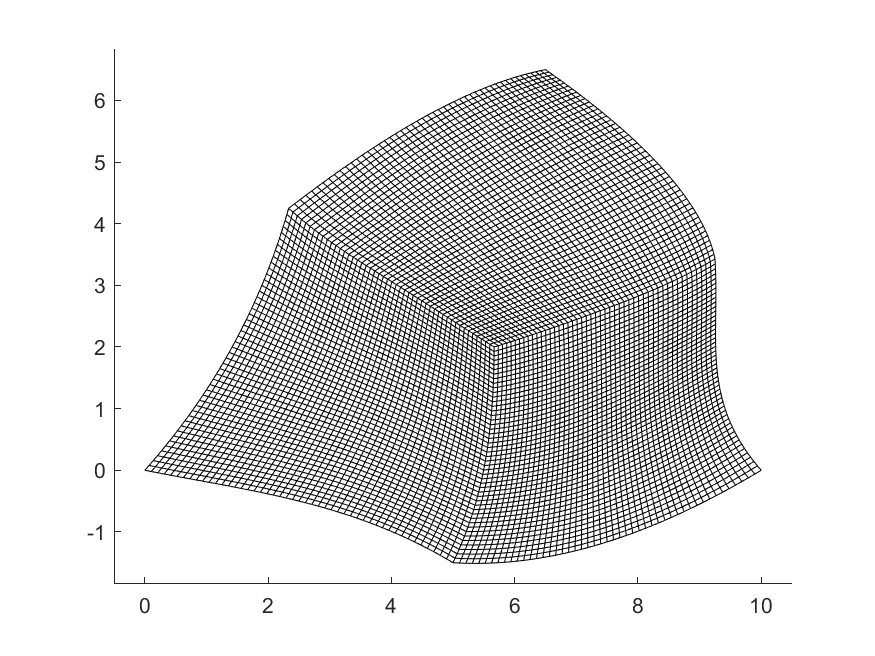}
%
\caption{\label{fig:initial_mesh_mp}Initial mesh for the multi-patch geometry (27,369 dofs) employed in the simulations of spinodal decomposition and nucleation in Figs. \ref{fig:mp_spinodal} and \ref{fig:mp_nucleation}.}
\end{figure}

Spinodal decomposition and nucleation tests are conducted using the same interface parameter $\lambda= 0.05$, the same initial mesh composed by four hierarchical levels ($N=4$) shown in Figure~\ref{fig:initial_mesh_mp}, and the same time discretization ($t \in \left[ 0, 50\right]$, $\Delta t=0.1$) but varying the initial average concentration $\bar{u}$.

\begin{figure}[p]
\begin{tabular}{cc}
(a)
\begin{minipage}[b]{.45\textwidth}
\includegraphics[width=\textwidth]{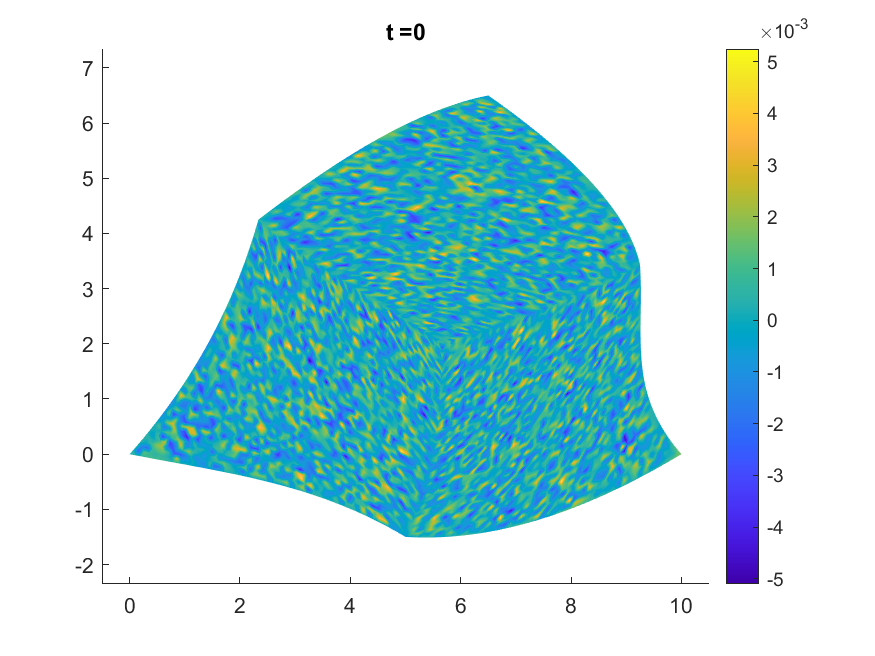} 
\end{minipage}
&
(b)
\begin{minipage}[b]{.45\textwidth}
\includegraphics[width=\textwidth]{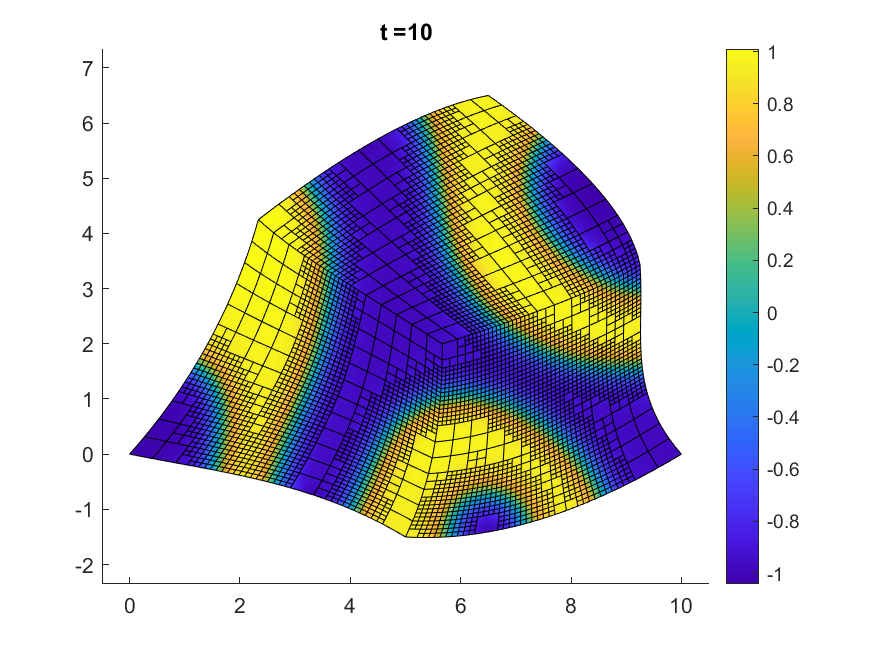} 
\end{minipage}
\\
(c)
\begin{minipage}[b]{.45\textwidth}
\centering
\includegraphics[width=\textwidth]{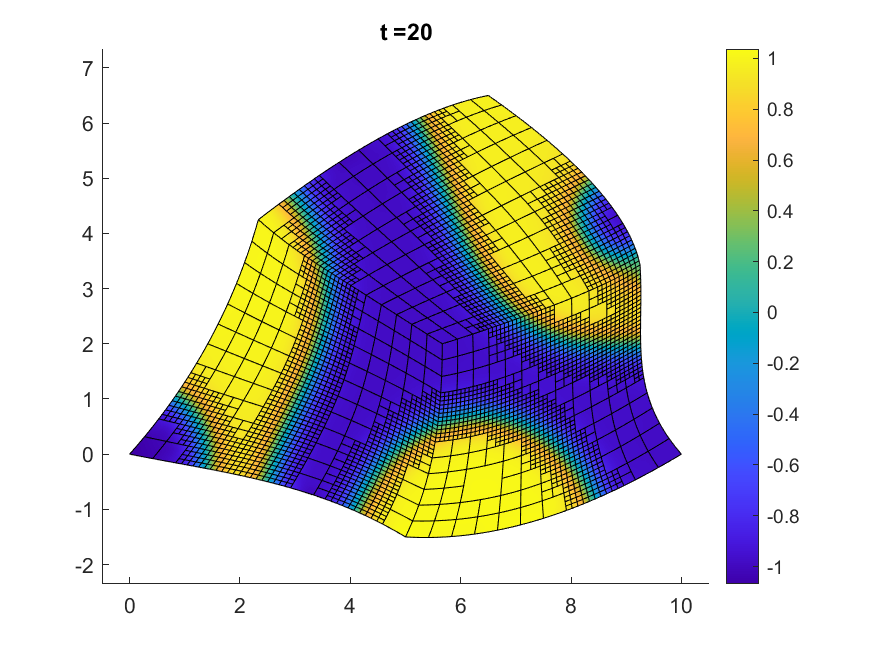} 
\end{minipage}
&
(d)
\begin{minipage}[b]{.45\textwidth}
\centering
\includegraphics[width=\textwidth]{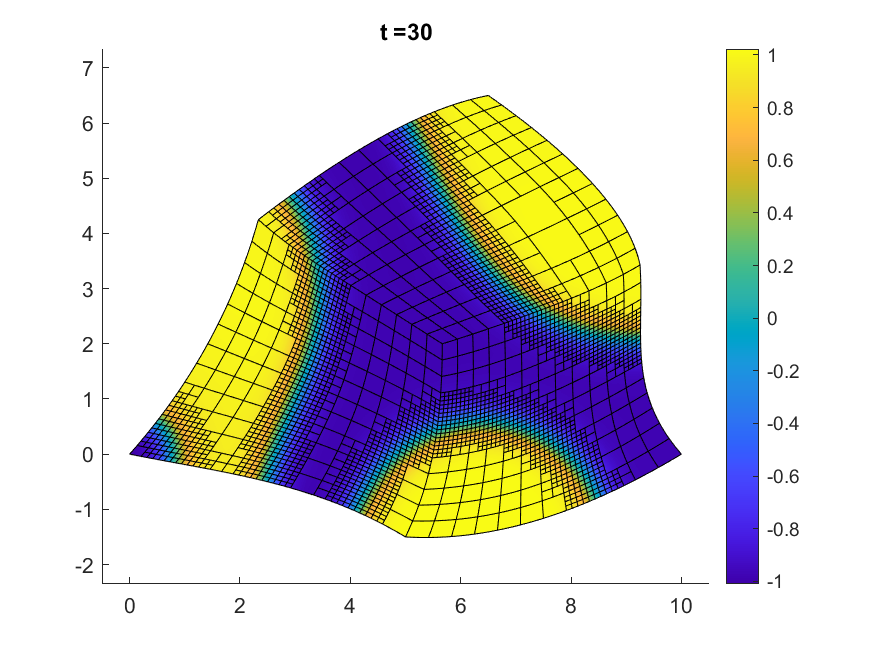} 
\end{minipage}
\\
(e)
\begin{minipage}[b]{.45\textwidth}
\centering
\includegraphics[width=\textwidth]{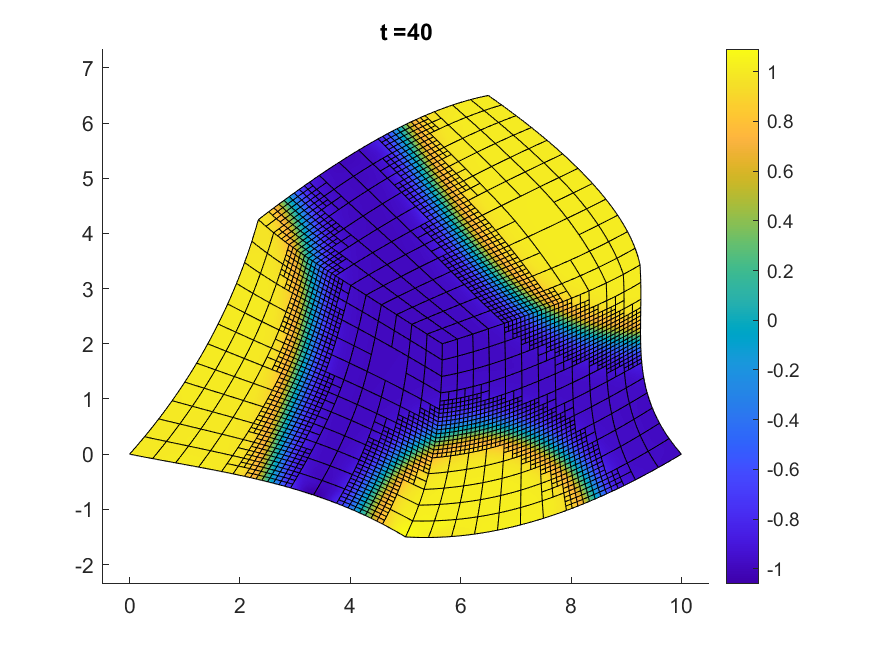} 
\end{minipage}
&
(f)
\begin{minipage}[b]{.45\textwidth}
\centering
\includegraphics[width=\textwidth]{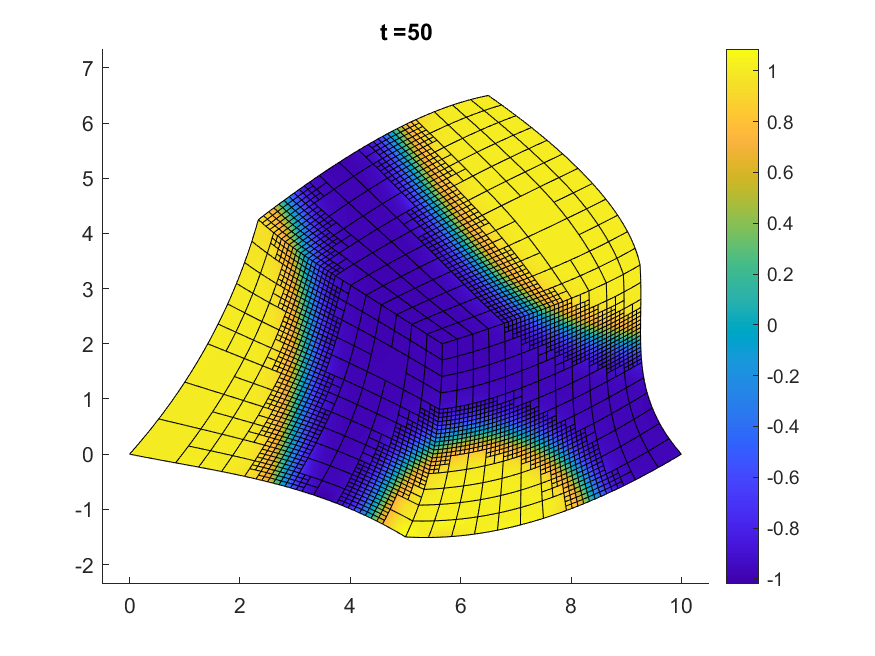}
\end{minipage}
\end{tabular}
\caption{\label{fig:mp_spinodal} Contour plots at different instants representing a spinodal decomposition in a multi-patch geometry. (a) Initial condition (the mesh, comprising 27,369 dofs, is shown in Figure \ref{fig:initial_mesh_mp}). (b)-(f) Evolution of the phase-field and corresponding adaptive mesh: (b)  16,149 dofs at time $t=10$, (c) 12,774 dofs at time $t=20$, (d)  10,673 dofs at time $t=30$, (e) 9,820 dofs at time $t=40$, and (f) 9,674 dofs at time $t=50$.}
\end{figure}

\begin{figure}[p]
\begin{tabular}{cc}
(a)
\begin{minipage}[b]{.45\textwidth}
\includegraphics[width=\textwidth]{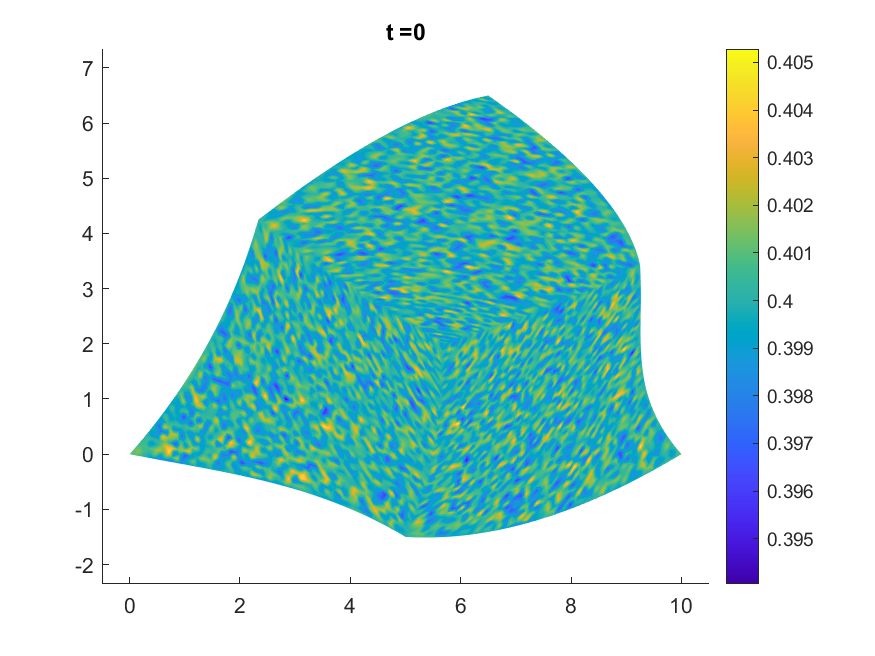} 
\end{minipage}
&
(b)
\begin{minipage}[b]{.45\textwidth}
\includegraphics[width=\textwidth]{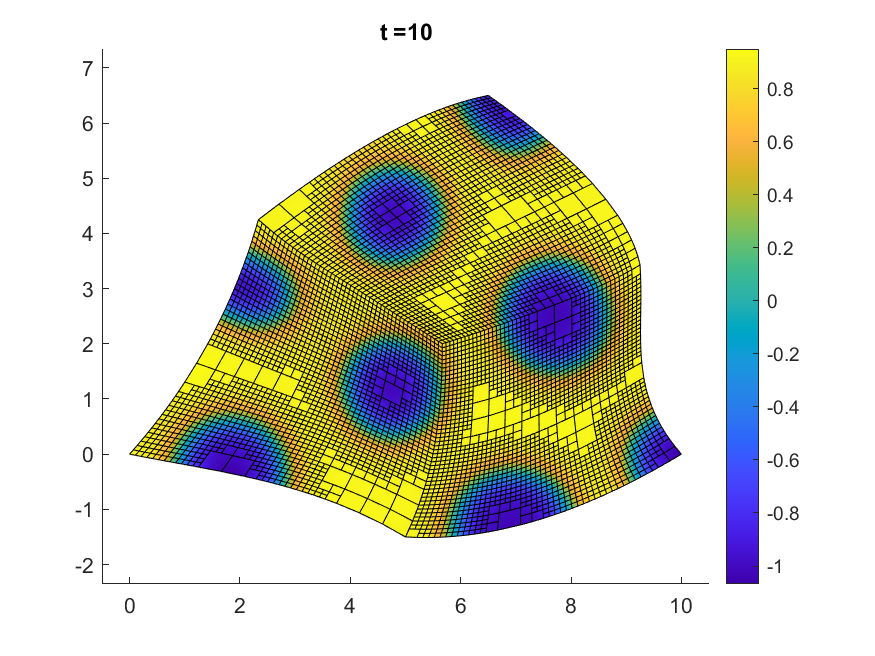} 
\end{minipage}
\\
(c)
\begin{minipage}[b]{.45\textwidth}
\centering
\includegraphics[width=\textwidth]{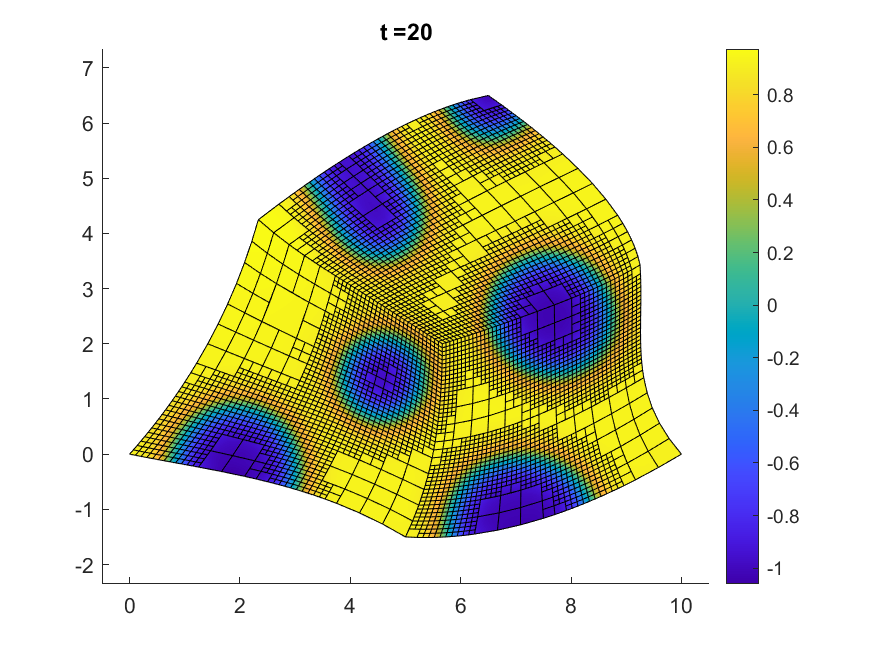} 
\end{minipage}
&
(d)
\begin{minipage}[b]{.45\textwidth}
\centering
\includegraphics[width=\textwidth]{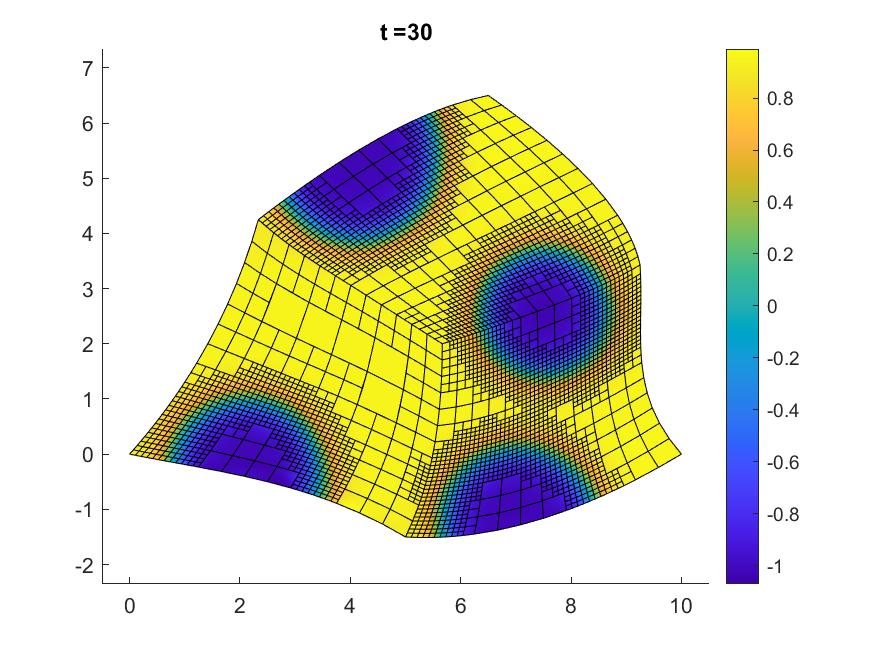} 
\end{minipage}
\\
(e)
\begin{minipage}[b]{.45\textwidth}
\centering
\includegraphics[width=\textwidth]{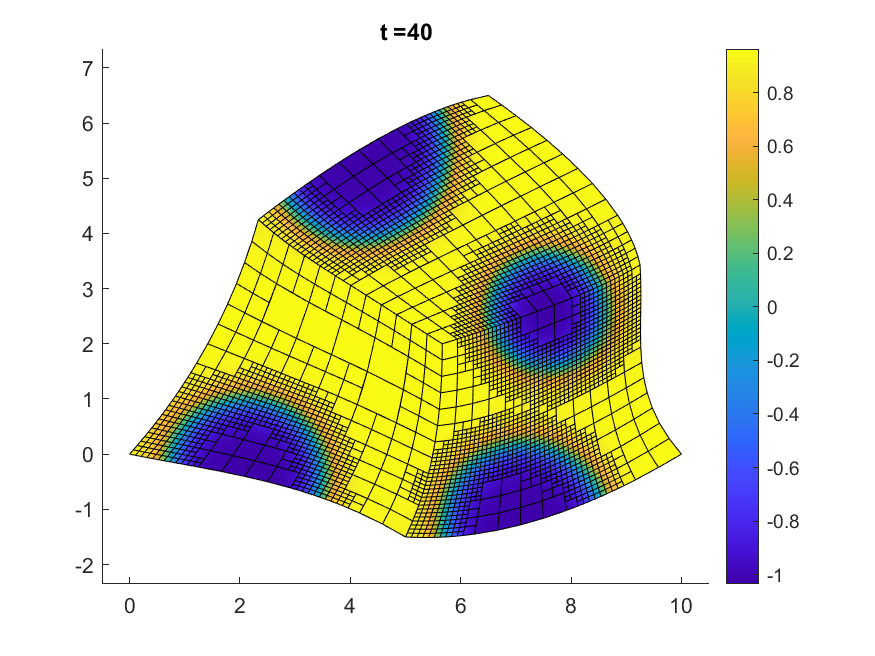} 
\end{minipage}
&
(f)
\begin{minipage}[b]{.45\textwidth}
\centering
\includegraphics[width=\textwidth]{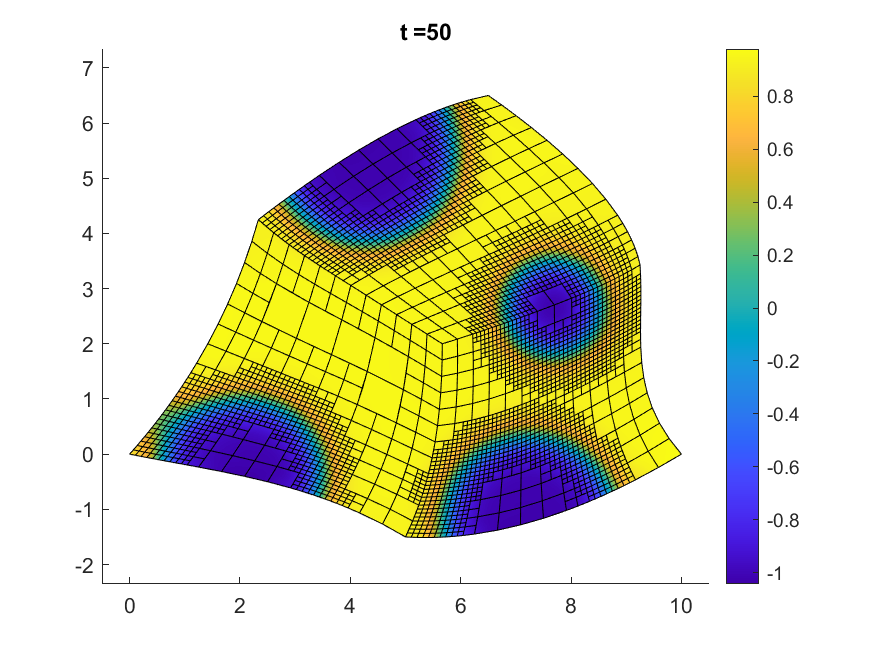}
\end{minipage}
\end{tabular}
\caption{\label{fig:mp_nucleation} Contour plots at different instants representing nucleation in a multi-patch geometry. (a) Initial condition (the mesh, comprising 27,369 dofs, is shown in Figure \ref{fig:initial_mesh_mp}). (b)-(f) Evolution of the phase-field and corresponding adaptive mesh: (b)  22,029 dofs at time $t=10$, (c) 17,987 dofs at time $t=20$, (d)  13,030 dofs at time $t=30$, (e) 12,800 dofs at time $t=40$, and (f) 12,566 dofs at time $t=50$.}
\end{figure}

The results of the spinodal decomposition, reported in Figure \ref{fig:mp_spinodal}, and nucleation process, reported in Figure \ref{fig:mp_nucleation}, obtained setting $\alpha_q=0.1$ show a substantial reduction of the number of degrees of freedom with respect to the initial uniform mesh. Moreover, the simulations confirm the applicability of the adaptive strategy with non-matching discretizations on the patch interfaces, where the $C^1$-continuity of the basis function is achieved through proper constructions of the hierarchical spaces.

\section{Conclusions}
\label{sec:closure}

The focus of the present work is the adaptive isogeometric analysis of phase-field models, herein exemplified by the fourth-order (nonlinear) Cahn-Hilliard equation describing the phase separation of immiscible fluids. In fact, the nature of phase-field modeling inherently calls for local refinement and coarsening, with the phase-field variable naturally providing an indicator to drive adaptivity.
The combination of isogeometric analysis and THB-splines indeed constitutes an ideal context for this, at least in the single-patch case, and, therefore, we have proposed and successfully tested a suitably graded hierarchical adaptive framework, that can be easily generalized to other phase-field modeling problems. In particular, besides showing the advantages that adaptivity can provide in terms of reduced number of degrees freedom, we have also discussed the effects of two different choices of indicator function, as well as the impact of admissibility in the refinement/coarsening process.
Moreover, since the solution in primal form of the fourth-order Cahn-Hilliard equation requires $C^1$-continuity, a suitable multi-patch strategy guaranteeing higher continuity across (non-matching) patch boundaries has to be employed whenever non-trivial geometries are involved. To this end, we have extended our adaptive isogeometric phase-field modeling framework to the multi-patch case adopting the strategy introduced in \cite{BrGiKaVa23}. Numerical tests\footnote{All proposed implementations and numerical examples will be made freely available within GeoPDEs \cite{geopdesv3} before publication of this paper.} have proven the effectiveness of our approach also in this context. 

The extension of the proposed adaptive framework to other multi-patch strategies available in the IGA literature, like, e.g., those relying on almost $C^1$ splines \cite{takacstoshniwal2023} or on weak patch coupling based on penalty or Nitsche's methods (see, e.g., \cite{proserpio2020framework,leonetti2020robust}), will be the object of future work. 
As already stated above, further research will be also carried out in the direction of adaptive isogeometric analysis of other interesting applications of phase-field models, like fracture \cite{proserpio2021phase} or tumor growth \cite{lorenzo2019}.
Clearly, the extension to 3D problems will be considered in the future as well. While this does not appear to be an issue for single-patch geometries, a proper strategy needs to be identified in the multi-patch case, most probably, at the present state, resorting to  weak patch-coupling methods.

\section*{Acknowledgements}
CB, CG, AR, and MT acknowledge 
the contribution of the National Recovery and Resilience Plan, Mission 4 Component 2 – Investment 1.4 – CN\_00000013 ``CENTRO NAZIONALE HPC, BIG DATA E QUANTUM COMPUTING'', spoke 6. RV has been partially supported by the Swiss National Science Foundation via the project n.200021\_188589. CB, CG, and RV are members of the INdAM research group GNCS. The INdAM-GNCS support is gratefully acknowledged.

%

\bibliographystyle{abbrv}

\bibliography{biblio}

\end{document}